\def\TITLE     {Computably Based Locally Compact Spaces}
\def\AUTHOR                {Paul Taylor}

\def\ABSTRACT {ASD (Abstract Stone Duality) is a re-axiomatisation of
  general topology in which the topology on a space is treated, not as
  an infinitary lattice, but as an exponential object of the same
  category as the original space, with an associated lambda-calculus.
  In this paper, this is shown to be equivalent to a notion of
  \emph{computable basis} for locally compact sober spaces or locales,
  involving a family of open subspaces and accompanying family of
  compact ones.  This generalises Smyth's \emph{effectively given
  domains} and Jung's \emph{strong proximity lattices}.  Part of the
  data for a basis is the inclusion relation of compact subspaces
  within open ones, which is formulated in locale theory as the
  \emph{way-below} relation on a continuous lattice.  The finitary
  properties of this relation are characterised here, including the
  \emph{Wilker condition} for the cover of a compact space by two open
  ones.  The real line is used as a running example, being closely
  related to Scott's \emph{domain of intervals}.  ASD does not use the
  category of sets, but the full subcategory of overt discrete objects
  plays this role; it is an \emph{arithmetic universe}
  (\emph{pretopos} with lists).  In particular, we use this
  subcategory to translate computable bases for classical spaces into
  objects in the ASD calculus.}

\newif\iflmcs\lmcstrue
\ifx\puredvitrue\undefined
   \expandafter\newif\csname ifpuredvi\endcsname\puredvitrue\fi
\ifx\ptsz\undefined\let\ptsz\empty\fi % default to 10pt

\iflmcs\documentclass{LMCS-Taylor}
\else  \expandafter\documentclass\ptsz{article}\usepackage{a4wide}
       \author{\AUTHOR}\title{\TITLE}
\fi

\usepackage{abbrevs}
\usepackage{amssymb}
\usepackage{eqnarray}
\usepackage[linktocpage]{hyperref}
\usepackage{letterlist}
\usepackage{logicsym}
\usepackage{prooftree}
  \def\wwwname{www.cs.man.ac.uk}
\usepackage{pthref}

\usepackage[PostScript=Rokicki,height=2.8em,width=3.5em,%
           tight,repositionpullbacks,nohug]{diagrams}

\usepackage{asdsymb}

\ifx\df\undefined
   \ifx\selectfont\undefined
      \let\df\bf
   \else
     \def\df{\fontfamily\rmdefault\fontseries\bfdefault\fontshape\itdefault
      \selectfont}
   \fi
\fi
\def\textdf#1{{\df #1}}
\ifx\mathbb\undefined\def\mathbb{\Bbb}\fi
\ifx\cal\undefined\def\cal{\mathcal}\fi
\def\AAA{{\mathbb A}}
\def\sizenine#1{{\hbox{\small$#1$}}}
\def\USet#1{{\left|{#1}\right|}}
\def\half{{\textstyle{\frac 1 2}}}
\def\sizenine#1{{\hbox{\small$#1$}}}

\def\trless{\mathrel{\triangleleft}}
\def\tridown    {\mathbin{\vtop{% relative slice category
                \kern-9pt\hbox{$\vert$}\nointerlineskip\kern-1pt
                \hbox{$\mkern-3mu\scriptscriptstyle\bigtriangledown$}}%
                \mkern-2mu}}%

\iflmcs
   \theoremstyle{definition}
   \theoremstyle{definition}
   \theoremstyle{definition}
 
   \def\newslthm#1#2{\newenvironment{r@#1}{\begin{#2}}{\end{#2}}}
   \def\newrmthm#1#2{\newenvironment{r@#1}{\begin{#2}}{\end{#2}}}

   \def\Proof{\proof\quad}

   \def\goalbreak#1{}
   \def\allowlines#1{}
   \def\closeupaline
       {\ifvmode\nobreak\vskip-\baselineskip
        \else\vadjust{\nobreak\vskip-\baselineskip}\fi}
   \def\WillHandleQED{}
   \def\HandleQED{}
   \def\qeds{\relax\hbox{$\square$}}
\else
   \usepackage[define-standard-theorems,theorems-as-commands]{QED}
   \def\newslthm#1#2{\newenvironment{r@#1}{\csname #1\endcsname}{}}
   \def\newrmthm#1#2{\newenvironment{r@#1}{\csname #1\endcsname}{}}
   \def\qeds{\qed}
\fi

\newslthm{Corollary}{cor}
\newslthm{Lemma}{lem}
\newslthm{Proposition}{prop}
\newslthm{Theorem}{thm}

\newrmthm{Axiom}{axm}
\newrmthm{Definition}{defi}
\newrmthm{Examples}{exas}
\newrmthm{Example}{exa}
\newrmthm{Notation}{notn}
\newrmthm{Remark}{rem}
\newrmthm{Warning}{wng}

\iflmcs
\def\doi{2 (1:1) 2006}
\lmcsheading%
{\doi}
{70}
{}
{}
{Sep.~\phantom{0}3, 2004}
{Mar.~\phantom{0}3, 2006}
{}
\fi

\begin{document}
\iflmcs
  \author[Paul Taylor]{\AUTHOR}
%  \title[\TITLE ]{\TITLE }
  \title{Computably Based Locally Compact Spaces}
  \address{University of Manchester}
  \email{pt@cs.man.ac.uk}
  \keywords{abstract Stone duality, %
    locally compact, %
    basis for topology, %
    effectively presented, %
    way below relation, %
    continuous lattice, %
    Sigma-split subspace}
  \subjclass{F.4.1}
  \amsclass{54D45, 03D45 (Primary), 06B35, 54D30, 68N18 (Secondary)}
  \begin{abstract}
    \noindent\ABSTRACT
  \end{abstract}
  \maketitle
\else
  \maketitle\begin{abstract}\noindent\ABSTRACT\end{abstract}
\fi

\abovedisplayskip=6pt plus 3pt minus 3pt
\abovedisplayshortskip=0pt plus 3pt
\belowdisplayskip=6pt plus 3pt minus 3pt
\belowdisplayshortskip=4pt plus 3pt minus 1pt

%\begin{center}
%  \small
%  \def\arraystretch{1}
%  \begin{tabular}{rlrcrlr}
%   \rlap{\textbf{Contents}}\phantom{00} &&&\quad&
%      %
%    9.& Basic corollaries                       &\pageref{corolls}\\
%    1.& Locally compact sober spaces            &\pageref{intro spaces}&&
%   10.& Primes and nuclei                       &\pageref{prime/nucleus}\\
%    2.& Locally compact locales                 &\pageref{intro locales}&&
%   11.& The way-below relation                  &\pageref{way-below}\\
%    3.& Axioms for abstract Stone duality       &\pageref{axioms I}&&
%   12.& Domain theory in ASD                    &\pageref{domain theory}\\
%    4.& Sets, unions and bases                  &\pageref{axioms II}&&
%   13.& The lattice basis on $\Sigma^N$         &\pageref{DLmonad}\\
%    5.& Compact subobjects                      &\pageref{cpct subsp}&&
%   14.& From the basis to the space             &\pageref{X from <<}\\
%    6.& Effective bases                         &\pageref{lcpct}&&
%   15.& The points of the new space             &\pageref{new points}\\
%    7.& $\Sigma$-split subobjects               &\pageref{subobjects}&&
%   16.& Morphisms as matrices                   &\pageref{matrix}\\
%    8.& Every definable object has a basis      &\pageref{bases}&\quad&
%   17.& Relating the two models                 &\pageref{relate}\\
%      %
%      & \kern .35\textwidth &&&& \kern .35\textwidth
%  \end{tabular}
%\end{center}

% ============================================================================
\vskip-\bigskipamount
%{\small
\tableofcontents
%}

\section*{Introduction}%\label{intro}
%\section{Introduction}%\label{intro}

A \emph{locally compact space} is one in which there is
a good interaction of \emph{open} and \emph{compact} subspaces.
In this paper we shall show how this interaction can be captured abstractly,
establishing an equivalence with a new recursive account of topology.

%\begin{r@Remark}
%\label{summary loc cpct}
We shall consider the following formulations of local compactness:
\begin{letterlist}
\item open and compact subspaces in a traditional sober topological space
  (Section~\ref{intro spaces});
\item disjoint open subspaces in a compact Hausdorff (regular) space
  (Example~\ref{cpct Hdf base}); and
\item the way-below relation $\ll$ in a continuous distributive lattice
  (Section~\ref{intro locales}).
\end{letterlist}

\smallskip

In the categories of locally compact sober spaces and locales,
the topology (lattice of open subspaces) of $X$, considered as another
space or locale equipped with the Scott topology,
is the exponential $\Sigma^X$, where $\Sigma$ is the Sierpi\'nski space.
We shall show how the three classical formulations above give rise to
the following abstract structure in this category:
\begin{letterlist*}
\item families of terms $\beta^n:\Sigma^X$ and $A_n:\Sigma^{\Sigma^X}$
  with a \textdf{basis expansion} (Section~\ref{lcpct})
  $$ \phi x \;\eq\; \Some n.A_n\phi\land \beta^n x\hbox{;} $$
\item a $\Sigma$-split inclusion $i:X\splitinto\Sigma^N$ and
  $I:\Sigma^X\retract\Sigma^{\Sigma^N}$ with $\Sigma^i\cdot I=\id_{\Sigma^X}$
  (Section~\ref{subobjects});  
\item an idempotent $\E\equiv I\cdot\Sigma^i$ on $\Sigma^{\Sigma^N}$
  (Section~\ref{prime/nucleus});
\item an abstract way-below relation $n\waybelow m\equiv A_n\beta^m$ on $N$
  (Section~\ref{way-below}); and
\item an Eilenberg--Moore algebra for the monad arising from the adjunction
  $\Sigma^\blank\adjoint\Sigma^\blank$ \cite{TaylorP:subasd}.
\end{letterlist*}

\smallskip

We will give complete axiomatisations of each of these structures,
and define the translations amongst them,
doing this in a new $\lambda$-calculus (Abstract Stone Duality).
We consider several special cases, in particular the real line,
the interval domain (Definition~\ref{intdom}),
continuous dcpos and various kinds of domains (Section~\ref{domain theory}).

%\end{r@Remark}

\medskip

It may help you to find your way around this lengthy paper
if I explain that Section~\ref{lcpct} was the first to be written.
It shows how the classical notions of locally compact space and locale
can be expressed in a $\lambda$-calculus.

However, the present paper was the first to put this calculus to work
to do topology.
The previous ones \cite{TaylorP:sobsc,TaylorP:subasd,TaylorP:geohol}
had examined \emph{components} of the theory as abstract categorical ideas, 
without being aware of how many (in fact, few) more such components
would be needed to build a logically complete account of topology.
Sections \ref{axioms I}--\ref{axioms II} and
\ref{subobjects}--\ref{prime/nucleus}
discuss the way in which these components go together to make a viable theory.
This discussion has, in fact, spilled over into several additional papers,
in particular \cite{TaylorP:insema,TaylorP:eletvc}.
If you would prefer to read a dry summary of the axioms in a
``user manual'' style, please see \cite{TaylorP:dedras}.

Once all this structure is in place,
Section~\ref{way-below} begins to take it apart again,
characterising the finitary properties of the way below
relation in an abstract way.
This corresponds to the concrete treatment in Sections \ref{intro spaces}
and~\ref{intro locales}.
Sections~\ref{DLmonad}--\ref{matrix} construct an object of the ASD calculus
from any such abstract relation
(with a ``practice run'' for domains in Section~\ref{domain theory})
and the final section completes the circle with the traditional theory.

% ============================================================================
\section{Locally compact sober spaces}\label{intro spaces}

\begin{r@Definition}
\label{pt loc cpct}
The traditional definition of local compactness was generalised
from Hausdorff to sober spaces in \cite[p.~211]{HofmannKH:locccl}:
 Whenever a point is contained in an
open subspace ($x\in V$), there is a compact subspace $K$ and an open
one $U$ such that $x\in U\subset K\subset V$.
\end{r@Definition}

\smallskip

It is an easy exercise in the ``finite open sub-cover'' definition of
compactness to replace the point $x$ by another compact subspace:

\begin{r@Lemma}
\label{classical interpolation}
Let $L\subset V\subset X$ be compact and open subspaces of a
locally compact space. Then there are
$ L \subset U \subset K \subset V \subset X$
with $U$ open and $K$ compact. \qed
\end{r@Lemma}

\smallskip

We call this result the \textdf{interpolation property}.
Alternating inclusions of open and compact subspaces like this
will be very common in this paper,
and compact subspaces will be more important than points.

\begin{r@Notation}
\label{spatial <<} We write $U\waybelow V$ and $K\waybelow L$
if there is such an interpolating compact or open subspace, respectively. 
The second version, in which $K\supset L$, follows the usage of
\cite{HofmannKH:locccl,JungA:duacvo}, \cf our Theorem~\ref{cpct <<}.
\end{r@Notation}

\medskip

Now consider what we might mean by a \emph{computably defined}
locally compact space.

\smallskip

Suppose that you have some computational representation of a space.
It can only encode \emph{some} of the points and open and compact
subspaces, since in classical topology there are uncountably many of them
in any interesting case.
Hence your ``space'' cannot be literally sober, or have arbitrary
unions of open subspaces.
We understand the intended space to be the corresponding sober one,
in which arbitrary unions of opens have also been adjoined.
Those points, open and compact subspaces that have codes
are called \textdf{basic}.

\begin{r@Example}
\label{eg R1} A computable definition of $\realno$ as a locally compact
space might have
\begin{letterlist}
\item as \emph{basic points}, (encodings of) the rationals;
\item as \emph{basic open subspaces}, the (names of) open intervals
  $(q\pm\epsilon)\equiv(q-\epsilon,q+\epsilon)\equiv 
  \collect{x}{\left|x-q\right|\lt\epsilon}$
  with rational or infinite endpoints; and
\item as \emph{basic compact subspaces}, the closed intervals
  $[q\pm\delta]\equiv[q-\delta,q+\delta]\equiv
  \collect{x}{\left|x-q\right|\leq\delta}$
  with finite rational endpoints.
\end{letterlist}
Notice that \emph{both} open and closed intervals are used,
although treatments of exact real arithmetic
often just use one or the other, \cf Example~\ref{R-SN}.
Also, by ``rationals'' we might actually mean
all pairs $p/q$ with integers $p$ and $q\neq 0$,
or the dyadic rationals $p/2^n$, or continued fractions,
or whatever our favourite countable dense set of reals may be.
Unlike Dedekind cuts, this example readily generalises: for $\realno^3$
we use open and closed cuboids whose vertices have rational co-ordinates
(or, better, a system based on close packing of spheres
\cite{ConwayJH:sphplg}).
\end{r@Example}

In the Example, the intersection of two basic opens is again basic,
but for technical reasons we shall also need to extend the families
to include \emph{finite unions} of open (respectively, compact) intervals
or cuboids.
It's an exercise that's a little too complicated to be called algebra,
but easy programming,
to test inclusion and compute the representations
of such unions and intersections.

\begin{r@Definition}
\label{rec loc cpct}
A \textdf{computably based locally compact space} consists of
a set of codes for basic ``points'', ``open'' and ``compact'' subspaces,
together with an interpretation of these codes
in a locally compact sober space.
We require of the space that every open subspace be a union of basic ones.
We also want to be able to \emph{compute}
\begin{letterlist}
\item codes (that we shall just call $0$ and $1$) for the \emph{empty} set
  and the \emph{entire} space, considered as open and compact subspaces
  (if, that is, the entire space is in fact compact);
\item codes for the \emph{union} and \emph{intersection}
  of two open subspaces, and
  for the union of two compact ones, given their codes
  (we write $+$ and $\star$ instead of $\cup$ and $\cap$
  for these binary operations, to emphasise that they act on codes,
  rather than on the subspaces that the codes name);
\item whether a particular representable \emph{point} \emph{belongs}
  to a particular basic \emph{open} subspace, given their codes;
  but we only need a positive answer to this question if there is one,
  as failure of the property is indicated by non-termination;
\item more generally,
  whether an \emph{open} subspace \emph{includes} a \emph{compact} one,
  given their codes;
\item codes for $U$ and $K$ such that $L\subset U\subset K\subset V$,
  given codes for $L\subset V$ as above.
\item\label{UKpair} In fact,
  we shall require the basic compact and open subspaces to come in pairs,
  with $U^n\subset K^n$ as in \cite{JungA:duacvo},
  where the superscript $n$ names the pair,
  and we also need part (e) to yield such a pair as the interpolant.
\end{letterlist}
Extensional equivalence of computable functions is not captured
within the strength of the logic that we wish to study.
So, for the ``computations'' above,
we mean a particular program to be specified ---
at least up to provable equivalence,
which means that we don't have to nominate a programming language.
\end{r@Definition}

\begin{r@Definition}
\label{rec cts} A \textdf{computably continuous function}
between such spaces is a continuous function $f:X\to Y$ between the
topological spaces themselves, for which the binary relation
$$ f K^m \;\subset\; U^n $$
(between the \emph{codes} $(n,m)$ for a compact subspace $K^m\subset X$
and an open one $U^n\subset Y$) is~recursively enumerable, \cf part (c)
of the previous definition.

In particular, \textdf{computably equivalent bases} for the same
space are those for which the identities in both directions are
computably continuous functions.
This means that the relations $K^n\subset U^m$ and $K^m\subset U^n$
between $n$ and $m$ are recursively enumerable.
For example, whilst there are several choices for the ``rationals''
and intervals in Example~\ref{eg R1}, all of the reasonable ones are
computably equivalent.
\end{r@Definition}

\medskip

\begin{r@Remark}
\label{spc as ref} It is out of place in the definition of
something \emph{computable} to specify a topological space:
this was only included to guarantee topological consistency
of the computations of union, intersection, containment and interpolation.
This is also the reason why the codes for basic points
played no actual role in Definition~\ref{rec loc cpct} ---
there needn't be any of them.
Since we said that the spaces are sober,
continuous functions are determined entirely by their effect
on basic open and compact subspaces, as they are in locale theory.
Neither the (optional) basic points, nor the computable ones that
we recover in Section~\ref{new points},
nor the uncountably many points of the classical theory,
are actually needed to specify continuous functions,
or prove topological properties such as compactness.

At the technical level, therefore,
what we call a computably based locally compact ``space''
is really just a system of codes and programs acting on them,
for which \emph{there exists} a topological space 
in which these codes may be interpreted as compact and open subspaces.
Likewise, a computably ``continuous function'' is
just an equivalence class of programs,
for which \emph{there exists} a continuous function
that agrees with the program in the way that we have said.

The particular way in which computations are set up is of course a
matter for discussion: I still have to convince you of the merits of
my way, in comparison to the many other ways that have been used to
calculate in $\realno^\nn$ over the millennia,
and of course any theoretical picture is subject to technical optimisation.
However, these other ways share the same pattern, in that they have a
\emph{technical} definition that is (a)~used for applied computation but
(b)~justified by the \emph{existence} of an interpretation in some
``pure'' theory of (Euclidean) space.
The pure theory is not the computation,
at any rate if it is based on twentieth century topology.
So the practical situation is illustrated by Example~\ref{eg R1}:
we have a classical definition of a topological space,
equipped with a basis that is defined in some conventional way,
and which we want to use to obtain values in the space.
So long as the above features of the basis are computable,
the classical space guarantees the consistency conditions.

In Section~\ref{way-below} we shall formulate this consistency in terms of
finitary conditions on the transformations of encodings themselves,
eliminating the topological space,
\ie we prove necessity of the abstract conditions.
Then in Section~\ref{X from <<} we prove sufficiency,
\ie that \emph{any} such encoding (such as the one for $\realno$ above)
satisfying these consistency requirements does define a locally
compact sober space, and this is unique up to homeomorphism.

This construction is done, not in traditional topology itself,
which is not computable, but in a $\lambda$-calculus.
Therefore, \emph{the construction imports the classical data
into the computational world} (Section~\ref{relate}).
But again there is a distinction between the technical
and philosophical meaning:
whilst our theory (ASD) starts from (and is reducible back to)
\emph{computational} foundations,
it adapts these to give a new account of the \emph{mathematical} theory
that is ultimately intended to be usable by pure mathematicians
in place of the old one.
\end{r@Remark}

\begin{r@Remark}
\label{basis dist latt} The operations $\star$ and $+$ become
$\cap$ and $\cup$ when we interpret them \emph{via} the basic open
subspaces that they encode.
Amongst the abstract consistency requirements, therefore,
we would expect $(0,1,{+},{\star})$ to define a \emph{distributive lattice}.

However, we have only asked for the ability to test inclusion of
a compact subspace in an open one,
not inclusion or equality of two open or two compact subspaces,
nor of an open subspace in a compact one
(the inclusion $U^n\subset K^n$ is given, not tested).
Even if these happen to be possible,
it is computationally quite reasonable for different codes
to denote the same subspace,
but for this fact to be potentially undecidable.

On the other hand, as we want to stress the computable aspect
of the names of basic subspaces,
we shall often represent $+$ and $\star$ as concatenations of lists.
This clumsiness actually serves an expository purpose,
keeping this ``imposed'' structure on codes separate in our minds
from the ``intrinsic'' structure in the topology on~$X$.
(We shall regard this topology as another space.)
If we used the notation and equations of a distributive lattice
for the set $N$ of codes,
it would be all too easy to lapse into confusing it with
the actual topology on the space.
This would in fact make logical assumptions
that amount to a solution of the Halting Problem (or worse).
\end{r@Remark}

\smallskip

The topological information is actually contained,
not in the quasi-lattice structure $(0,1,{+},{\star})$ on~$N$,
but in the inclusion relation between compact and open subspaces.
This satisfies some easily verified properties:

\begin{r@Lemma}
\label{classical 0+<<}\label{cpct<<*}\quad $\emptyset\subset V$,
$$
\begin{prooftree}
  K\subset L\subset U\subset V
  \justifies
  K\subset V
\end{prooftree}
\qquad\qquad
\begin{prooftree}
  K\subset V \qquad L\subset V
  \justifies
  K\cup L \subset V
\end{prooftree}
\qquad\qquad
\begin{prooftree}
  K\subset U \qquad K\subset V
  \justifies 
  K\subset U\cap V
\end{prooftree}
\eqno\qEd
$$
\end{r@Lemma}

\medskip

Finally, there is a property similar to Lemma~\ref{classical interpolation}
that concerns binary \emph{unions}.
Easy enough though this property may be to prove --- when you see it ---
it is not something whose significance one would identify in advance.
Various forms of it were originally
studied by Peter Wilker \cite{WilkerP:adjphf}.

\begin{r@Lemma}
\label{classical Wilker}
Let $K$ be a compact subspace covered by two open subspaces
of a locally compact sober space $X$, that is, $K\subset U\cup V$.
Then there are compact subspaces $L$ and $M$ and open ones $U'$ and $V'$
such that
$$K\subset U'\cup V' \qquad U'\subset L\subset U
  \quad\hbox{and}\quad V'\subset M\subset V.$$
\end{r@Lemma}

\Proof Classically, $K\setminus V$ is a closed subspace of a compact space,
and is therefore compact too, whilst $K\setminus V\subset U$,
so by the interpolation property (Lemma~\ref{classical interpolation}) we have
$$ K\setminus V \subset U' \subset L \subset U \subset X $$
for some $U'$ open and $L$ compact.
Then $K\setminus U' \subset K\setminus (K\setminus V) \subset V$ so
$$ K\setminus U' \subset V' \subset M \subset V \subset X $$
for some $V'$ open and $M$ compact.
Finally,
$K = (K\cap U') \cup (K\setminus U') \subset U'\cup V'$. \qed

\medskip

We didn't mention the \emph{intersection} of two compact subspaces in
Definition~\ref{rec loc cpct},
because there are spaces in which this need not be compact.

\begin{r@Definition}
\label{stab loc cpct sp} A locally compact sober space
is called \textdf{stably locally compact} if the whole space is compact
and the intersection of any two compact subspaces is again compact.
\end{r@Definition}

\begin{r@Examples}
\label{not stably loc cpct}
\begin{letterlist}
\item Consider two copies of the real unit interval $[0,1]$ identified
  on their interiors (or, if you prefer, an interval with duplicated
  endpoints).  Then the two copies of the interval are compact
  subspaces, but their intersection is not
  \cite[Problem~5~B(a)]{KelleyJL:gent}.
  %\cite{JohnstonePT:stos,JohnstonePT:viells}.

\item\label{PCA} A \textdf{combinatory algebra}
  has constants $k$ and $s$ and a (non-associative) binary operation $\cdot$
  such that $(k\cdot x)\cdot y=x$ and
  $\big(((s\cdot x)\cdot y)\cdot z\big) = (x\cdot z)\cdot(y\cdot z)$.
  In the free such algebra $A$, terms can be enumerated,
  and proved equal, using these rules.
  But, since this structure can be used to encode computation
  \cite{BarendregtHP:lamcss},
  proving \emph{in}equality in $A$ is like solving the Halting problem.
  (Similar \emph{unsolvable word problems} can also be set up in
  other algebraic theories, such as groups.)
  \begin{center}
    Hence $A$ is discrete but not Hausdorff.
  \end{center}
  Since it is discrete, its compact subspaces are
  the \textdf{Kuratowski-finite} ones,
  \ie those that can be listed, possibly with repetition.
  In particular, singletons $\setof x$ are open and compact, but not closed.
  However, the intersection $\setof x\cap\setof y$ is compact
  (Kuratowski-finite, listable) iff it is either empty ($x\neq y$)
  or a singleton ($x=y$), \ie iff equality is decidable,
  which it isn't.
\qed
\end{letterlist}
\end{r@Examples}

% ============================================================================
\section{Locally compact locales}\label{intro locales}

Points disappeared from the discussion right at the beginning,
and we saw in Example~\ref{eg R1} that
it is easier to specify $\realno$ with the Euclidean topology
using open and compact subspaces than using open subspaces and points.
Arguably, topology should be axiomatised in this way,
just as traditional geometry was axiomatised in terms of lines and circles
that were entities in themselves, rather than being sets of points.

\textdf{Locale theory} reduces the description further,
to one involving \emph{open} subspaces alone.
To do this for locally compact spaces,
we must represent \emph{compact} subspaces in terms of
their systems of neighbourhoods.
We also have to characterise the situation
$(U\waybelow V)\equiv\Some K.(U\subset K\subset V)$.
We shall see that the way in which this is done for locales
is subtly different from that for spaces,
especially in the non-stably locally compact case
(Remark~\ref{filter/lattice choice}).

\goalbreak{5\baselineskip}

\begin{r@Definition}
 Let $L$ be a complete lattice.
\begin{letterlist}
\item\label{def dirsup} A family $(\psi_s)\subset L$ is called
  \textdf{directed}%
 \footnote{The letter $s$ stands for semilattice, but see Definition~\ref{DJ}.}
  if it is inhabited, and whenever $\psi_r$ and
  $\psi_s$ belong to the family, there is some $\psi_t\geq\psi_r,\psi_s$.
  The join of the family is written $\dirsup\psi_s$.
\item\label{def lattice ll} Now, for $\beta,\phi\in L$,
  we write $\beta\ll\phi$ (\textdf{way-below}) if,
  whenever $\phi\leq\dirsup_s\psi_s$,
  there is already some $s$ for which $\beta\leq\psi_s$.
  (So $\beta\ll\phi$ implies $\beta\leq\phi$.)
\item\label{def cts latt} Then $L$ is a \textdf{continuous lattice}
  \cite{GierzGK:comcl,HofmannKH:locccl} if,
  for all $\phi\in L$,\quad $\phi = \dirsup \collect \beta {\beta\ll\phi}$.
\end{letterlist}
\end{r@Definition}

\begin{r@Proposition}
\label{open <<} The topology of any locally compact space is a
distributive continuous lattice, in which $U\waybelow V$ iff $U\ll V$
\cite[p.~212]{HofmannKH:locccl}.
\end{r@Proposition}

\Proof $U\waybelow V$ implies $U\ll V$ by compactness of $K$ with
$U\subset K\subset V$, and
$$ V = \bigcup \collect W {W\waybelow V} $$
by Definition~\ref{pt loc cpct}.
This union is directed by Lemma~\ref{classical 0+<<},
so it may be used in the definition of $U\ll V$,
giving $U\subset W\waybelow V$ for some $W$, but then $U\waybelow V$ too.
Hence $U\ll V$ iff $U\waybelow V$, but notice that the proof does not
supply an interpolating compact subspace $U\subset K\subset V$.\qed

\smallskip

Conversely, every distributive continuous lattice is the lattice of
open subspaces of some locally compact sober space.
However, this result relies on the axiom of choice,
and even then it is not a trivial matter to prove it 
(Remark~\ref{lacuna}, \cite{HofmannKH:spetcl}).

\medskip

Definition~\ref{rec loc cpct} for spaces has a simpler analogue for
locales, since it's all lattice theory.

\begin{r@Definition}
\label{rec cts latt}
A \textdf{computable basis} $(N,0,1,{+},{\star},{\waybelow})$
for a continuous distributive lattice $L$
is a set $N$ with constants $0,1\in N$,
computable binary operations ${+},{\star}:N\times N\to N$,
a~recursively enumerable binary relation~$\waybelow$
and an interpretation $\beta^\blank:N\to L$
that takes $(0,1,{+},{\star})$ to the lattice structure in $L$, such that
$n\waybelow m$ iff $\beta^n\ll\beta^m$ and
$$ \hbox{for each }\phi\in L, \quad
   \phi \;=\; \dirsup\collect{\beta^n}{\beta^n\ll\phi}.$$
If $L_1$ and $L_2$ have bases $(\beta^m)$ and $(\gamma^n)$
then $H:L_2\to L_1$ is a \textdf{computable frame homomorphism}
if $H$ preserves $\top$, $\land$ and $\bigvee$,
and the relation $(\beta^m\ll H\gamma^n)$ is recursively enumerable.
The interested reader may like to translate the abstract conditions
on this relation that are set out in Section~\ref{matrix}
into locale theory, and thereby recover the frame homomorphism~$H$.
\end{r@Definition}

Once again we seek to remove the locale or continuous lattice
from the definition,
this time with the goal of eliminating the infinitary joins
in favour of finitary properties of the way-below relation.
Of course, since $\ll$ was itself defined using directed joins,
in Section~\ref{way-below}
it will have to be replaced by an abstract relation $\waybelow$.
This will satisfy axioms based on the following properties,
which are the analogues of Lemmas
\ref{classical interpolation}, \ref{classical 0+<<} and~\ref{classical Wilker}:
    
\begin{r@Lemma}
\label{locale monotone} If $\beta'\leq\beta\ll\phi\leq\phi'$ then
$\beta'\ll\phi'$. \qed
\end{r@Lemma}

\begin{r@Lemma}
\label{locale interpolation} The relation $\ll$ is
\textdf{transitive} and \textdf{interpolative}:
if $\alpha\ll\beta\ll\gamma$ then $\alpha\ll\gamma$,
and conversely given $\alpha\ll\gamma$, there is some $\beta$ with
$\alpha\ll\beta\ll\gamma$. \qed
\end{r@Lemma}

\begin{r@Lemma}
\label{locale 0+<<} $\bot\ll\phi$, and
if $\alpha\ll\phi$ and $\beta\ll\phi$ then $(\alpha\lor\beta)\ll\phi$. \qed
\end{r@Lemma}

\medskip

The Wilker property in Lemma~\ref{classical Wilker} used
excluded middle, but its analogue for continuous lattices is both
intuitionistic and very simple:
it follows from the observation that binary joins distribute
over joins of inhabited, and in particular directed, families.

\begin{r@Lemma}
\label{locale Wilker} In any continuous lattice, if
$\alpha\ll\beta\lor\gamma$ then $\alpha\ll\beta'\lor\gamma'$ for some
$\beta'\ll\beta$ and $\gamma'\ll\gamma$.
\end{r@Lemma}

\Proof Since any directed set is inhabited,
$$ \beta\lor\gamma
  \;=\; \dirsup\collect{\beta'\lor\gamma}{\beta'\ll\beta} \;=\;
  \dirsup\collect{\beta'\lor\gamma'}{\beta'\ll\beta,\;\gamma'\ll\gamma}.
$$
Then if $\alpha\ll\beta\lor\gamma$, we have $\alpha\ll\beta'\lor\gamma'$ for
some term in this join. \qed

\medskip

The relationship between $\ll$ and $\land$ is more subtle.

\begin{r@Lemma}
\label{locale <<*}
If $\alpha\ll\beta\ll\phi$ and $\beta\ll\psi$ then $\alpha\ll\phi\land\psi$.
\qed
\end{r@Lemma}

\smallskip

\begin{r@Definition}
\label{stab loc cpct loc}
A \textdf{stably locally compact locale} is one in which $\top\ll\top$,
and if $\beta\ll\phi$ and $\beta\ll\psi$ then $\beta\ll\phi\land\psi$.
\end{r@Definition}

Examples~\ref{not stably loc cpct} can be adapted to yield
locally compact locales that are not stably locally compact,
and therefore only obey the weaker rule in the Lemma.

\smallskip

\begin{r@Remark}
Stably locally compact objects enjoy many superior properties,
illustrating the duality between compact and open subspaces.
Jung and S\"underhauf \cite{JungA:duacvo} set out ``consistency
conditions'' for them that are similar to ours, 
except that they choose to make $(0,1,{+},{\star})$
a genuine (``strong proximity'') lattice.

However, as Examples~\ref{not stably loc cpct} illustrated,
not all of the locally compact objects that we wish to consider
are stably so, either in geometric topology or recursion theory.
Most obviously for the former,
in $\realno$, the whole space (the trivial intersection) is not compact.

Another reason why we consider the more general situation
is that it corresponds to the monadic Axiom~\ref{monad}
that was the fundamental idea behind the research programme
of which this paper is a part.
This correspondence, which has no counterpart in \cite{JungA:duacvo},
is the main technical goal of this paper;
the consistency conditions are just an intermediate step.

We shall see at the end of Section~\ref{lcpct} that the non-stable situation
also highlights an interesting difference (separate from the usually
mentioned ones of constructivity) between locales and sober spaces
as ways of presenting topological information.

Applying G.H.~Hardy's test \cite{HardyGH:mata},
we may wonder which of the stable and non-stable theories is ``beautiful''
and which is ``ugly''.
My suspicion is that stably locally compact spaces and the
``relational'' morphisms that they describe
(\cite{JungA:stacsc} and our \ref{AD=DA}, \ref{upd}, \ref{upd rk}
and \ref{upd prop})
play a different role in the bigger picture,
whilst we are right to study ``functional'' morphisms in the non-stable case.

The logic that we use is a very weak computational one,
but on closer examination, we see that a great deal of the work
that has been done in domain theory, using many notions of ``basis''
or ``information system'' could actually be formulated in such a logic.
\end{r@Remark}

%============================================================================
\section{Axioms for abstract Stone duality}\label{axioms I}

In this paper we develop a computable account of locally compact sober
spaces and locales, but using a $\lambda$-calculus in place of the
usual infinitary lattice theory, which conflicts with computable ideas.
This calculus, called Abstract Stone Duality,
exploits the fact that, for any such space $X$,
its lattice of open subspaces provides the exponential $\Sigma^X$
in the category.
Here $\Sigma$ is the Sierpi\'nski space
(which, classically, has one open and one closed point),
and the lattice $\Sigma^X$ is equipped with the Scott topology.
Its relationship to the ``consistency conditions'' in the previous sections
and in \cite{JungA:duacvo} will be examined
in Sections~\ref{way-below}--\ref{matrix}.

In this section and the next we present the axiomatisation of ASD
in a reflective fashion that is intended to examine the reasons
for our choices.
For a straightforward summary of the calculus intended for applications,
see \cite{TaylorP:dedras} instead.

\begin{r@Remark}
\label{two models}
We are primarily interested in two particular models of the calculus:
\begin{letterlist}
\item as a source of topological intuition, the classical ones,
  namely the categories of locally compact sober spaces ($\LKSp$)
  and of locally compact locales ($\LKLoc$),
  maybe over an elementary topos other than $\Set$;
\item for computation, the term model.
\end{letterlist}

\smallskip

The classical side has a wealth of concepts motivated by geometry and analysis,
but its traditional foundations are logically very strong, 
being able to define many functions that are neither continuous
nor computable, besides many other famous pathologies.
Having created such a wild theory, we have to rein it back in again,
with a double bridle.
The topological bridle is constructed with infinitary lattice theory,
whilst in recursion theory
we are reduced to using G\"odel numberings of manipulations
of codes for basic elements. % \cite{SmythMB:effgd}.
Abstract Stone duality avoids all of this
by only introducing computably continuous functions in the first place.
(We pay a logical price for this in not being able to define objects
and functions anything like so readily as in set theory.)

On the other hand, logically motivated discussions often read the foundational
aspects more literally than their topological authors ever intended.
They are so bound up in their own questions of what constitutes
``constructivity'' that they lose sight of the conceptual structure
behind the mathematics itself.
For one example, in several constructive approaches to topology,
the closed real interval fails to be compact.
For another,
whilst the \emph{mathematician} in hot pursuit of a proof typically postulates
a \emph{least} counterexample, in order to rule it out,
it is impertinent of the \emph{logician} to emphasise that this uses
excluded middle, since proofs of this kind (once found)
can very often be recast in terms of the induction scheme
\cite[Section 2.5]{TaylorP:prafm}.
As a result, foundational work often fails to reach ``ground level''
in the intended construction.

The lesson that we draw from this is that logic and topology readily drift
apart, if ever we let go of either of them for just one moment.
So we have to tell their stories in parallel.
This means that we must often make do with rough-and-ready versions of
parts of one, in order to make progress with the other.
For example, we use concepts such as compactness to motivate our
$\lambda$-calculus, which itself underlies the technical machinery
that eventually justifies the correspondence with traditional topology
and the use of its language.
We shall make a point of explaining how each step (Definition, Lemma, \etc)
of the new argument corresponds to some idea in traditional topology.

This means that you will have to be prepared for some sudden
switches between contexts.
This paper is highly experimental.
It introduces $\lambda$-calculus formulations of topological concepts
and arguments,
which it has to try out in the traditional categories for topology,
even though these are defective for the purpose.
When we are sure that we have \emph{all} of the axioms in place
(which was certainly not the case before this investigation began,
but is close to being so after it)
the best way to handle the relationship will not be
to interpret ASD in locale theory, but \viceversa~\cite{TaylorP:eletvc}.

In this section we have to set up some \emph{logical} structure
to which the \emph{topologically} minded reader may prefer to return
after first reading Section~\ref{lcpct}.
We catch up on some of the justifications in Sections \ref{corolls}
and~\ref{prime/nucleus},
for example characterising our morphisms in terms of preservation
of meets and joins.
In this section too, we interleave the axioms and definitions of the
$\lambda$-calculus with the technical issues that they are intended to
handle.
\end{r@Remark}

\medskip

\begin{r@Remark}
\label{classify open}
Recall first the universal properties of the Sierpi\'nski space, $\Sigma$.
Any open subspace $U\subset X$ is \textdf{classified}
by a continuous function $\phi:X\to\Sigma$,
in the sense that $U=\phi^{-1}(\top)$,
where $\top$ is the open point of $\Sigma$,
and $\phi$ is unique with this property.
This is summed up by the pullback diagram
\begin{diagram}
  U \SEpbk & \rTo & \terminalobj \\
  \dOpeninto && \dTo>\top \\
  X & \rTo^\phi & \Sigma
\end{diagram}
In our calculus, we shall use the $\lambda$-term $\phi:\Sigma^X$
instead of the subset $U\subset X$.
A similar diagram, using the closed point $\bot\in\Sigma$ instead of $\top$,
classifies (or, as we shall say, \textdf{co-classifies}) the closed subspace,
that, classically, is complementary to~$U$.
We follow topology rather than logic in retaining the bijection between
open and closed subspaces, even though they are not actually complementary
in the sense of constructive set theory.
This is because it is not sets but topological spaces that we wish to capture.

Unions and intersections of open and closed subspaces
make the topology $\Sigma^X$, and in particular the object $\Sigma$,
into distributive lattices ---
honest ones now, not consisting of codes as in Remark~\ref{basis dist latt}.
We therefore have to axiomatise the exponential and the lattice structure.
\end{r@Remark}

\begin{r@Axiom}
\label{monad axiom}
The category $\S$ of ``spaces'' has finite products and an object $\Sigma$
of which all exponentials $\Sigma^X$ exist.
Then the adjunction $\Sigma^\blank\adjoint\Sigma^\blank$
that relates $\opp\S$ to $\S$ is to be \emph{monadic}.
This categorical statement has an associated symbolic form, consisting of
\begin{letterlist}
\item\label{restr lcalc} the simply typed $\lambda$-calculus,
  except that we may only introduce types
  of the form $\Sigma^{X\times Y\times\cdots}$
  (or $X\to Y\to\cdots\to\Sigma$ if you prefer),
  and therefore $\lambda$-abstractions whose bodies already have such types;
\item\label{sober} an additional \emph{term}-forming operation, $\focus$,
  which may only be applied to a term $P$ of type $\Sigma^{\Sigma^X}$
  that satisfies a certain \textdf{primality} equation capturing
  the situation $P=\Lamb\phi.\phi a$, and then $\focus P=a$;
  the use of $\focus$ makes the space $X$ \textdf{sober} \cite{TaylorP:sobsc};
\item\label{monad}
  and an additional \emph{type}-forming operation (with associated term
  calculus) that provides a \emph{formal} subspace
  $i:X\equiv\collect Y E\splitinto Y$ with $I:\Sigma^X\retract\Sigma^Y$
  such that $\Sigma^i\cdot I=\id_{\Sigma^X}$
  given an endomorphism $E:\Sigma^Y\to\Sigma^Y$
  satisfying a certain equation, and then $E=I\cdot\Sigma^I$;
  members of the formal subspace are then the \textdf{admissible} terms $a:Y$,
  \ie those satisfying yet another equation,
  for which we may then introduce $\admit a:X$ \cite{TaylorP:subasd}.
\end{letterlist}
\end{r@Axiom}

\begin{r@Remark}
\label{formulate prime/nucleus}
The equations required in parts (b) and (c) can each be formulated
in two different ways:
an abstract one that states the underlying monadic idea in $\lambda$-notation,
and a lattice-theoretic one that takes advantage of the additional
``topological'' structure of the calculus.
One of our many tasks in this paper (Section~\ref{prime/nucleus})
is to prove the equivalence of these formulations.

In~(b), we shall find that $P:\Sigma^{\Sigma^X}$ is prime iff
it is a lattice homomorphism:
$$ P\top\eq\top \qquad P\bot\eq\bot \qquad
   P(\phi\land\psi)\eq P\phi\land P\psi 
   \quad\hbox{and}\quad
   P(\phi\lor\psi)\eq P\phi\lor P\psi,
$$
whilst for (c), the term $E$ is of the form $I\cdot\Sigma^i$ iff
$$ E(\phi\land\psi) = E(E\phi\land E\psi)
   \quad\hbox{and}\quad
   E(\phi\lor\psi) = E(E\phi\lor E\psi).
$$
We call such a term $E$ a \textdf{nucleus},
shamelessly appropriating this word from locale theory,
in which a nucleus is a monotone endofunction $j:L\to L$ of a frame
such that $\id\leq j=j\cdot j$.
A~nucleus $E$ in our sense must be Scott continuous,
but need not be order-related to the identity,
but the senses agree when $E\geq\id$ and $j$ is Scott continuous.

We discuss $\Sigma$-split subspaces in Section~\ref{subobjects},
and give the admissibility equation in Definition~\ref{def admit},
which is where we first use it in this paper.
\end{r@Remark}

\begin{r@Remark}
\label{normalisation}
 We have said that we are interested in the term model of the calculus.
Now, this is powerful enough to encode computations,
because domain theory can be developed in it (Section~\ref{domain theory})
and programming languages such as Plotkin's PCF can be interpreted
using standard techniques of denotational semantics \cite{PlotkinGD:lcfcpl}.
This means that we cannot hope to give an explicit description
of the term model: the best we can do is to indicate how
the topological features of ASD can be ``compiled out'',
to yield programs in some recognisable programming language.

In \cite[\S\S 9--11]{TaylorP:subasd}, arbitrarily complicated
combinations of $\Sigma$-split subspaces and exponentials
are reduced to a single $\Sigma$-split subspace of
a (still complicated) exponential.
In this paper (Section~\ref{bases}) we shall reduce this even further,
to a $\Sigma$-split subspace of $\Sigma^N$,
showing that this description is equivalent
to giving an $N$-indexed basis for~$X$.
The terms also have a corresponding normalisation.

This leaves $\focus$, the quantifiers, recursion \etc on $\Sigma^\natno$.
The idea of $\focus$ is that it recovers a point of a space
from its Scott-open prime filter of neighbourhoods,
which can be characterised using either infinitary lattice theory
\cite[\S 2]{TaylorP:sobsc}
or the $\lambda$-calculus arising from the monad
\cite[\S 4]{TaylorP:sobsc}.
In this calculus, $\focus$ is pushed to the outside of any term,
and is shown to be redundant at all objects except $\natno$
\cite[\S 8]{TaylorP:sobsc}.

Using the (abstract) lattice structure and primitive recursion,
$\focus$ on $\natno$ is inter-definable with
definition by description \cite[\S\S 9--10]{TaylorP:sobsc}.
General recursion can also be defined.
Finally, \cite[\S 11]{TaylorP:sobsc} sketches how
everything apart from disjunction and recursion strongly normalises,
to what is essentially a \PROLOG\ clause.
Restoring those two features yields a (parallel) \PROLOG\ program.
On the other hand, normalising the $\lambda$-applications
may not be such a good idea, so it's better to think of ASD terms
as parallel \LPROLOG\ programs.

The free model of the calculus therefore satisfies the Church--Turing thesis
(\ie it is equivalent to thousands of other ways of writing programs)
so we do not need to introduce Kleene-style notation with G\"odel numbers
in order to justify calling it a \emph{computable} account of topology.
In particular, the ``computations'' referred to in
Definitions \ref{rec loc cpct}, \ref{rec cts} and~\ref{rec cts latt}
may be defined by terms of our calculus.
\end{r@Remark}

\begin{r@Axiom}
\label{Phoa}
Returning to the Sierpi\'nski space itself,
$(\Sigma,\top,\bot,{\land},{\lor})$
is an internal distributive lattice in $\S$,
which also satisfies the \textdf{Phoa principle},
$$  F:\Sigma^\Sigma,\; \sigma:\Sigma \ \proves\
   F\sigma \;\eq\; F\bot \lor \sigma \land F\top. $$
This equation (which is bracketable either way, by distributivity)
is used to ensure that terms of type $\Sigma^X$ yield data
for the open or closed subspace of $X$,
as required by the monadic axiom \cite[Sections 2--3]{TaylorP:geohol},
and hence the pullback in Remark~\ref{classify open}.
It thereby enforces the bijective correspondences amongst
open and closed subspaces and terms of type~$\Sigma^X$.
We shall assert an ``infinitary'' generalisation of the Phoa principle shortly.
\end{r@Axiom}

\begin{r@Definition}
\label{order}
The lattice structure on $\Sigma$ and $\Sigma^X$ defines
an \textdf{intrinsic order},~$\leq$, where
$\phi\leq\psi$ iff $\phi=\phi\land\psi$ iff $\phi\lor\psi=\psi$.
In the case of terms of type $\Sigma$ (which we call \textdf{propositions}),
we shall write $\Implies$ and $\eq$ instead of $\leq$ and~$=$.
The order is inherited by other objects:
$$ \Gamma\ \proves\ a\leq b:X \quad\hbox{if}\quad
   \Gamma\ \proves\ 
   (\Lamb\phi:\Sigma^X.\phi a)\;\leq\;(\Lamb\phi.\phi b):\Sigma^{\Sigma^X}.$$
There are other ways of defining an order, 
but sobriety and the Phoa principle make them equivalent,
and also say that all maps are monotone \cite[Section 5]{TaylorP:geohol}.

In $\LKLoc$ this order on hom-sets arises from the order on the objects,
considered as lattices, 
whilst in $\LKSp$ it is the \textdf{specialisation order},
$$ x \leq y \;\equiv\;
   (\All U\subset X \hbox{ open}.x\in U\Implies y\in U)
   \;\equiv\; x \in \overline{\setof y},
$$
where $\overline{\setof y}$ is the smallest closed subspace containing~$y$.
\end{r@Definition}

\begin{r@Lemma}
\label{subsp quant} Let $i:U\rOpeninto X$ and $j:C\rClosedinto X$
be the open and closed subspaces (co)classified by $\phi:\Sigma^X$.
Then, with respect to the intrinsic orders on $\Sigma$ and $\Sigma^X$,
there are adjoints
$$ \exists_i\adjoint\Sigma^i \quad\hbox{and}\quad \Sigma^j\adjoint\forall_j
   \vadjust{\nobreak}
$$
that behave like quantifiers \cite[Section 3]{TaylorP:geohol}.
(Note the different tails on the arrows.)
\qed
\end{r@Lemma}

\goodbreak

\begin{r@Remark}
\label{contrav <<}
The order $\leq$ on $\Sigma^X$ presents an important problem
for the axiomatisation of the join in Definition~\ref{def cts latt},
$$ \phi \;=\; \dirsup\collect\beta{\beta\ll\phi}. $$
The condition on the right of the ``$\vert$'' is monotone (covariant) in $\phi$
(indeed, Example~\ref{upup} will use the fact that it is Scott-continuous),
but \emph{contravariant} with respect to $\beta$:
$$ \hbox{if}\quad (\beta'\leq\beta) \quad\hbox{then}\quad
   \big((\beta'\ll\phi)\Impliedby(\beta\ll\phi)\big). $$
We want to regard the topology $\Sigma^X$ as another \emph{space},
but the subset $\collect\beta{\beta\ll\phi}$ is not a sober subspace
(in the Scott topology), since it isn't closed under $\dirsup$.
\end{r@Remark}

\begin{r@Remark}
\label{USet}
This means that we have, after all, to make a distinction between the
exponential space $\Sigma^X$ and the \emph{set} (albeit structured)
of open subspaces of~$X$. We shall write $\USet{\Sigma^X}$ for the
latter, since it is the set of points of the space $\Sigma^X$.
As there are \apriori no ``sets'' in ASD, we have to explain 
what special (properties of) spaces play their role in the theory.
As we shall see, one of these properties is an ``equality test''.

Although locale theory plays down the underlying set functor
$\USet{-}:\Loc\to\Set$, since it is not faithful,
this functor nevertheless exists,
and the subject makes intensive use of $\USet{\Sigma^X}$ in particular,
this being the frame corresponding to the locale~$X$.
The functor $\USet{-}$ may be characterised as the
right adjoint to the inclusion $\Set\to\Loc$ that equips any set
with its discrete topology, \ie the powerset considered as a frame.
In fact, this right adjoint is precisely what we have to add to the
computably motivated axioms of abstract Stone duality given in this paper,
in order to make them agree with
the ``official'' theory of locally compact locales
over an elementary topos \cite{TaylorP:eletvc}
(which writes $\Omega X$ or ${\mathsf U}\Sigma^X$ for $\USet{\Sigma^X}$).
In other words, it distinguishes between the 
two leading models in Remark~\ref{two models}.

This adjunction says that there is a map $\eps:\USet{\Sigma^X}\to\Sigma^X$
that is \emph{couniversal} amongst maps $\beta:N\to\Sigma^X$ from the
objects $N$ of a certain full subcategory,
so $\beta$ factors uniquely as $N\to\USet{\Sigma^X}\to\Sigma^X$.
However,
the couniversal object $\USet{\Sigma^X}$ cannot exist in the computable theory:
besides being uncountable, its equality test would solve the Halting Problem.
So we have to develop alternatives to it.
In traditional language, all we need is that any $\phi\in\Sigma^X$
be expressible as a directed join of $\beta^n$s,
as in Definition~\ref{rec cts latt}.

\medskip

Returning to the problem of Definition~\ref{def cts latt},
we must use \emph{codes} for (basic) open sets,
since we cannot define $\ll$ in abstract Stone duality
using the open set $\beta:\Sigma^X$ itself.
Thus Remark~\ref{basis=cts} will replace the formula $\beta\ll\phi$
(with variables $\beta$ and $\phi$) by
$$ n:N,\; \phi:\Sigma^X \ \proves\ (\beta^n\ll\phi):\Sigma, $$
where $\beta^n$ is the basic open subspace with code~$n$.
\end{r@Remark}

%============================================================================
\section{Sets, unions and bases}\label{axioms II}

Marshall Stone showed us that we should always look for the topology
on mathematical objects, and ASD was called after him partly because
this is its fundamental message too.
However, our discussion of infinitary joins and computable bases
in Definitions \ref{rec loc cpct} and~\ref{def cts latt}
has shown that this topology may still be the discrete one,
\ie that we need sets after all.
Since we have no \emph{sets} as such,
we think of these objects as ``discrete'' spaces.
In ASD we take this particular word to mean that
there is an internal notion of \emph{equality},
but the logical structure that we need to express the join is
\emph{existential quantification}.

\begin{r@Definition}
 An object $N\in\ob\S$ is \textdf{discrete} if there is a pullback
\begin{diagram}
  N \SEpbk & \rTo & \terminalobj \\
  \dOpeninto<\Delta && \dTo>\top \\
  N\times N & \rTo^{(=_N)} & \Sigma
\end{diagram}
\ie the diagonal $N\openinto N\times N$ is \emph{open}.
We may express this symbolically by the rule
$$\hbox{for } \Gamma\proves n,m:N, \qquad
\begin{prooftree}
  \Gamma \proves n = m : N
  \Justifies
  \Gamma \proves (n =_N m) \eq \top : \Sigma
\end{prooftree}
$$
The $(=)$ above the line denotes externally provable equality of terms
of type $N$, whilst $(=_N)$ below the line is the new internal structure.
Categorically,
$(n,m):\Gamma\to N\times N$ and $!:\Gamma\to\terminalobj$
provide a cone for the pullback iff $(n =_N m) \eq \top$,
and in just this case $(n,m)$ factors through $\Delta$,
\ie $n=m$ as morphisms.

We also write
$\setof n:\Sigma^N$ for the singleton \emph{predicate} $\Lamb m.(n=_N m)$
but $\listof n:\List N$ for the singleton \emph{list}.

Beware that this notion of discreteness says that the diagonal and singletons
are open, but not necessarily that \emph{all} subspaces of $N$ are open.
For example, the G\"odel numbers for non-terminating programs
form a closed but not open subspace of $\natno$
\cite{TaylorP:nonagr,TaylorP:pcfasd},
because the classically continuous map $\natno\to\Sigma$ that classifies
it is not computable (Definition~\ref{rec cts}).
\end{r@Definition}

\begin{r@Definition}
\label{def Hdf}
Similarly, an object $H$ is \textdf{Hausdorff} if $H\rClosedinto H\times H$
is \emph{closed}, \ie co-classified by $(\neq_H):H\times H\to\Sigma$, so
$$\hbox{for } \Gamma\proves a,b:H, \qquad
\begin{prooftree}
  \Gamma \proves a = b : H
  \Justifies
  \Gamma \proves (a \neq_H b) \eq \bot : \Sigma.
\end{prooftree}
$$
Notice that equality of the terms $a,b:H$ corresponds to a sort of
doubly negated internal equality, so this definition carries the scent 
of excluded middle.
Like that of a closed subspace,
it was chosen on the basis of topological rather than logical intuition.
As we saw in Example~\ref{PCA} discrete spaces need not be Hausdorff.
\end{r@Definition}

\medskip

\begin{r@Definition}
\label{def overt}
An object $N\in\ob\S$ is \textdf{overt} if $\Sigma^{!_N}$
(where $!_N$ is the unique map $N\to\terminalobj$)
has a left adjoint, $\exists_N:\Sigma^N\to\Sigma$,
with respect to the intrinsic order (Definition~\ref{order}).
Then
$$\hbox{for } \Gamma\proves \sigma:\Sigma,\; \phi:\Sigma^N, \qquad
\begin{prooftree}
  \Gamma,\; x:N \ \proves \phi x \;\Implies\; \sigma:\Sigma
  \Justifies
  \Gamma \ \proves \Some x.\phi x \;\Implies\; \sigma:\Sigma
\end{prooftree}
$$
where we write $\Some x:N.\phi x$ or just $\Some x.\phi x$
for $\exists_N(\Lamb x:N.\phi x)$.
The Frobenius law,
$$ \sigma\land\Some x.\phi x \;\eq\; \Some x.\sigma\land\phi x $$
may be derived from the Phoa principle (Axiom~\ref{Phoa}),
and Beck--Chevalley is also automatic \cite[Section~8]{TaylorP:geohol}.
The lattice dual of this definition is the subject of the next section.
\end{r@Definition}

\begin{r@Axiom}
\label{N axiom}
  $\S$ has a natural numbers object $\natno$, \ie an object that
  admits primitive recursion at all types, and so is discrete and
  Hausdorff.  We also require $\natno$ to be overt.
  
  In order to prove the results that we expect by induction,
  equational hypotheses must be allowed in the context $\Gamma$
  of each judgement.
  See \cite[\S 2]{TaylorP:insema} for a discussion of this.
\end{r@Axiom}

\begin{r@Remark}
\label{use of N} Since the object $N$ over which the codes
range is discrete, its \emph{intrinsic} order $\leq$ is trivial
\cite[Lemma 6.2]{TaylorP:geohol}.
However, the structure $(N,0,1,{+},{\star})$ of an abstract basis is,
at least morally, supposed to be
that of a distributive lattice (Remark~\ref{basis dist latt}).
But this is an ``imposed'' structure,
\ie one that is only defined by the explicit specification of
$(0,1,{+},{\star})$,
rather than by its relationship to the other objects in the category.
We are completely at liberty to consider functions that preserve,
reverse or disregard the associated \emph{imposed} order relation $\baseleq$
that may be defined from $+$ and~$\star$ (Definition~\ref{leq from +*}).
Indeed, the need to use conditions that are contravariant with respect to
the natural order was precisely the reason for distinguishing between
$N$ (or $\USet{\Sigma^X}$) and $\Sigma^X$ in Remark~\ref{contrav <<}.

We typically use the letter $N$ for an overt discrete object, as its
\emph{topological} properties are like those of the natural numbers ($\natno$),
though the foregoing remarks do not give it either arithmetical
or recursive structure.
So when we write $(\beta^n)$ for the basis of open subspaces,
we intend a \emph{family}, not a \emph{sequence}.
The use of the letter $N$ is merely a convention, like $K$ for compact spaces;
the letter $I$ (for \emph{indexing set}) is often used elsewhere
in the situations where we use $N$ below,
but it has acquired another conventional use in abstract Stone duality
(\cf Lemma~\ref{base->subsp}).
The notation and narrative may give the impression of countability,
but bases in the classical model (Remark~\ref{two models}(a))
may be indexed by \emph{any} set, however large you please.
Nevertheless, I make no apology for this impression,
as I~consider $\aleph_1$ and the like to have no place in topology.
I also suspect that occurrences of ``sequences'' and
``countable sense subsets'' in the subject betray the influence of objects 
whose significant property is overtness and not recursion.
\end{r@Remark}

\begin{r@Theorem}
 The full subcategory of overt discrete spaces is a \emph{pretopos},
\ie we may form products, equalisers and stable disjoint unions of
them, as well as quotients by open equivalence relations
\cite[Section 11]{TaylorP:geohol}.
If the ``underlying set'' functor in Remark~\ref{USet}
exists, as in the classical models $\LKSp$ and $\LKLoc$,
then the overt discrete spaces form
an elementary \emph{topos}~\cite{TaylorP:eletvc}. \qed
\end{r@Theorem}

\medskip

The combinatorial structures of most importance to us,
however, are the following:

\begin{r@Theorem}
\label{K+List} Assuming a ``linear fixed point'' axiom
(that any $F:\Sigma^U\to\Sigma^V$ preserves joins of ascending sequences),
every overt discrete object $N$ generates
a free semilattice $\Kur N$ and a free monoid $\List N$ in $\S$,
which satisfy primitive recursion and equational induction schemes,
and are again overt discrete objects \cite{TaylorP:insema}. \qed
\end{r@Theorem}

\medskip

This result is easy to see in the two cases of primary interest,
namely the classical and term models (Remark~\ref{two models}).
In the classical ones ($\LKSp$ and $\LKLoc$), overt discrete spaces are
just sets with the discrete topology, and form a topos.
In this case the general constructions of $\List\blank$ and
$\Kur\blank$ are well known:
$\Kur N$ is often called the \textdf{finite powerset}.

\begin{r@Remark}
\label{N enough} In the term model, on the other hand,
we shall find in Section~\ref{bases} that
$\natno$ itself is adequate to index the bases of all definable objects.
Moreover, any definable overt discrete space $N$ is in fact the
subquotient of $\natno$ by a some open partial equivalence relation
(Corollary~\ref{ovdisc=subquot}).
This both allows us to construct $\List(N)$ and $\Kur(N)$, 
and also to extend any $N$-indexed basis to an $\natno$-indexed one.
There is, therefore, no loss of generality in taking all bases in this
model to be indexed by $\natno$, if only as a method of
``bootstrapping'' the theory.

To construct $\List(\natno)$, we could use encodings of
pairs, lists and finite sets of numbers as numbers.
However, it is much neater to replace $\natno$ with the set $\treeno$
of binary trees.
Like $\natno$, $\treeno$ has one constant ($0$) and one operation,
but the latter is binary (pairing) instead of unary (successor),
and the primitive recursion scheme is modified accordingly.
Hence $\treeno\isomo\setof0+\treeno\times\treeno\isomo\List(\treeno)$,
whereas $\natno\isomo\setof0+\natno$.
The encoding of lists in $\treeno$ has been well known to declarative
programmers since \textsc{Lisp}:
$0$~denotes the empty list, and the ``cons'' $h::t$ is a pair.

Membership of a list is easily defined by list recursion,
as are existential and universal quantification,
\ie finite disjunction and conjunction.
$$
\begin{array}{rclcrcl}
  (\Lamb m.m\in 0)n &\equiv& \bot &&
  (\Lamb m.m\in h::t)n &\equiv& (h=n)\lor(\Lamb m.m\in t)n \\
  \All m\in 0.\phi m &\equiv& \top &&
  \All m\in {h::t}.\phi m &\equiv& \phi h\land\All m\in t.\phi m \\
  \Some m\in 0.\phi m &\equiv& \bot &&
  \Some m\in {h::t}.\phi m &\equiv& \phi h\lor\Some m\in t.\phi m
\end{array}
$$
\end{r@Remark}

\begin{r@Notation}
\label{N notn}
The notation that we actually use in this paper conceals the
preceding discussion.
The letter $N$ may just stand for $\treeno$ in the term model,
but may denote any overt discrete space.
So, in the classical models, $N$ is a set, or an object of the base topos,
equipped with the discrete topology.

We use $\Fin(N)$ to mean $\Kur(N)$ or $\List(N)$ ambiguously.
Since they are respectively the free monoids on $N$ with and without
the commutative and idempotent laws,
this is legitimate so long as their interpretation also obeys these laws.
But this is easy, as the interpretation is usually in $\Sigma^X$,
with either $\land$ or $\lor$ for the associative operation.
In fact,
these two interpretations in $\Sigma^{\Sigma^N}$ are jointly faithful
\cite{TaylorP:insema}.
\end{r@Notation}

\medskip

We are now ready to state the central assumption of this paper.

\begin{r@Axiom}
\label{Scott} The \textdf{Scott principle}:
for any overt discrete object $N$,
$$ \Phi:\Sigma^{\Sigma^N},\; \xi:\Sigma^N \ \proves\
  \Phi\xi \;\eq\; 
   \Some\ell:\Fin(N).\Phi(\Lamb n.n\in\ell) \land \All n\in\ell. \xi n.
$$
Notice that the Phoa principle (Axiom~\ref{Phoa}) is the special case
with $N=\terminalobj$ and so $\Kur(N)=\setof{\initialobj,\terminalobj}$,
whilst the ``linear fixed point'' axiom in Theorem~\ref{K+List}
can easily be shown to be equivalent to the case $N=\natno$.
More generally, the significance of this axiom is that
it forces every object to have and every map to preserve directed joins
(which we have yet to define),
and so captures the important properties that are characteristic
of topology and domain theory.
\end{r@Axiom}

\begin{r@Remark}
\label{exists for bigvee} The lattice $\Sigma^X$ has
intrinsic $M$-indexed joins, for any overt discrete object~$M$.
These are given by $\Lamb x.\Some m:M.\phi^m x$,
and are preserved by any $\Sigma^f$ \cite[Corollary~8.4]{TaylorP:geohol}.

In speaking of such ``infinitary'' joins in $\Sigma^X$,
we are making no additional assertion about lattice completeness:
there are as many joins in each $\Sigma^X$ as there are overt objects,
no more, no fewer.
In particular, there are not ``enough'' to justify impredicative definitions
such as the interior of a subspace, or Heyting implication
(though these can be made in the context of the ``underlying set''
assumption in Remark~\ref{USet} \cite{TaylorP:eletvc}).

Moreover, we use the symbol $\exists$ to emphasise that our joins
are \emph{internal} to our category $\S$,
whereas those in locale theory (written $\bigvee$) are external to $\LKLoc$,
involving the topos ($\Set$) over which it is defined.
This distinction will, unfortunately, become a little blurred
because of the need to compare the ideas of abstract Stone duality
with those of traditional topology and locale theory.
This happens in particular in
the definition of compactness (Definition~\ref{asd def cpct}),
the way-below relation 
(Definition~\ref{def lattice ll} and Remark~\ref{cts=basis})
and the characterisation of finite spaces (Theorem~\ref{Kfinsp})
\end{r@Remark}

\begin{r@Remark}
\label{dep bigvee} When we use $M$-indexed joins,
we shall need $M$ to be a \emph{dependent type},
given, in traditional comprehension notation
(\emph{not} that of \cite{TaylorP:subasd}), by
$$ M \;\equiv\; \collect{n:N}{\alpha_n}\;\subset\; N, $$
where $\alpha_n$ selects the subset of indices $n$
for which $\phi^n:\Sigma^X$ is to contribute to the join.
In practice, this subset is always open, so $\alpha_n:\Sigma$.
The sub- and super-script notation here (and in \cite{TaylorP:eletvc})
indicates that $\phi^n$ typically varies covariantly and $\alpha_n$
contravariantly with respect to an imposed order on~$N$.
Indeed there would be no point in using $\alpha_n$ to select which of
the $\phi$s to include in the join if this had to be an \emph{upper}
subset, as the result would always be the greatest element,
whilst it is harmless to close the subset \emph{downwards}.

This means that, when using the existential quantifier,
we can avoid introducing dependent types by defining
$$ \Some m:\collect{n:N}{\alpha_n}.\phi^m \quad\hbox{as}\quad
   \Some n:N.\alpha_n\land\phi^n.
$$
In Section~\ref{lcpct} we shall refer to
terms like $\alpha_n$ of type $\Sigma$ as \textdf{scalars} and
those like $\phi^n$ of type $\Sigma^X$ as \textdf{vectors}.
\end{r@Remark}

\begin{r@Definition}
\label{DJ} A pair of families
$$ \Gamma,\; s:S \ \proves\ \alpha_s:\Sigma,\quad \phi^s:\Sigma^X $$
indexed by an overt discrete object $S$
is called a \textdf{directed diagram}, and the corresponding
$$ \Some s.\alpha_s\land\phi^s $$
is called a \textdf{directed join} (\cf Definition~\ref{def dirsup}), if
\begin{letterlist}
\item $\alpha_s\eq\top$ for some $s:S$ that we call $0$
 (but see Lemma~\ref{inhab dir}), and
\item $\alpha_{s+t}\eq\alpha_s\land\alpha_t$ and $\phi^{s+t}\geq\phi^t\lor\phi^s$
  for some binary operation ${+}:S\times S\to S$. 
\end{letterlist}
In this, $\alpha_s\land\alpha_t$ means that both $\phi^s$ and $\phi^t$
contribute to the join, so for directedness in the informal sense,
we require some $\phi^{s+t}$ to be above them both (covariance),
and also to count towards the join, for which $\alpha_{s+t}$ must be true.

Although the letter $S$ stands for semilattice, in order to allow
concatenation of lists to serve for $+$ (and the empty list for $0$),
we do not require this operation
to be commutative or idempotent (or even associative).
\end{r@Definition}

\medskip

\begin{r@Remark}
\label{Scott cts rk}
 By the Scott principle, any $\Gamma\proves F:\Sigma^{\Sigma^X}$
preserves directed joins, in the sense that
$$ F(\Some s.\alpha_s\land\phi^s) \;\eq\; \Some s.\alpha_s\land F\phi^s. $$
Notice that $F$ is attached to the ``vector'' $\phi^s$ and not to the scalar
$\alpha_s$,
since the join being considered is really that over the
subset $M\equiv\collect{s}{\alpha_s}$.

This is proved in Theorem~\ref{F:SSX cts}, so it is one of the points on
which the logical proofs lag some way behind the topological intuitions.
Of course, this result could have been used in place of Axiom~\ref{Scott},
but I feel that I made an important point in \cite{TaylorP:fixpps}
by showing how sobriety (actually the slightly weaker notion of repleteness)
transmits Scott continuity from the single object $\Sigma^N$ to
the whole category.
\end{r@Remark}

\medskip

\begin{r@Remark}
\label{+=x} After the separation of directed joins from
\emph{finite} ones ($\bot$, $\lor$),
the behaviour of the latter corresponds much more closely to that
of meets ($\top$, $\land$).
So, although Scott continuity, strictly speaking, breaks the 
lattice duality that we enjoyed in \cite{TaylorP:geohol,TaylorP:nonagr},
we shall still often be able to treat meets and joins at the same time.
We sometimes use the symbol $\sqr$ for either of them,
for example in Lemma~\ref{Hmn+*}.
This means that we can try to transform arguments about open,
discrete, overt, existential, ... things into their lattice duals
about closed, Hausdorff, compact, universal, ... things.
Such lattice dual results seem to be far more common that anyone brought
up with intuitionistic logic or locale theory might expect.
For example, there is a dual basis expansion,
in which $\exists\land$ are replaced by $\forall\lor$,
but we must leave that to later work.
\end{r@Remark}

\begin{r@Lemma}
\label{inhab dir}
For Scott continuity, it suffices that the family be \emph{inhabited},
in the sense that $\Some s:S.\alpha_s\eq\top$,
and have a binary operation $+$ with the properties stated.
\end{r@Lemma}

In this paper, $0$ will typically be the empty list,
but in \cite{TaylorP:dedras,TaylorP:lamcra}
we shall use $\ratno$ with $\max$ or $\min$,
which have no semilattice unit.

\Proof Define $S'\equiv S+\setof 0$, $\alpha_0\equiv\top$ and
$\phi^0\equiv\bot$. Then
$$ \Some s:S'.\alpha_s\land\phi^s \;\eq\;
   \alpha_0\land\phi^0 \;\lor\; \Some s:S.\alpha_s\land\phi^s \;\eq\;
   \Some s:S.\alpha_s\land\phi^s,
$$
but we need $\Some s:S.\alpha_s\eq\top$ for
$$ F\bot \;\eq\; F\bot\land \Some s:S.\alpha_s
 \;\eq\; \Some s:S.\alpha_s\land F\bot
 \;\Implies\; \Some s:S.\alpha_s\land F\phi^s.
 \vadjust{\goodbreak}
$$
Then\closeupaline
\begin{eqnarray*}
   F(\Some s:S.\alpha_s\land\phi^s)
   &\eq& F(\Some s:S'.\alpha_s\land\phi^s) &above\\
   &\eq& \Some s:S'.\alpha_s\land F\phi^s &Remark~\ref{Scott cts rk}\\
   &\eq& \alpha_0\land F\phi^0 \lor \Some s:S.\alpha_s\land F\phi^s \\
   &\eq& \Some s:S.\alpha_s\land F\phi^s. &\qeds
\end{eqnarray*}

Inhabitedness is \emph{all} that we require of a join for it to satisfy
the dual distributive or Frobenius law, \cf Lemma~\ref{locale Wilker}.

\begin{r@Lemma}
\label{or distrib dirsup}
If $\Some n.\alpha_n\eq\top$ then
$ \phi\lor(\Some n.\alpha_n\land\psi^n) =
\Some n.\alpha_n\land(\phi\lor\psi^n)$. \qed
\end{r@Lemma}

The most dramatic effect that Scott continuity has on the presentation
of general topology is the way in which it simplifies the notion of
compactness, to which we now turn.

%============================================================================
\section{Compact subspaces}\label{cpct subsp}

Returning from logic to topology, recall from Definition~\ref{rec loc
cpct} that we need to represent points, open subobjects and compact
subobjects in our calculus.  These will be given by terms of type $X$,
$\Sigma^X$ and $\Sigma^{\Sigma^X}$ respectively, but whilst the first
two correspond exactly, the question of \emph{which} second order
terms denote compact subspaces is more complicated.
 
We pick up the treatment of compact \emph{objects} from \cite[Sections
7--8]{TaylorP:geohol}, and show how this leads to the representation
of compact \emph{sub}objects as quantifiers or modal operators.
Following \cite{HofmannKH:locccl}, these should preserve finite meets
and directed joins, but the basic idea of ASD is that the directed
joins come for free.  A compact subobject is therefore represented by
a term $A:\Sigma^{\Sigma^X}$ for which $A\top\eq\top$ and
$A(\phi\land\psi)\eq A\phi\land A\psi$.
 
This abstract idea works very well as far as the relationship between
compact and closed subspaces of a compact Hausdorff space.  But then
we run into some trouble, because at that point we actually need to
use Scott continuity, but our logical structure will not be strong
enough to do this until Section~\ref{corolls}.  There is an
alternative treatment in \cite{TaylorP:lamcra} that avoids this
problem by restricting to the Hausdorff case, and in particular
$\realno$ and $\realno^\nn$; you may prefer to read that instead of
this section, or in parallel with it.
 
In the experimental spirit,  we try to go further down the path laid by
\cite{HofmannKH:locccl,JungA:duacvo},
but we stumble into more and more potholes, showing how much more work
needs to be done in ASD to make it into a full theory of general
topology.  But Jung and S\"underhauf do have something interesting to
show us when we make this excursion, namely a very close duality
between compact and open subobjects of stably locally compact objects.

The abstract analogue of the way in which locale theory treats local
compactness is much simpler than this.  In that case \emph{any} term
$A:\Sigma^{\Sigma^X}$ may be used --- it doesn't have to preserve
meets.

\begin{r@Remark}
\label{trad cts}
 Traditionally%
   \footnote{In fact Bourbaki \cite[I~9.3]{BourbakiN:topg}
     relegated this formulation to Axiom C$'''$,
     also calling it ``the axiom of Borel--Lebesgue''.
     The older intuitions from analysis involve the existence of
     cluster points of sequences or nets of points, or
     of filters of subsets.},
a topological space $K$ has been defined to be
\textdf{compact} if \emph{every open cover},
\ie family $\collect{U_s}{s\in S}$ of open subspaces such that $K=\bigcup_s U_s$,
\emph{has a finite subcover}, $F\subset S$ with $K=\bigcup_{s\in F} U_s$.
If the family is directed (Definition~\ref{def dirsup})
then $F$ need only be a singleton,
\ie there is already some $s\in S$ with $K=U_s$.

The Scott topology on the lattice $\USet{\Sigma^K}$
offers a simpler way of saying that $K$ is compact.  In this lattice,
$\top$ denotes the whole of $K$, so compactness says that if we can
get into the subset $\setof\top\subset\USet{\Sigma^K}$ by a directed join
$\dirunion_s U_s$ then some member $U_s$ of the family was already
there%
   \footnote{This is therefore an externally defined join,
     \cf Remark~\ref{exists for bigvee}.}.
In other words, $\setof\top\subset\Sigma^K$ is an open subset in the
Scott topology on the lattice.
\end{r@Remark}

\begin{r@Definition}
\label{asd def cpct} In abstract Stone duality we say that
an object $K$ is \textdf{compact} if there is a pullback
\begin{diagram}
  \terminalobj\SEpbk & \rTo & \terminalobj \\
  \dTo<\top && \dTo>\top \\
  \Sigma^K & \rDotsto^{\forall_K} & \Sigma
\end{diagram}
Using the fact that $\Sigma$ classifies open subobjects
(Remark~\ref{classify open}),
together with its \emph{finitary} lattice structure
(but not the Scott principle),
\cite[Proposition~7.10]{TaylorP:geohol} shows that $\forall_K$ exists
with this property iff it is right adjoint to $\Sigma^{!_K}$
(where $!_K$ is the unique map $K\to\terminalobj$).
It is then demonstrated that this map does indeed behave
like a universal quantifier in logic
(the corresponding existential quantifier was
given by Definition~\ref{def overt}).
\end{r@Definition}

\medskip

\begin{r@Lemma}
\label{0+ cpct} If $K$ and $L$ are compact objects then $K+L$
is also compact, as is $\emptyset$.
\begin{diagram}[loose,l>=5em]
  \Sigma^{K+L} \isomo
  \Sigma^K\times\Sigma^L &
  \pile{\rTo^{\forall_K\times\forall_K}\\\top\\
    \lTo_{\Sigma^{!_K}\times\Sigma^{!_L}}}
  & \Sigma\times\Sigma &
  \pile{\rTo^\land\\\top\\\lTo_\Delta} & \Sigma
\end{diagram}
\end{r@Lemma}

\Proof $\forall_{K+L}(\phi,\psi)\eq\forall_K\phi\land\forall_L\psi$ and
$\forall_\emptyset=\top:\Sigma^\emptyset\isomo\terminalobj\to\Sigma$. \qed

\bigbreak

%\goalbreak{5\baselineskip}

The following result is the well known fact that the \emph{direct}
image of any compact subobject is compact (whereas \emph{inverse}
images of open subobjects are open).

\begin{r@Lemma}
\label{Kepi}
Let $K$ be a compact object and $p:K\onto X$ be \textdf{$\Sigma$-epi}%
   \footnote{In \emph{our} category, where every object is a subobject
     of some $\Sigma^Y$, $\Sigma$-epi is the same as epi.
     The name has been inherited from synthetic domain theory,
     a model of which is a topos, only \emph{some} of whose objects
     (namely the predomains) may be so embedded.},
\ie a map for which $\Sigma^p:\Sigma^X\to\Sigma^K$ is mono.
Then $X$ is also compact, with quantifier $\forall_X=\forall_X\cdot\Sigma^p$.
\end{r@Lemma}

\Proof The given quantifier, $\Sigma^{!_K}\adjoint\forall_K$,
satisfies the inequalities
$$ \id \;\leq\; \forall_K\cdot\Sigma^p\cdot\Sigma^{!_X}
    \;=\; \forall_K\cdot\Sigma^{!_K} $$
$$ \Sigma^p\cdot\Sigma^{!_X}\cdot\forall_K\cdot\Sigma^p
   \;=\; \Sigma^{!_K}\cdot\forall_K\cdot\Sigma^p
   \;\leq\; \Sigma^p, $$
from which we deduce $\Sigma^{!_X}\adjoint\forall_K\cdot\Sigma^p$,
since $\Sigma^p$ is mono. \qed

\medbreak

What becomes of the quantifier $\forall_K$ when $K$ no longer stands
alone but is a subobject of some other object $X$?
We see that any compact subobject $K\subset X$ defines a map
$A:\Sigma^X\to\Sigma$ or $A:\Sigma^{\Sigma^X}$
that preserves $\top$ and $\land$.
We shall call this a \textdf{necessity} ($\square$) \textdf{modal operator},
rather than a quantifier, since we have lost
the instantiation or $\forall$-elimination rule $A\phi\Implies\phi x$.

\begin{r@Lemma}
\label{K->A} Let $i:K\to X$ be any map, with $K$ compact.
Then $A\equiv\forall_K\cdot\Sigma^i:\Sigma^X\to\Sigma$ preserves $\top$
and $\land$.
Moreover $A\phi\equiv\forall_K(\Sigma^i\phi)\eq\top$ iff $\Sigma^i\phi=\top_K$.
\end{r@Lemma}

\Proof $\Sigma^i$ preserves the lattice operations
and $\forall_K$ is a right adjoint. \qed

\begin{r@Remark}
\label{A=nhd} In traditional language, for any open set $U\subset X$
classified by $\phi:\Sigma^X$, the predicate $A\phi$ says whether $K\subset U$.
Just as $\forall_K:\Sigma^K\to\Sigma$ classifies $\setof\top\subset\Sigma^K$,
the map $A:\Sigma^X\to\Sigma$ classifies a Scott-open family $\F$
of open subspaces of $X$, namely the family of \emph{neighbourhoods} of~$K$.
Preservation of $\land$ and $\top$ says that
this family is a \textdf{filter}.

Notice that we have
an example of the alternation $K\subset U$ of compact and open subobjects
that we saw in Lemma~\ref{classical interpolation}.
Also, recall from Definition~\ref{rec loc cpct} that we wanted to test
$x\in U$ and $K\subset U$;
these are expressed by the $\lambda$-applications $\phi x$ and $A\phi$
respectively.
\end{r@Remark}

\begin{r@Lemma}
\label{01+ cpct}
\begin{letterlist}
\item If $K=\emptyset$ then $A=\Lamb\phi.\top$,
  so in particular $A\bot\eq\top$.
\item If $K=\setof p$ then $A=\Lamb\phi.\phi p\equiv\eta_X(p)$,
  which is prime (Axiom~\ref{sober})
  and preserves all four lattice operations, in particular $A\bot\eq\bot$.
\item\label{cup cpct}
  If $K\cup L$ exists then $A_{K\cup L}\phi\eq A_K\phi\land A_{L}\phi$.
\item\label{AK contrav} If $K\subset L$ then $A_K\geq A_{L}$.
\item\label{KcapL} Even when $K\cap L$ is compact,
  $A_K\lor A_{L}$ need not be its modal operator.
\end{letterlist}
\end{r@Lemma}

\Proof
\begin{letterlist}
\item $\Sigma^K=\terminalobj$ and $\forall_K=\Lamb\psi.\top$,
  so $A\equiv\forall_K\cdot\Sigma^i=\Lamb\phi.\top$.
\item $\Sigma^K=\Sigma$, $\forall_K=\id$ and $\Sigma^i=\Lamb\phi.\phi p$.
\item Since $K+L\rOnto K\cup L\rSplitinto X$,
  we have $A_{K\cup L}=A_{K+L}=A_K\land A_L$
  by Lemmas \ref{0+ cpct} and~\ref{Kepi}.
\item This follows from the previous part, since $K\cup L=L$.
\item Classically, of course, $K\cap L\subset U$ does not imply
  $(K\subset U\lor L\subset U)$.
  For example, let $K=\setof p$ and $L=\setof q$,
  where $p,q\in X$ are incomparable points in the specialisation order
  (Definition~\ref{order}), so $K\cap L=\emptyset$.
  Then $A_{K\cap L}\bot\eq\top$ but $\bot\eq A_K\bot\lor A_{L}\bot$.
\qed\end{letterlist}

\nobreak

\begin{r@Corollary}
\label{AKbot} $A_K\bot\eq\top$ if $K$ is empty, $\bot$ if it's inhabited. \qed
\end{r@Corollary}

\begin{r@Remark}
\label{inhab cpct} Classically, any compact (sub)space that satisfies
$A_K\bot\eq\bot$ must therefore contain a point, by excluded middle,
but Choice is still needed to \emph{find} it
\cite[Exercise VII 4.3]{JohnstonePT:stos}.
In the localic formulation, it is a typical use of the maximality principle
to find a prime filter containing $\F$,
%\cite[Lemma~III~1.9]{JohnstonePT:stos},
so $A\leq P$ in our notation. \qed
\end{r@Remark}

\bigskip

The situation in which open and compact subobjects interact extremely
well is that of a compact Hausdorff object
(Definitions \ref{def Hdf} and~\ref{asd def cpct}).

\begin{r@Lemma}
 Let $X$ be a compact Hausdorff object and $\phi:\Sigma^X$. Then
\begin{eqnarray*}
  (x\neq y)\lor\phi x &\eq& (x\neq y)\lor\phi y \\
  \phi x &\eq& \All y.(x\neq y)\lor\phi y \\
  (x\neq x) &\eq& \bot\\
  (x\neq y) &\eq& (y\neq x)\\
  (x\neq z) &\Implies& (x\neq y)\lor(y\neq z)
\end{eqnarray*}
\end{r@Lemma}

\Proof These are the lattice duals of
$(x=y)\land\phi x \eq (x=y)\land\phi y$, $\phi x \eq \Some y.(x=y)\land\phi y$,
reflexivity, symmetry and transitivity in any overt discrete object,
and the following proof is just the dual of that in 
\cite[Lemma~6.7]{TaylorP:geohol}.

The closed subobject $\Delta:X\subset X\times X$ is co-classified by $(\neq_X)$,
so the corresponding nucleus,
in the senses of both locale theory ($j\equiv\Delta_*\cdot\Delta^*$)
and of abstract Stone duality ($E\equiv\forall_\Delta\cdot\Sigma^\Delta$,
Axiom~\ref{monad axiom}(c) and Lemma~\ref{subsp quant}) takes
$$ \phi:\Sigma^{X\times X} \quad\hbox{to}\quad
   \Lamb x y.(x\neq y)\lor\phi y. $$
Consider in particular $\psi\equiv\Sigma^{p_0}\phi$,
or $\psi(x,y)\equiv\phi x$; then
$$ \forall_\Delta\phi \;=\; 
   \forall_\Delta\cdot\Sigma^\Delta\cdot\Sigma^{p_0}\phi \;=\;
   \forall_\Delta\cdot\Sigma^\Delta \psi \;=\;
   (\forall_\Delta\bot)\lor\Sigma^{p_0}\phi.$$
The same thing in $\lambda$-calculus notation, and its analogue for $p_1$, are
$$ (\forall_\Delta\phi)(x,y) \;\eq\; (x\neq y)\lor \phi x
    \;\eq\; (x\neq y)\lor \phi y. $$
Now apply $\forall_{p_0}\equiv\forall y$, so 
$ \phi \;=\; \forall_{p_0}\cdot\forall_\Delta \phi
  \;=\; \forall_{p_0}(\forall_\Delta\bot\lor\Sigma^{p_0}\phi)$,
\ie
$ \phi x \;\eq\; \All y.(x\neq y)\lor\phi y $.
Inequality is irreflexive by definition.
We deduce symmetry by putting $\phi\equiv\Lamb u.(y\neq u)$
and the dual of the transitive law  with $\phi\equiv\Lamb u.(u\neq z)$. \qed

\bigskip
   
The idea of the following construction is that
$K\subset U$ iff $U\cup V=X$,
and conversely $x\in V$ iff $K\subset X\setminus\setof x$, where
$K$ is compact and $V$ is its complementary open subobject,
encoded by $A$ and $\psi$ respectively.
The lattice dual of this result --- that open and overt subobjects
coincide in an overt discrete object --- was proved for $\natno$
in \cite[Section 10]{TaylorP:sobsc}.

\begin{r@Proposition}
\label{cpct=closed} In any compact Hausdorff object $K$,
there is a retraction $\Sigma^K\retract\Sigma^{\Sigma^K}$ given by
\begin{eqnarray*}
  \psi &\mapsto& \Lamb\phi.\All x.\psi x \lor \phi x \\
  A &\mapsto& \Lamb x.A(\Lamb y.x\neq y),
\end{eqnarray*}
where, if $A$ is so defined from $\psi$ then it preserves $\top$ and $\land$.
\end{r@Proposition}

Later we shall show that closed and compact subobjects agree exactly.

\Proof $\psi\mapsto A\mapsto\Lamb x.A(\Lamb y.\psi y\lor x\neq y)=\psi$,
by the second part of the Lemma.
Then $A\top\eq\top$ easily,
whilst $A(\phi_1\land\phi_2)\eq A\phi_1\land A\phi_2$ by distributivity.

On the other hand,
$A\mapsto\psi\mapsto\Lamb\phi.\All x.A(\Lamb y.\phi y\lor x\neq y)$,
whilst the first part of the Lemma says that
$A=\Lamb\phi.A(\Lamb y.\All x.\phi y\lor x\neq y)$,
so for the bijection between closed and compact objects we need
$A$ to commute with $\forall x$. 
To show this, we need to know about the Tychonov product topology
on $X\times X$ (Remark~\ref{quantifiers commute}),
and to use Scott continuity (Theorem~\ref{F:SSX cts}). \qed

\bigskip

As we said, it is not obligatory that the $A$s that we shall use to define
locally compact objects in the next section actually satisfy these conditions,
\ie preserve $\top$ and $\land$.
Consider the localic situation.

\begin{r@Example}
\label{upup} For $\beta\in L$ in a continuous lattice
(Definition~\ref{def cts latt}), the subset
$$ \upup\beta\;\equiv\;\collect\phi{\beta\ll\phi}\;\subset\; L $$
is Scott-open, and therefore classified by some $A:\Sigma^L$
(Remark~\ref{classify open}).
However, $\upup\beta$ need not be a filter
(\cf Definition~\ref{stab loc cpct loc}),
so $A$ need not preserve $\top$ or $\land$.  \qed
\end{r@Example}

\medbreak

We obtain similar behaviour even in traditional topology.
The following may seem a strange thing to do,
but it will fall into place as we start to use the $\lambda$-calculus.

\begin{r@Notation}
\label{Hayo notation}
In \cite{TaylorP:sobsc} we found it useful to regard \emph{any} map
$F:\Sigma^X\to\Sigma^Y$ (not necessarily a homomorphism for the monad 
or of frames, but nevertheless Scott-continuous) as a ``second class''
map $\hayo F:Y\Hto X$,
and to write $\Hayo\S$ for the category composed of such maps.
The work cited there explains how they interpret ``control effects''
such as jumps in programming languages.
We shall in particular meet $I:\Sigma^X\to\Sigma^Y$
such that $\Sigma^i\cdot I=\id$, where $i:X\rSplitinto Y$.
\end{r@Notation}

\begin{r@Remark}
\label{image compact} Just as in Lemma~\ref{K->A}
we formed the modal operator corresponding to the direct image $i K$
of $K$ along $i:K\to X$ as the composite $A=\forall_K\cdot\Sigma^i$,
so we may form the direct image $\hayo F A$
along a second class map $\hayo F:K\Hto X$ as $\forall_K\cdot F$.
Its open neighbourhoods are given, as in Remark~\ref{A=nhd}, by
$$ (\hayo F K\subset\phi) \;\eq\; 
   (\forall_K\cdot F)\phi \;\eq\;
    A(F\phi)  \;\eq\; 
   (K\subset F\phi).
$$
However, $\forall_K\cdot F$ need only be a filter when $F$ preserves
$\top$ and $\land$.

\bigbreak

Even when $A$ does preserve meets, and so classifies a Scott-open filter
of open subspaces, it need not correspond to a (locally) compact subspace.
(It does in a compact Hausdorff object, but even there we do not yet
have to tools to prove it, \cf Proposition~\ref{cpct=closed}.)
We make a brief excursion into classical topology to illustrate the
duality of open and compact subspaces,
and their alternating inclusions (Notation~\ref{spatial <<}).
Recall from Definition~\ref{pt loc cpct}
that any open subspace is the union of the compact subspaces inside it:
the first result answers the dual question of when a compact subspace is
the \emph{intersection} of the open subspaces that contain~it.
\end{r@Remark}

\begin{r@Lemma}
\label{cpct sat} Using excluding middle in classical topology
\cite[p221 def 2.1]{HofmannKH:locccl},
$$ \collect{y\in X}{\Some x\in K.x\leq y}
   \;=\; \dirinter\collect{U\subset X \hbox{ open}}{K\subset U}, $$
where $\leq$ is the specialisation order (Definition~\ref{order}).
We call this the \textdf{saturation} of $K$.
Hence compact subspaces of $X$ can be recovered from their modal operators
iff $X$ is ${\mathrm T}_1$ (when $\leq$ is trivial).
For a specific non-example, let $X\equiv\Sigma$ and $K\equiv\setof\bot$,
so $A\equiv\Lamb\phi.\phi\bot$ and the saturation of $K$ is $\Sigma$. \qed
\end{r@Lemma}

\medskip

What we might hope to recover from the modal operator $A$, therefore, is
the saturation of $K$, as Karl Hofmann and Michael Mislove did for
sober spaces \cite[Theorem~2.16]{HofmannKH:locccl},
and Peter Johnstone did for locales \cite{JohnstonePT:viells}.
This result shows that we have identified enough of the properties 
of modal operators, at least in the classical model.

\begin{r@Lemma}
\label{HM0} Let $\F\subset\USet{\Sigma^X}$ be a Scott-open
filter of open subspaces of a sober (but not necessarily locally
compact) space.  Then, assuming the axiom of choice,
$K\equiv\dirinter\F\subset X$ is a compact subspace,
and $\F=\collect{U}{K\subset U}$. \qed
\end{r@Lemma}

\medskip

\begin{r@Lemma}
\label{codir meet cpct} 
\cite[p221 thm 2.16]{HofmannKH:locccl}
Let $K\equiv\dirinter_s K_s$
be a co-directed intersection of compact saturated subspaces of a sober space.
Then $K$ is also compact saturated, and $A_K=\dirsup A_{K_s}$.
If all of the $K_s$ are inhabited then so is $K$. \qed
\end{r@Lemma}

\begin{r@Theorem}
\label{cpct <<}
\cite[Problem 5 F(a)]{KelleyJL:gent}
Compact saturated subspaces of a locally
compact sober space form a continuous preframe under reverse inclusion.
That is, for any compact saturated subspace $K\subset X$,
$$ K \;=\; \dirinter\collect{L \hbox{ compact saturated}}{L\waybelow K}. $$
Here $L\waybelow K$ means that there is an open subspace $U$
with $K\subset U\subset L$ (\emph{sic} --- Notation~\ref{spatial <<}),
but this is equivalent to $L\ll K$,
the order-reversed analogue of Definition~\ref{def lattice ll}, \ie
$$ K \supset \dirinter_s M_s \ \Implies\ \Some s.\ K\supset M_s. \eqno\qEd $$
\end{r@Theorem}

\begin{r@Corollary}
\label{classical dual Wilker}
Using Choice in classical topology,
stably locally compact sober spaces enjoy the dual Wilker property
(\cf Lemma~\ref{classical Wilker}):
if $K\cap L\subset U$ then there are open subspaces $V$ and $W$
and compact ones $K'$ and $L'$ such that
$K'\cap L'\subset U$,  $K\subset V\subset K'$ and $L\subset W\subset L'$
\cite[Theorem 23]{JungA:duacvo}.
\end{r@Corollary}

\Proof Use Lemma~\ref{locale Wilker} in the continuous lattice of
compact saturated subspaces under reverse inclusion.
It is significant in this that the abstract joins in the lattice be
given by intersections of subspaces. \qed

\medskip

\begin{r@Remark}
\label{A not loc cpct} There is still a problem.
Even when $X$ itself is locally compact,
and we have a filter $\F$ that is Scott-open and therefore classified
by a term $A$ of type $\Sigma^{\Sigma^X}$,
the compact subspace $K$ that they define need not be \emph{locally} compact.
Therefore they need not be definable as actual objects
in the theory that we describe in this paper.
However, an extension to and beyond the category of all locales is envisaged
for later work.
Preliminary investigations in this suggest
that any $A:\Sigma^{\Sigma^X}$ that preserves $\top$ and $\land$
can be expressed as the universal quantifier of a compact subobject.

However, as we said at the beginning of this section,
the ``dual bases'' that we shall introduce in the next section
use the $A$s and not the~$K$s.
It is a \emph{desirable} property of $A$ that it be a filter,
\ie preserve $\top$ and $\land$,
because then we can fully exploit the intuitions of traditional topology.
But it is not \emph{necessary}:
locale theory gives rise to bases that don't have this property,
and both ways of doing things have their advantages.
\end{r@Remark}

%============================================================================
\section{Effective bases}\label{lcpct}

This section introduces the first of our abstract characterisations
of local compactness %(Remark~\ref{summary loc cpct}(d))
((d) in the Introduction)
and explores its relationship to the traditional definitions
in point-set topology and locale theory.

\begin{r@Remark}
  Recall that we call a system $(U^n)$ of vectors in a \emph{vector} space
a \emph{basis} if any other vector $V$ can be expressed as a sum of
scalar multiples of basic vectors.  Likewise, we say that a system
$(U^n)$ of open subspaces of a \emph{topological} space is a basis
if any other open set $V$ can be expressed as a ``sum'' (union or
disjunction) of basic opens.

How do we find out which basis elements contribute to the sum, and (in
the case of vector spaces) by what scalar multiple?  By applying the
\emph{dual basis} $A_n$ to the given element~$V$, giving $A_n\cdot V$.
Then
$$V  \;=\;  \sum_n  A_n\cdot V  *  U^n$$
\vskip-\baselineskip\noindent
where
\begin{letterlist}
\item $\sum$ denotes linear sum, union, disjunction or existential
  quantification;
\item ``scalars'' in the case of topology range over the Sierpi\'nski space;
\item the dot denotes
  \begin{itemize}
  \item inner product of a dual vector with a vector to yield a scalar,
  \item that $V$ is an element of the family classified by $A_n$, or
  \item $\lambda$-application; and
  \end{itemize}
\item `$*$'  denotes multiplication by a scalar of a vector, or conjunction.
\end{letterlist}
\end{r@Remark}

\smallskip

\noindent In abstract Stone duality,
since the application of $A_n$ to a predicate $V:\Sigma^X$ yields a scalar,
it must have type $\Sigma^{\Sigma^X}$.
In the previous section we saw that (some) such terms play the role of
compact subobjects.

\begin{r@Definition}
\label{eff basis} An \textdf{effective basis} for an object
$X$ is a pair of families
$$ n:N\;\proves \;\beta^n:\Sigma^X  \qquad 
   n:N\;\proves\; A_n:\Sigma^{\Sigma^X},  $$
where $N$ is an overt discrete object (Section~\ref{axioms II}),
such that every ``vector'' $\phi:\Sigma^X$ has a \textdf{basis decomposition},
$$ \phi \;=\; \Some n. A_n\phi \land \beta^n. $$
This is a join $\bigvee\collect{\beta^n}{\alpha_n}$ over a subset,
with $\alpha_n\equiv A_n\phi$, in the sense of Remark~\ref{dep bigvee}.
In more traditional topological notation, this equation says that
$$ \hbox{for all } V\subset X \hbox{ open},\qquad
  V \;=\; \bigcup\collect{U^n}{V\in\F_n}, $$
where $\beta^n,\phi:\Sigma^X$ and $A_n:\Sigma^{\Sigma^X}$
classify open $U^n,V\subset X$ and $\F_n\subset\Sigma^X$
in the sense of Remark~\ref{classify open}.

We call $(\beta^n)$ the \textdf{basis} and $(A_n)$ the \textdf{dual basis}.
The reason for saying that the basis is ``effective'' is that it is
accompanied by a dual basis, so that the coefficients are given
effectively by the above formula ---
it is not the computational aspect that we mean to stress at this point.
The sub- and superscripts indicate the co- and contravariant behaviour
of compact and open subobjects respectively with respect to continuous maps.
\end{r@Definition}

\begin{r@Remark}
\label{ext def loc cpct}
We shall find in Section~\ref{bases} that every object that is definable 
in the ASD calculus
(as we set it out in Sections~\ref{axioms I}--\ref{axioms II})
has an effective basis.
However, there are plans to extend this calculus to and beyond general locales,
where those results will no longer apply.
Then we shall \emph{define} a \textdf{locally compact} object 
to be one that has an effective basis in the above sense.
\end{r@Remark}

\medskip

The first observation that we make about this definition expresses the
inclusion $U^n\subset K^n$ (\cf Definition~\ref{UKpair}).
After that we see some suggestion of the role of compact subobjects,
although Lemma~\ref{Anbn=1} is too specific to be of much use,
unless $X$ happens to have a basis of compact open subobjects,
\ie it is a \emph{coherent object} \cite[Section II~3]{JohnstonePT:stos}.

\begin{r@Lemma}
\label{cpct base} If $\Gamma\proves\phi:\Sigma^X$ satisfies
$\Gamma\proves A_n\phi\eq\top$ then $\Gamma\proves \beta^n\leq\phi$.
\end{r@Lemma}

\Proof Since $A_n\phi\eq\top$,
the basis decomposition for $\phi$ includes $\beta^n$ as a disjunct. \qed

\begin{r@Lemma}
\label{Anbn=1} If $\proves A_n\beta^n \eq \top$ then $\beta^n$
classifies a compact open subobject.
\begin{diagram}
  \terminalobj\SEpbk && \rTo && \terminalobj \\
  \dTo<{\beta^n} &&&& \dTo>\top \\
  \Sigma^K=\Sigma^X\commacat\beta^n \shift -2em & \rSplitinto
      & \Sigma^X  & \rTo^{A_n} & \Sigma
\end{diagram}
\end{r@Lemma}

\Proof The equation $A_n\beta^n \eq \top$ says that the square commutes.
Any test map $\phi:\Gamma\to\Sigma^X\commacat\beta^n$ that
(together with ${!}:\Gamma\to\terminalobj$)
also makes a square commute must satisfy
$ \Gamma\proves A_n\phi\eq\top $ and $ \Gamma\proves\phi\leq\beta^n $,
but then $\phi=\beta^n$ by the previous result.
Hence the square is a pullback, whilst $\beta^n=\top_{\Sigma^K}$,
so the lower composite is $\forall_K$, making $K$ compact. \qed

\bigbreak

\goalbreak{8\baselineskip}

The following jargon will be useful:

\begin{r@Definition}
\label{basis jargon}
An effective basis $(\beta^n,A_n)$ for an object $X$ is called
\begin{letterlist}
\item a \textdf{directed} or $\lor$-\textdf{basis} if
  there is some element (that we call $0\in N$) such that
  $$ \beta^0=\Lamb x.\bot \quad\hbox{and}\quad A_0=\Lamb\phi.\top $$
  (though $A_0=\top$ entails $\beta^0=\bot$ by Lemma~\ref{cpct base})
  and a binary operation ${+}:N\times N\to N$ such that
  $$ \beta^{n+m} \;=\; \beta^n\lor\beta^m \qquad\hbox{and}\qquad
  A_{n+m}  \;=\; A_n\land A_m;
  $$
  this definition is designed to work with Definition~\ref{DJ};
  it is used first in Lemma~\ref{SX basis} and then
  extensively in Sections \ref{way-below}--\ref{matrix};
\item an \textdf{intersection} or \textdf{$\land$-basis} if
  $\beta^1=\Lamb x.\top$ for some element (that we call $1\in N$),
  and there is a binary operation $\star$ such that
  $$ \beta^{n\star m} \;=\; \beta^n\land\beta^m \qquad
     A_n\;\leq\; A_{n\star m} \quad\hbox{and}\quad
     A_m\;\leq\; A_{n\star m}, $$
  so the intersection of finitely many basic opens is basic;
  this is a positive way of saying that we do have a basis
  instead of what is traditionally known as a \textdf{sub-basis};
\item a \textdf{lattice basis} if it is both $\land$ and $\lor$;
\item a \textdf{filter basis} if each $A_n$ preserves $\land$ and~$\top$,
  and so corresponds in classical topology
  to a compact saturated (though not necessarily locally compact)
  subspace $K^n$, by Lemma~\ref{HM0};
\item a \textdf{prime basis} if each $A_n$ of the form $A_n\phi\eq\phi p^n$
  for some $p^n:X$ (\cf Axiom~\ref{sober}),
  the corresponding compact subobject being $K^n=\setof{p^n}$.
\end{letterlist}
\end{r@Definition}

Any effective basis can be ``up-graded'' to a lattice basis
by formally adding unions and intersections (Lemma~\ref{make or basis}ff).
It can be made into a filter $\lor$-basis instead (Remark~\ref{lacuna}).
Some of the other terminology is only applicable in special situations:
if $A_1\top\eq\top$ then the object is compact, with $\forall_X=A_1$,
whilst a sober topological space has a prime basis
iff it is a continuous dcpo with the Scott topology
(Theorem~\ref{prime basis}).
Even when the intersection of two compact subobjects is compact,
there is nothing to make $A_{n\star m}$ correspond to it,
but we shall rectify this in Proposition~\ref{ASD dual Wilker}.

\begin{r@Examples}
\label{N+SN} Let $N$ be an overt discrete object. Then
\begin{letterlist}
\item\label{N basis} $N$ has an $N$-indexed prime basis given by
  $$ \beta^n\;\equiv\;\setof n\;\equiv\;\Lamb x.(x=_N n)
  \quad\hbox{and}\quad
  A_n\;\equiv\;\eta_N(n)\;\equiv\;\Lamb\phi.\phi n,
  $$
  because
  $ \Some n.A_n\phi\land\beta^n x
    \;\equiv\; \Some n.\eta n\phi\land\setof n x
  \;\eq\; \Some n.\phi n\land (x=n) \;\eq\; \phi x$; and
\item\label{SN basis}
  $\Sigma^N$ has a $\Fin(N)$-indexed prime $\land$-basis given by
  $$ B^\ell \equiv \Lamb\xi.\All m\in\ell.\xi m
  \quad\hbox{and}\quad
  \A_\ell \equiv \Lamb \Phi.\Phi(\Lamb m.m\in\ell),
  $$
  because $\Phi\xi \eq \Some\ell.\Phi(\Lamb m.m\in\ell)\land\All m\in\ell.\xi m$
  by Axiom~\ref{Scott}.
  (The convention that superscripts indicate covariance
  (Remarks \ref{contrav <<} and~\ref{dep bigvee})
  means that the imposed order on $\Fin(N)$ here
  is \emph{reverse} inclusion of lists.)
\end{letterlist}
\end{r@Examples}

\bigskip

\goalbreak{5\baselineskip}

We devote the remainder of this section to showing that every locally
compact sober space or locale has an effective basis in our sense.
We start with traditional topology.

\begin{r@Proposition}
\label{basis=K}
Any locally compact sober space has a filter basis.
\end{r@Proposition}

\Proof Definitions \ref{pt loc cpct} and~\ref{UKpair} provide
families $(K^n)$ and $(U^n)$ of compact and open subspaces such that
$$ \hbox{for each open } V\subset X, \qquad
    V \;=\; \bigcup \collect{U^n}{K^n \subset V}. $$
As the subspace $K^n$ is compact, Remark~\ref{A=nhd} defines
$A_n:\Sigma^X\to\Sigma$ such that $K^n\subset V$ iff $A_n\phi\eq\top$,
where $\beta^n$ and $\phi:\Sigma^X$ classify $U^n$ and $V\subset X$,
so $ \phi = \Some n.A_n\phi\land\beta^n$. \qed

\medskip

Notice how the basis decomposition ``short changes'' us, for
individual basis elements: we ``pay'' $K^n\subset V$ but only
receive $U^n\subset K^n$ as a contribution to the union.
Nevertheless, the interpolation property (Lemma~\ref{classical interpolation})
ensures that we get our money back in the end.

In many examples, $U^n$ may be chosen to be \emph{interior} of $K^n$,
and $K^n$ the closure of $U^n$.
However, this may not be possible if we require an $\lor$-basis.
For example, such a basis for $\realno$ would have as one of its members
a pair of touching intervals, $(0,1)\cup(1,2)$,
which is not the interior of $[0,1]\cup[1,2]$.

\begin{r@Remark}
\label{fbasis=K} It is difficult to identify a substantive
Theorem by way of a converse to this in traditional topology,
since the $\lambda$-calculus can only by interpreted in topological spaces
if they are already locally compact, and therefore have effective bases.
Nevertheless, we can show that a filter basis $(\beta^n,A_n)$ 
can only arise in the way that we have just described.

By Lemma~\ref{HM0},
each $A_n$ corresponds to some compact saturated subspace $K^n$, where
$$ (K^n\subset V)\iff A_n\phi\;\Longrightarrow(U^n\subset V), $$
$\beta^n$ and $\phi$ being the classifiers for $U^n$ and $V$ as usual.
Since $K^n=\dirinter\collect{V}{K^n\subset V}$, we must have $U^n\subset K^n$.
Then, given $x\in V$, so $\phi x\eq\top$,
the basis decomposition $\phi x\eq\Some n.A_n\phi\land\beta^n x$,
means that $x\in U^n\subset K^n\subset V$,
as in Definition~\ref{pt loc cpct}. \qed
\end{r@Remark}

\bigskip

\goalbreak{5\baselineskip}

\begin{r@Example}
\label{eg R2} The closed real unit interval has a filter
$\land$-basis with
\begin{eqnarray*}
  \beta^{q\pm\epsilon} &\equiv& (q\pm\epsilon)\;\equiv\;
  \Lamb x.{{\left|x-q\right|\lt\epsilon}} \\
  A_{q\pm\delta} &\equiv& [q\pm\delta]\;\equiv\;
  \Lamb\phi.\All x.{{\left|x-q\right|\leq\delta}}\Implies\phi x
\end{eqnarray*}
where $\epsilon,\delta\gt0$ and $q$ are rational, and we re-deploy the
interval notation of Example~\ref{eg R1} in our $\lambda$-calculus.
We also write $<q\pm\delta>$ for a variable that ranges over the codes,
as opposed to the open $(q\pm\delta)$ and compact $[q\pm\delta]$
intervals that it names.
The imposed order is given by comparison of the radii $\epsilon$ or $\delta$.
\end{r@Example}

\Proof Let $x\in U\subset[0,1]$; then, for some $\epsilon\gt 0$,
$\All y.{{\left|x-y\right|\lt\epsilon}}\Implies y\in U$.
So with $\delta\equiv\half\epsilon$, and
$q$ rational such that ${\left|x-q\right|\lt\delta}$,
$$ \All y.{{\left|y-q\right|\leq\delta}}\Implies y\in U\quad\land\quad
   x\in\collect y{{{\left|y-q\right|\lt\delta}}}, $$
which is $\Some<q\pm\delta>.A_{q\pm\delta}\phi\land\beta^{q\pm\delta}x$
in our notation.  \qed

\begin{r@Example}
\label{cpct Hdf base}
Recall that any compact Hausdorff space $X$ has a stronger property
called \textdf{regularity}: if $C\subset X$ is closed with $x\notin C$
then there are $x\in U\disjoint V\supset C$ with $U,V\subset X$ open and
disjoint.
Writing $K=X\setminus V$ and $W=X\setminus C$ for the complementary
compact and open subspaces, this says that given $x\in W$,
we can find $x\in U\subset K\subset W$, as in Definition~\ref{pt loc cpct}.
Let $(\beta^n,\gamma_n)$ classify a sufficient computable family
$(U^n\disjoint V_n)$ of disjoint open pairs,
and put $A_n\equiv\Lamb\phi.\All x.\phi x\lor\gamma_n x$,
which corresponds to the compact complement of~$V_n$
(Proposition~\ref{cpct=closed}).
Then $(\beta^n,A_n)$ is a filter basis for~$X$.
It is a lattice basis if the families $(\beta^n)$ and $(\gamma_n)$
are sublattices, with $\gamma_0=\top$, $\gamma_1=\bot$, 
$\gamma_{n+m}=\gamma_n\land\gamma_m$ and
$\gamma_{n\star m}=\gamma_n\lor\gamma_m$. \qed
\end{r@Example}

\medskip

\begin{r@Remark}
 In these idioms of topology, where we say that there
``exists'' an open or compact subobject within certain bounds,
that subobject may typically be chosen to be \emph{basic},
and the existential quantifier in the assertion ranges over
the set (overt discrete object) $N$ of \emph{indices} for the basis,
rather than over the topology $\Sigma^X$ itself.
For example, this is the case with the Wilker property in
Lemmas \ref{classical Wilker} and~\ref{locale Wilker}.
\end{r@Remark}

\bigskip

Now we turn to locale theory.

\begin{r@Proposition}
\label{basis=cts} Any locally compact locale has a lattice
basis.
\end{r@Proposition}

\Proof The localic definition is that the frame $L$ be a distributive
continuous lattice (Definition~\ref{def cts latt}), so
$$ \hbox{for all } \phi\in L, \qquad
   \phi \;=\; \bigvee \collect{\beta\in L}{\beta\ll\phi},
   \vadjust{\nobreak}
$$
and by Example~\ref{upup}, $\upup\beta\equiv\collect\phi{\beta\ll\phi}$
is classified by some $A:\Sigma^L$.

\goodbreak

This means that there is a basis decomposition 
$$ \phi \;=\; \Some n. \beta^n \land A_n\phi, 
   \qquad\hbox{where}\quad
   A_n \;\equiv\;\Lamb\phi.(\beta^n\ll\phi), $$
so $(\beta^n,A_n)$ is a lattice basis.

Recall, however, from Remark~\ref{contrav <<} that
we must consider this basis to be indexed by the \emph{underlying set},
$\USet L$, of the frame $L$.
Thus $N\equiv\USet L$, whilst $\beta^\blank:N\to L$ is the
couniversal map from an overt discrete object $N$ to $L$.

In fact it is enough for the image of $N$ to generate $L$
under directed joins (Definition~\ref{rec cts latt}).
There is no need for $N$ to be the \emph{couniversal} way of doing this,
and so no need for the underlying set functor $\USet-$. \qed

\medskip

\begin{r@Remark}
\label{cts=basis}
The first part of the converse to this is Lemma~\ref{Anbn=1},
but with the relative notion $\ll$ in place of compactness itself:
if $A_n\phi\eq\top$ then $\beta^n\ll\phi$
in the sense of Definition~\ref{def lattice ll}.
For suppose that $\phi\leq\dirsup_s\theta_s$, so
$\top \;\eq\; A_n\phi \;\Implies\; A_n\dirsup_s\theta_s
\;\eq\; \dirsup_s A_n\theta_s$.
Then%
  \footnote{This involves the same confusion of internally and
  externally defined joins as in Definition~\ref{asd def cpct},
  \cf Remark~\ref{exists for bigvee}.}
$A_n\theta_s\eq\top$ for some $s$,
so $\beta^n\leq\theta_s$ by Lemma~\ref{cpct base}.

Now suppose that a locale carries an effective basis in our sense. Then
$$ \phi \;=\; \Some n.A_n\phi\land\beta^n \;\equiv\;
   \bigvee\collect{\beta^n}{A_n\phi} \;\leq\;
   \dirsup\collect{\beta^n}{\beta^n\ll\phi} \;\leq\; \phi, $$
in which the second join is directed, by Lemma~\ref{locale 0+<<}.
Hence the frame $L$ is continuous. \qed
\end{r@Remark}

\begin{r@Remark}
\label{beta not inj}
Notice that $\beta^m\leq\beta^n$ does not
imply any relationship between $n,m\in N$,
because the function $\beta^\blank:N\to L$ need not be injective.
This is the reason why $A_n$ need not be exactly $\Lamb\phi.\beta^n\ll\phi$,
\cf Remark~\ref{basis dist latt}.
\end{r@Remark}

\begin{r@Proposition}
\label{slc=latt filt}
 A locale is stably locally compact iff it has a lattice
filter basis.
\end{r@Proposition}

\Proof $[\Implies]$ Let $A_n=\Lamb\phi.\beta^n\ll\phi$ as before.
$[\Impliedby]$ As we have an $\lor$-basis,
the basis expansion is a directed join,
to which we may apply the definition of $\beta^m\ll\phi$.
In this case there is some $n$ with $\beta^m\leq\beta^n$ and $A_n\phi\eq\top$,
so $\beta^m\ll\top$.
Also, if $\beta^m\ll\phi$ and $\beta^m\ll\psi$ then, for some $p$, $q$,
$$ \beta^m\leq\beta^p,\quad \beta^m\leq\beta^q,\quad
   A_p\phi\eq\top \hbox{ and } A_q\psi\eq\top, $$
so $\beta^m\leq\beta^p\land\beta^q\equiv\beta^{p\star q}$.
As we have a filter $\land$-basis,
$A_{p\star q}\phi\land A_{p\star q}\psi\eq
A_{p\star q}(\phi\land\psi)\eq\top$.
Hence, with $n=p\star q$, $\beta^m\leq\beta^n$ and $A_n(\phi\land\psi)\eq\top$,
so $\beta^m\ll\phi\land\psi$. \qed

\medskip

\begin{r@Remark}
\label{filter/lattice choice}
In summary, the dual basis $A_n\phi$ essentially says that
there is a compact subobject $K^n$ lying between $\beta^n$ and $\phi$,
but $K^n$ seems to play no actual role itself,
and the localic definition in terms of $\ll$ makes it redundant,
\cf Proposition~\ref{open <<}.
Nevertheless, each definition actually has its technical advantages:
\begin{letterlist}
\item in the localic one, $\beta$ ranges over a lattice, but
  $\upup\beta$ need not be a filter;
\item in the spatial one we have filters, but the basis need only be
  indexed by a semilattice.
\end{letterlist}
Effective bases in our sense can be made to behave in \emph{either}
fashion, though we shall only consider lattice bases (\cf the localic
situation) in this paper
(Remark~\ref{lacuna} sketches how filter bases may be obtained).
Stably locally compact objects have lattice filter bases,
whose properties will be improved in Proposition~\ref{ASD dual Wilker}
to take advantage of the intersections of compact subobjects.
\end{r@Remark}

%============================================================================
\section{Sigma-split subobjects}
\label{subobjects}

A basis for a vector space is exactly (the data for) an isomorphism
with $\realno^N$, where $N$ is the dimension of the space.
It is not important for the analogy that the field of scalars be~$\realno$,
or that the dimension be finite.
The significance of $\realno^N$ is that
it carries a \emph{standard structure}
(in which the $n$th basis vector
has a $1$ in the $n$th co-ordinate and $0$ elsewhere),
which is transferred by the isomorphism
to the chosen structure on the space under study.
The standard object in our case is the object $\Sigma^N$
(or the corresponding algebra $\Sigma^{\Sigma^N}$),
for which Axiom~\ref{Scott} defined a basis.

Bases for lattices are actually more like spanning sets than
(linearly independent) bases for vector spaces,
since we may add unions of members to the basis as we please,
as we do in Lemma~\ref{make or basis} below.
Consequently, instead of \emph{isomorphisms} with the standard structure,
we have $\Sigma$-split \emph{embeddings} $X\rSplitinto\Sigma^N$.
We shall see that these embeddings capture several well known
constructions involving $\realno$ and locally compact objects.

\begin{r@Definition}
 $i:X\rSplitinto Y$ is a \textdf{$\Sigma$-split subobject}
if (it is the equaliser of some pair \cite{TaylorP:subasd} and)
there is a map $I:\Sigma^X\to\Sigma^Y$ such that
$\Sigma^i\cdot I=\id_{\Sigma^X}$.
Using Notation~\ref{Hayo notation}, we write
\begin{diagram}
  X& \pile{\rSplitinto^i\\\\\lHto_{\hayo I}} & Y &&
  \Sigma^X& \pile{\lOnto^{\Sigma^i}\\\\\rSplitinto_I} &\Sigma^Y
\end{diagram}
The effect of this is that $X$ carries the subspace topology
inherited from $Y$, in a canonical way:
for an open subobject $\phi$ of $X$,
it provides a particular open subobject $I\phi$ of~$Y$
for which the restriction $\Sigma^i(I\phi)$ is~$\phi$.

In the case where $X$ is an open subobject of $Y$ (Lemma~\ref{subsp quant}),
$I\phi\equiv\exists_i\phi$ is the least such extension,
whilst for a closed subobject, $I\phi\equiv\forall_i\phi$ is the greatest one,
but in general it need not be either of these.

The computational significance of $\Sigma$-split embeddings is that
any observation (computation of type $\Sigma$)
on the subobject extends canonically, though not uniquely, to the whole object.
\end{r@Definition}

\begin{r@Warning}
  Let $Y\rightrightarrows Z$ be a parallel pair of continuous functions
  between locally compact sober spaces
  and $X\equiv\collect{y:Y}{f y=g y:Z}$ the subset of points that satisfy
  the equation.
  We give $X$ the subspace topology,
  \ie its open subspaces are the restrictions of those of~$Y$.
  This is the equaliser $X\splitinto Y\rightrightarrows Z$
  in the category of sober topological spaces,
  but $X$ need not be locally compact (or have an effective basis).
  Equalisers also exist in the category of locales,
  but again need not be locally compact, or even spatial.

  $\Sigma$-split subspaces are therefore a very special case.
\end{r@Warning}

\begin{r@Remark}
For this reason,
the formulation of ASD that we gave in Axiom~\ref{monad axiom}(c)
does not provide general equalisers.
It \emph{formally adjoins} \cite[\S\S 4--6]{TaylorP:subasd}
a $\Sigma$-split subobject $i:X\equiv\collect Y E\splitinto Y$ of $Y$
with $I:\Sigma^X\retract\Sigma^Y$
such that $\Sigma^i\cdot I=\id_{\Sigma^X}$ and $I\cdot\Sigma^i=E$,
given any endomorphism $E:\Sigma^Y\to\Sigma^Y$
that satisfies the appropriate equation.
We shall consider the formulation of this equation
in Section~\ref{prime/nucleus}.
The $\lambda$-calculus developed in \cite[\S\S 8--10]{TaylorP:subasd}
allows arbitrarily complicated combinations of $\Sigma$-split subobjects
and exponentials,
but in this section and the next we reduce them to
$\Sigma$-split subobjects of $\Sigma^N$.
\end{r@Remark}

\medskip

The following results link three (d--f) of the abstract characterisations
of local compactness in the Introduction. %Remark~\ref{summary loc cpct}(d--f).

\begin{r@Lemma}
\label{base->subsp} Any object $X$ that has an effective basis
$(\beta^n,A_n)$ indexed by $N$ is a $\Sigma$-split subobject of
$\Sigma^N$.
\end{r@Lemma}

\Proof Using the basis $(\beta^n,A_n)$, define
\begin{eqnarray*}
  i: X \to \Sigma^N  &\hbox{by}&  x \mapsto \xi\equiv\lambda n. \beta^n x \\
  I: \Sigma^X \to \Sigma^{\Sigma^N}
      &\hbox{by}& \phi\mapsto\Phi\equiv\lambda\xi.\Some n.A_n\phi\land\xi n.
\end{eqnarray*}
Then
$\Sigma^i(I\phi) \;=\; \Lamb x.(I\phi)(i x) \;=\;
  \Lamb x.\Some n. A_n\phi\land\beta^n x \;=\; \phi$.
We also recover
$$ x \;=\; \focus (\Some n.A_n\land\xi n). \eqno\qEd $$

\medskip

Incidentally, notice the use of letters here:
we write $x:X$, $\phi:\Sigma^X$ and $F:\Sigma^{\Sigma^X}$
for the object under study,
but $n:N$, $\xi:\Sigma^N$ and $\Phi:\Sigma^{\Sigma^N}$ for its basis,
the idea being that $x$ and $\phi$ are represented by $\xi$ and $\Phi$
under the embedding.

Conversely,
any $\Sigma$-split subobject inherits the basis of the ambient object,
using the inverse images of the basic open subobjects along~$i:X\rSplitinto Y$.
However, for the compact subobjects, we use their \emph{direct} images
along the ``second class'' map $\hayo I:Y\rHto X$,
in the sense of Remark~\ref{image compact}.
Since $I$ need not preserve meets,
nor need the modal operator $\Sigma^I A\equiv A\cdot I$.
This is why we find bases in which $A_n$ need not preserve $\top$ and $\land$.

\begin{r@Lemma}
\label{subobject basis} Let $(\beta^n,A_n)$ be an effective basis
for $Y$ and $i:X\rSplitinto Y$ a $\Sigma$-split subobject.
Then $(\Sigma^i\beta^n,\Sigma^I A_n)$ is an effective basis for $X$.
If an $\lor$- or $\land$-basis was given, the result is one too.
If $A_n$ is a filter and $I$ preserves $\top$ and $\land$
then $\Sigma^I A_n$ is also a filter.
\end{r@Lemma}

In other words, $\Sigma$-split subobjects of locally compact objects are
again locally compact.

\Proof For $\phi:\Sigma^X$, $I\phi:\Sigma^Y$ has basis decomposition
$$ I\phi \;\eq\; \Some n.A_n(I\phi)\land\beta^n
  \;\equiv\; \Some n.(\Sigma^I A_n)\phi\land\beta^n. $$
Since $\Sigma^i$ is a homomorphism, it preserves scalars, $\land$ and
$\exists$, so
$$ \phi \;=\;
   \Sigma^i(I\phi) \;=\;
   \Sigma^i\big(\Some n.A_n(I\phi)\land\beta^n\big)  \;=\;
   \Some n.A_n(I\phi)\land\Sigma^i\beta^n. \eqno\qEd
$$

\begin{r@Corollary}
\label{basis=<SN}
An object has an effective basis iff it is a $\Sigma$-split
subobject of some $\Sigma^N$.
\end{r@Corollary}

In the extension of ASD, we shall take this as the definition of a
locally compact object.

\Proof The Scott principle (Axiom~\ref{Scott}) defined a basis on
$\Sigma^N$. \qed

\medskip

These two constructions are not inverse:
given an $N$-indexed effective basis on $X$,
we obtain a $\Sigma$-split embedding $X\splitinto\Sigma^N$,
and then a basis indexed by $\Fin(N)$.
In Lemma~\ref{basis from <<} we shall want to recover the original,
$N$-indexed, basis.

\begin{r@Lemma}
\label{embed->basis}
An embedding $X\rSplitinto\Sigma^N$ arises from a basis
according to Lemma~\ref{base->subsp} iff each $I\phi$ preserves joins.
It's then a filter basis iff, for each $n$,
$\Lamb\phi.I\phi\setof n$ also preserves finite meets.
\end{r@Lemma}

\Proof For any basis, $\xi\mapsto I\phi\xi\equiv\Some n. A_n\phi\land\xi n$
preserves joins.
Conversely, with
$$ \beta^n\;\equiv\;\Lamb x.i x n
    \quad\hbox{and}\quad
    A_n\;\equiv\;\Lamb\phi.I\phi\setof n
$$
we recover
$ \phi x \;\equiv\; (I\phi)(i x) \;\eq\; \Some n.I\phi\setof n\land i x n
   \;\eq\; \Some n. A_n\phi\land\beta^n x
$
so long as $I\phi$ preserves the join $i x=\Some n.i x n\land\setof n$. \qed

\bigbreak

In the rest of this section we consider the classical interpretations
of the $\Sigma$-split embedding that arises from an effective basis,
starting with traditional topology.
Recall from \cite[Theorem II 1.2]{JohnstonePT:stos} that the free frame
on $N$ is $\Upsilon\Kur N$ (the lattice of upper subsets of $\Kur N$),
and that this is isomorphic to the lattice of Scott-open subsets of
the powerset $\powerset(N)$.

\begin{r@Theorem}
\label{lcpct<PN} Let $X$ be a locally compact sober space with
$N$-indexed basis $(U^n,K^n)$. Then $X$ is a $\Sigma$-split subspace
of $\powerset(N)$.
\end{r@Theorem}

\Proof The embedding in Lemma~\ref{base->subsp} takes
\begin{eqnarray*}
   x\in X &\hbox{ to }& \collect n{x\in U^n}\in\powerset(N) \\
   V\subset X  &\hbox{ to }&
   \collect{\ell}{\Some n\in\ell.K^n\subset V}\in\Upsilon\Kur N.
\end{eqnarray*}
The second map, $I$, is Scott-continuous because it takes $\dirunion_s V_s$ to
$$ \collect{\ell}{\Some n\in\ell.K^n\subset\dirunion_s V_s}
  \;=\; \collect{\ell}{\Some n\in\ell.\Some s.K^n\subset V_s}
  \;=\; \dirunion_s\collect{\ell}{\Some n\in\ell.K^n\subset V_s}. $$
The composite $\Sigma^i\cdot I$ takes $V\subset X$ to
$ \bigcup \collect{\bigcap_{n\in\ell}U^n}{\Some n\in\ell.K^n\subset V}$,
and by Definition~\ref{pt loc cpct},
this contains $V=\collect x{\Some n.x\in U^n\subset K^n\subset V}$.
For the converse, if $x\in\Sigma^i(I V)$ then
$\Some\ell.\All n\in\ell. x\in U^n\land \Some m\in\ell.K^n\subset V$,
so $\Some n.x\in U^n\subset K^n\subset V$. \qed

\smallskip

\begin{r@Example}
 A compact Hausdorff space has a basis determined by a family
of disjoint pairs $(U^n\disjoint V_n)$ of open subspaces.
In this case, the embedding is
\begin{eqnarray*}
  x &\mapsto &\collect n{x\in U^n}  \\
  W &\mapsto &\collect\ell{\Some n\in\ell.V_n\cup W=X} &\qeds
\end{eqnarray*}
\end{r@Example}

\bigskip

Consider in particular the embedding of $\realno$ in $\Sigma^N$,
where $N$ indexes a basis of open and closed intervals
(Examples \ref{eg R1} and~\ref{eg R2}).
This is closely related to one of the first examples that Dana Scott
used to show how continuous lattices could be used as a model of
computation \cite{ScottDS:outmtc}.

\begin{r@Definition}
\label{intdom}
 The \textdf{lattice of intervals}, $\domint\realno^\top$,
of $\realno$ is the set of convex closed subspaces
(including the empty one),
ordered by \emph{reverse} inclusion, and given the Scott topology.
The \textdf{domain of intervals}, $\domint\realno\subset\domint\realno^\top$,
consists of the inhabited such subspaces.
Amongst these, points of $\realno$ are identified with
the \emph{maximal} elements.
Classically, the members of $\domint\realno$ are of the form
$[r\pm\delta]$, with $r\in\realno$ and $0\leq\delta\leq\infty$.
\end{r@Definition}

\begin{r@Remark}
 In ASD, closed subobjects under reverse inclusion 
correspond to their co-classifiers under the (forward) intrinsic order,
and to the ``complementary'' open subobjects under inclusion.
When the closed subobject is convex, the open one is $\delta\lor\upsilon$,
where $\delta$ is lower and $\upsilon$ upper in the arithmetical order,
so the pair $(\delta,\upsilon)$ is almost a Dedekind cut~\cite{TaylorP:dedras},
except that it is disjoint but need not be ``located''.
It has this last property exactly when the closed subobject is a singleton.

Intervals $[p,q]$ or $[r\pm\delta]$ with specified endpoints are of this form,
but (in ASD and other constructive forms of analysis),
not every inhabited convex closed bounded subobject need have endpoints.
It does so iff it is also overt \cite{TaylorP:lamcra}.
In particular, a directed join in $\domint\realno$
of subobjects with endpoints need not itself have them.

It is, however, not necessary to appreciate this issue
in order to see the relationship between $\domint\realno$
and our representation of $\realno$ as a $\Sigma$-split subobject of $\Sigma^N$.
\end{r@Remark}

\begin{r@Proposition}
\label{R-SN} In $\LKSp$,
$\realno \rSplitinto \domint\realno \rClosedinto \domint\realno^\top 
\pile{\rSplitinto\\\lOnto}\Sigma^N$,
where $N$ is the set of pairs, written $<q\pm\epsilon>$, with
$\epsilon\gt0$ and $q$ rational.
\end{r@Proposition}

(Recall that the $\rClosedinto$ indicates a \emph{closed} inclusion.)

\Proof The embedding takes $r\in\realno$ to $[r\pm0]$ and then to
$\Lamb<q\pm\epsilon>.r\in(q\pm\epsilon)$,
which is (the exponential transpose of) a continuous function.
The retraction is defined by intersection:
\begin{eqnarray*}
  \intdom\realno^\top &\pile{\rSplitinto\\\lOnto}& \Sigma^N \\
  C &\rMapsto&\Lamb<q\pm\epsilon>.C\subset(q\pm\epsilon) \\\relax
  [r\pm\delta] &\rMapsto&
     \Lamb<q\pm\epsilon>.q-\epsilon\lt r-\delta\leq r+\delta\lt q+\epsilon \\
  \emptyset &\rMapsto&\Lamb<q\pm\epsilon>.\top \\
  \realno &\rMapsto&\Lamb<q\pm\epsilon>.(\epsilon=_N\infty) \\
  \!\!\!\!\bigcap_{\phi<q\pm\epsilon>}[q\pm\epsilon] &\lMapsto& \phi
\end{eqnarray*}
Any compact interval of the ${\mathrm T}_1$ space $\realno$ is saturated
in the sense of Lemma~\ref{cpct sat}.
It is therefore the intersection of its open neighbourhoods,
amongst which open intervals suffice.
Hence the composite is $\domint\realno^\top\to\Sigma^N\to\domint\realno^\top$
is the identity.

The projection $\domint\realno^\top\lOnto\Sigma^N$ is Scott-continuous
because it clearly takes directed unions of sets of codes to
codirected intersections of compact subspaces.
Classically, Lemma~\ref{codir meet cpct} showed that such intersections
correspond to unions of neighbourhood filters,
whilst in ASD they correspond to the unions
of the complementary open subobjects $\delta\lor\upsilon$.
Hence the inclusion $\domint\realno^\top\rSplitinto\Sigma^N$
is also Scott-continuous.

The inverse image of $\top$ under $\domint\realno^\top\lOnto\Sigma^N$
is the open subobject classified by the \emph{inconsistency predicate}
$$ \InCon(\phi) \;\equiv\;
  \Some<q_1\pm\epsilon_1><q_2\pm\epsilon_2>.
  (q_1+\epsilon_1\lt q_2-\epsilon_2)
  \land\phi<q_1\pm\epsilon_1>\land\phi<q_2\pm\epsilon_2>.
$$
The complementary closed subobject of $\Sigma^N$ is of course not classified,
as it's not open,
but when we restrict attention to its (overt discrete collection of)
finite elements,
we find that \emph{consistency} is characterised by the decidable formula
$$ \Con(\ell) \;\equiv\;
  \Some x:\ratno.\All<q\pm\epsilon>\in\ell.x\in<q\pm\epsilon>,
$$
so\closeupaline
$$ \InCon(\Lamb<q\pm\epsilon>.<q\pm\epsilon>\in\ell) \;=\; \lnot\Con(\ell).
   \eqno\qEd
$$

\begin{r@Remark}
\label{discuss intdom}
The idea behind the domain of intervals is not hard to generalise.
Indeed, we may embed any locally compact sober space as a subspace of
its continuous preframe of compact saturated subspaces (Theorem~\ref{cpct <<}),
each point being represented by its saturation in
the sense of Lemma~\ref{cpct sat}.
That the image consists of the \emph{maximal} points (excluding $\emptyset$)
plainly depends on starting with a ${\mathrm T}_1$-space,
so can't be an essential feature of the construction.

Another interesting aspect of the general construction is
the decidable consistency predicate on finite elements.
Scott developed these into what became a standard form of domain theory
\cite{ScottDS:domds}.
The subobject of $\Sigma^N$ that is determined by the consistent elements
is closed, but also overt \cite{TaylorP:pcfasd}.

We shall develop our own construction of the real line \via Dedekind cuts,
and also discuss the interval domain further, in \cite{TaylorP:dedras}.
\end{r@Remark}

\bigbreak

Now we give the localic version of Theorem~\ref{lcpct<PN}.
This result is the most efficient way of establishing the connection 
between locally compact locales and ASD.

\begin{r@Theorem}
\label{ctsdist<UKN}
Let $L$ be any $N$-based continuous distributive lattice.
Then there is a frame homomorphism $H$ and a Scott-continuous function $I$
with $H\cdot I=\id_L$, as shown:
\begin{diagram}[midshaft]
  L & \pile{\lOnto^H\\\\\rSplitinto_I} & \Upsilon\Kur N
\end{diagram}
Conversely, any lattice $L$ that admits such a pair of functions
is continuous and distributive.
\end{r@Theorem}

\Proof $[\Implies]$ Let $(\beta^n)$ be a basis for the continuous
lattice $L$ and $A_n\equiv\Lamb\phi.(\beta^n\ll\phi)$.
Let $H:\Upsilon\Kur N\to L$ be the unique frame homomorphism that
extends $\beta^\blank:N\to L$, so
$$ \hbox{for } \U\in\Upsilon\Kur N, \quad
     H\U \;=\; \bigvee_{\ell\in\U}\bigwedge_{n\in\ell}\beta^n \in L, $$
and define $I:L\to\Upsilon\Kur N$ by
$ I\phi\;=\;\collect{\ell}{\Some n\in\ell.\beta^n\ll\phi}$.
This is Scott-continuous by a similar argument to that in
Theorem~\ref{lcpct<PN}, with $\beta^n\waybelow\phi_s$ instead of
$K^n\subset V_s$. Then
\begin{eqnarray*}
  \beta^n\ll H(I\phi)
  &\iff& \Some\ell.\beta^n\ll\bigwedge_{m\in\ell}\beta^m \ \land\
            \Some m\in\ell.\beta^m\ll\phi\\
  &\iff& \Some m.\beta^n\ll\beta^m\ll\phi
  \;\iff\; \beta^n\ll\phi,
\end{eqnarray*}
from which we deduce $H(I\phi)=\phi$, because $(\beta^n)$ is a basis.

$[\Impliedby]$ Conversely, if such a diagram exists then $I\cdot H$ is a
Scott-continuous idempotent on a continuous lattice,
so its splitting $L$ is also continuous \cite[Lemma VII 2.3]{JohnstonePT:stos}.
As $H$ preserves joins, it has a right adjoint, $H\adjoint R$,
so $\id\leq R\cdot H\equiv j=j\cdot j$,
and $R$ preserves meets but not necessarily directed joins.
Since $H$ also preserves finite meets, so does $j$,
and this is a nucleus in the sense of locale theory
\cite[Section II 2.2]{JohnstonePT:stos}, so its splitting $L$ is a frame.\qed

\medskip

Any locally compact locale is therefore determined by a
Scott-continuous idempotent $\E$ on $\Upsilon\Kur N$.
It is not just an idempotent, however,
since the surjective part of its splitting must be a frame homomorphism.
Since the latter preserves $\top$ and $\bot$ by monotonicity,
and $\dirsup$ as $\E$ is Scott-continuous,
it is enough to identify the condition on $\E$ the ensures
preservation of the two binary lattice connectives,
which we may treat exactly alike.
(This is where Remark~\ref{formulate prime/nucleus} came from.)

\begin{r@Lemma}
\label{E@} Let $I$ and $H$ be monotone functions between two
semilattices, with $H\cdot I=\id$. Then $H$ preserves $\land$ iff
$E\equiv I\cdot H$ satisfies the equation $E(\phi\land\psi) =
E(E\phi\land E\psi)$.
\end{r@Lemma}

\Proof If $H$ preserves $\land$ then
\begin{eqnarray*}
  E(\phi\land\psi) &\equiv& I\big(H(\phi\land\psi)\big)  \\
  &=& I(H\phi\land H\psi) &hypothesis\\
  &=& I(H\cdot I\cdot H\phi\land H\cdot I\cdot H\psi) &$H\cdot I=\id$\\
  &=& I\cdot H(I\cdot H\phi\land I\cdot H\psi) &hypothesis\\
  &\equiv& E(E\phi\land E\psi) 
\end{eqnarray*}
For the converse, note first that we have
$H(\phi\land\psi)\leq H\phi\land H\psi$ and
$I(\phi'\land\psi')\leq I\phi'\land I\psi'$
by the definition of $\land$. Then
\begin{eqnarray*}
  I(H\phi\land H\psi)
  &=& E\big(I(H\phi\land H\psi)\big) &$I=I\cdot H\cdot I=E\cdot I$\\
  &\leq& E(I\cdot H\phi\land I\cdot H\psi) &above\\
  &=& E(E\phi\land E\psi) &$E=I\cdot H$\\
  &=& E(\phi\land \psi) &hypothesis\\
  &=& I\cdot H(\phi\land \psi) &$E=I\cdot H$\\
  H\phi\land H\psi &\leq& H(\phi\land \psi) &$H\cdot I=\id$%
\end{eqnarray*}
so $H(\phi\land\psi)\leq H\phi\land H\psi\leq H(\phi\land\psi)$. \qed

\begin{r@Corollary}
 Any $N$-based locally compact locale is determined by
a Scott-continuous endofunction of the frame $\Upsilon\Kur N$
(or by a locale endomorphism of $\Sigma^{\Sigma^N}$) such that
$$ E(\phi\land\psi) = E(E\phi\land E\psi)
   \quad\hbox{and}\quad
   E(\phi\lor\psi) = E(E\phi\lor E\psi).
   \eqno\qEd
$$
\end{r@Corollary}

We now have to concentrate on the logical development within abstract
Stone duality, and will only return to the connection with traditional
topology in Section~\ref{relate}.
The first task is to show that every definable object has an effective
basis, and is therefore a $\Sigma$-split subobject of $\Sigma^N$.
Such subobjects are determined by idempotents $\E$ on $\Sigma^{\Sigma^N}$
satisfying the equation that we have just identified, 
along with its counterpart for~$\lor$.
However, even that characterisation depends on the use of bases
(Lemma~\ref{nucleus->E@}ff).

%============================================================================
\section{Every definable object has a basis}\label{bases}

In Section~\ref{lcpct}, we justified the notion of effective basis in
the classical models, \ie for locally compact sober spaces and locales.
This section considers the term model, showing by structural recursion
that every definable object has an effective basis.
We have already dealt with the base cases ($\natno$, $\Sigma^\natno$)
in Example~\ref{N+SN},
and with $\Sigma$-split subobjects in Lemma~\ref{subobject basis}.
So we consider binary products first
(leaving the reader to define the $\terminalobj$-indexed basis
for~$\terminalobj$),
but devote most of the section to the exponential $\Sigma^X$.

\begin{r@Lemma}
\label{prod basis} If $X$ and $Y$ have effective bases then so
does $X\times Y$, given by Tychonov rectangles.
\end{r@Lemma}

\Proof We are given $(\beta^n,A_n)$ and $(\gamma^m,D_m)$ on $X$ and $Y$.
Then, on $X\times Y$, we define
\begin{eqnarray*}
  \epsilon^{(n,m)} &\equiv&\Lamb x y.\beta^n x\land\gamma^m y \\
  F_{(n,m)} &\equiv&
  \Lamb\theta:\Sigma^{X\times Y}.D_m\big(\Lamb y.A_n(\Lamb x.\theta(x,y))\big),
  \\
  \llap{Then\qquad\qquad}
  \theta(x,y)
  &\eq& \Some n.A_n\big(\Lamb x'.\theta(x',y)\big)\land\beta^n x \\
  &\eq& \Some n m.
     D_m\big(\Lamb y'.A_n(\Lamb x'.\theta(x',y'))\big)
    \land\gamma^m y\land\beta^n x \\
  &\eq& \Some (n,m). F_{(n,m)}\theta\land\epsilon^{(n,m)}(x,y) &$\qeds$%
\end{eqnarray*}

\smallskip

Notice that the formula for $F_{(n,m)}$ is not symmetrical in $X$ and
$Y$, though we have learned to expect properties of binary products to
be asymmetrical and problematic \cite{TaylorP:sobsc}.
In fact, if $A_n$ and $D_m$ are filter bases, we have another example
of the same problem that held us back in Proposition~\ref{cpct=closed},
along with the core of its solution.

\begin{r@Lemma}
\label{finite AD=DA}
If $A:\Sigma^{\Sigma^X}$ and $D:\Sigma^{\Sigma^Y}$ both preserve
\emph{either} $\top$ and $\land$ \emph{or} $\bot$ and $\lor$ then
$$ A\big(\Lamb x.D(\Lamb y.\theta x y)\big) \;\eq\;
   D\big(\Lamb y.A(\Lamb x.\theta x y)\big) $$
whenever $\theta:\Sigma^{X\times Y}$ is a finite union of rectangles.
\end{r@Lemma}

\Proof Applying the Phoa principle (Axiom~\ref{Phoa}) to a single rectangle,
\begin{eqnarray*}
  A\big(\Lamb x.D(\Lamb y.\phi x\land\psi y)\big)
     \hidewidth\\\qquad
  &\eq& A(\Lamb x.D\bot\lor\phi x\land D\psi) &Phoa for $D$ wrt $\phi x$\\
  &\eq& A(\Lamb x.\phi x\land D\psi)\lor (D\bot\land A\top)
      &Phoa for $A$ wrt $D\bot$\\
  &\eq& A\bot\lor(D\psi\land A\phi)\lor(D\bot\land A\top)
      &Phoa for $A$ wrt $D\psi$.%
\end{eqnarray*}
This would have the required symmetry if we had
$$ A(D\bot) \;\equiv\; A\bot\lor(D\bot\land A\top)
   \;\eq\; (A\bot\land D\top)\lor D\bot  \;\equiv\; D(A\bot).
$$
If $A\top\eq\top\eq D\top$ then both sides are $A\bot\lor D\bot$,
whilst if $A\bot\eq\bot\eq D\bot$ then they are both~$\bot$.
Under either hypothesis, the lattice dual argument shows the similar result
$$ A\big(\Lamb x.D(\Lamb y.\phi x\lor\psi y)\big) 
  \;\eq\; D\big(\Lamb y.A(\Lamb x.\phi x\lor\psi y)\big) $$
for a cross.
The result uses equational induction \cite[\S 2]{TaylorP:insema},
but we shall just illustrate this with the case where
$\theta$ is a union of three rectangles,
$$ \theta x y
  \;\eq\; (\phi_1 x\land\psi_1 y) \,\lor\, (\phi_2 x\land\psi_2 y) \,\lor\,
  (\phi_3 x\land\psi_3 y). $$
If $A$ and $D$ preserve $\bot$ and $\lor$ then $A D\theta$ is also a union
of three terms, to each of which the first part applies.

We may also use distributivity of $\lor$ over $\land$
to re-express the union $\theta$ as an intersection of eight crosses,
\begin{eqnarray*}
  \theta x y
  &\eq& (\phi_1 x\lor\phi_2 x\lor\phi_3 x) \,\land\,
      (\phi_1 x\lor\phi_2 x\lor\psi_3 y) \,\land\,\\
  &&  (\phi_1 x\lor\psi_2 y\lor\phi_3 x) \,\land\,
      (\phi_1 x\lor\psi_2 y\lor\psi_3 y) \,\land\,\\
  &&  (\psi_1 y\lor\phi_2 x\lor\phi_3 x) \,\land\,
      (\psi_1 y\lor\phi_2 x\lor\psi_3 y) \,\land\,\\
  &&  (\psi_1 y\lor\psi_2 y\lor\phi_3 x) \,\land\,
      (\psi_1 y\lor\psi_2 y\lor\psi_3 y).
\end{eqnarray*}
So if $A$ and $D$ preserve $\land$ then $A D\theta$ is a conjunction of
eight factors, to each of which the second part applies,
for example with $\phi=\phi_1\lor\phi_2$ and $\psi=\psi_3$.
In both cases, $A D\theta\eq D A\theta$ as required. \qed

\begin{r@Remark}
\label{quantifiers commute} This still awaits
Theorem~\ref{F:SSX cts} on Scott continuity to extend finite unions of
rectangles to infinite ones, but once we have that we may draw the
corollaries that
\begin{letterlist}
\item if $X$ and $Y$ have filter bases then Lemma~\ref{prod basis}
  provides a filter basis for $X\times Y$, which is symmetrical in $X$ and $Y$;
\item Proposition~\ref{cpct=closed} yields a bijection between closed
  and compact subobjects of a compact Hausdorff object.
  Moreover, despite the other problems discussed in Section~\ref{cpct subsp},
  preserving \emph{finite} meets is enough to characterise $\square$
  modal operators.  \qed
\end{letterlist}
\end{r@Remark}

\bigskip

We shall need to be able to turn any effective basis into a
$\lor$-basis, which we do in the obvious way using finite unions of
basic open subobjects.
The corresponding \emph{unions} of compact subobjects give rise to
\emph{conjunctions} of $A$s by Lemma~\ref{cup cpct}.
Unfortunately, the result is topologically rather messy,
both for products and
for objects such as $\realno$ that we want to construct directly.

\begin{r@Lemma}
\label{make or basis} If $X$ has an effective basis indexed by
$N$ then it also has a $\lor$-basis indexed by $\Fin(N)$.  If we were
given a filter basis, the result is one too.
\end{r@Lemma}

\Proof Given any basis $(\beta^n,A_n)$, define
$$ \gamma^\ell \;\equiv\; \Lamb x.\Some n\in\ell.\beta^n x \qquad
    D_\ell \;\equiv\; \Lamb\phi.\All n\in\ell.A_n\phi. $$
Then
$\phi\;=\;\Some n.A_n\phi\land\beta^n\;\leq\;
\Some\ell.D_\ell\phi\land\gamma^\ell$
using singleton lists.
Conversely, 
\begin{eqnarray*}
  \Some\ell.D_\ell\phi\land\gamma^\ell
  &=& \Some\ell.(\All n\in\ell.A_n\phi)\land(\Some m\in\ell.\beta^m) \\
  &=& \Some\ell.\Some m\in\ell.(\All n\in\ell.A_n\phi)\land\beta^m \\
  &\leq& \Some\ell.\Some m\in\ell.A_m\phi\land\beta^m \;=\; \phi.
\end{eqnarray*}  
Then $(\gamma^\ell,D_\ell)$ is a $\lor$-basis,
using list concatenation for~$+$.
The imposed order $\baseleq$ on $\Fin(N)$ is list or subset inclusion,
$$ (\ell\baseleq\ell') \;\equiv\; (\ell\subset\ell')
    \;\equiv\; \All n\in\ell.n\in\ell'. $$
Finally, if the $A_n$ were filters then so are the $D_\ell$,
since $\forall m\in\ell$ preserves $\land$ and~$\top$. \qed

\begin{r@Remark}
\label{embed in SSN}
Computationally, it may be better to see this as an embedding
$i:X\splitinto\Sigma^{\Sigma^N}$ rather than into $\Sigma^{\Fin N}$,
as this gives a very simple representation of the basis and dual basis.
Using the notation $\pi_\ell\equiv\Lamb n.n\in\ell$ and
$[\ell]\equiv\All n\in\ell.\psi n$ from \cite{TaylorP:insema},
which also provides a fixed-point equation for $\exists\ell$, this is:
\begin{eqnarray*}
  i x &=& \Lamb\xi.\Some n.\beta^n x \land \xi n \\
  I \phi &=&
  \Lamb\Phi.\Some\ell.\Phi(\Lamb n.n\in\ell) \land \All n\in\ell.A_n\phi\\
  \gamma^\ell x &\eq& 
  \Some n\in\ell.\beta^n x \;\eq\; i x (\Lamb n.n\in\ell) \;\eq\; i x\pi_\ell \\
  D_\ell\phi &\eq&
  \All n\in\ell.A_n\phi \;\eq\; I\phi (\Lamb\psi.\All n\in\ell.\psi n)
  \;\eq\; I\phi[\ell] \\
  \phi x &\eq& \Some\ell.I\phi[\ell]\land i x \pi_\ell &$\qeds$%
\end{eqnarray*}
\end{r@Remark}

\begin{r@Lemma}
\label{make lattice basis} If an $\land$-basis was given,
Lemma~\ref{make or basis} yields a lattice basis.
\end{r@Lemma}

\Proof We are given $\beta^1=\top$ and $\beta^n\land\beta^m=\beta^{n\star m}$.

Let $\ell\star\ell'$ be the list (it doesn't matter in what order)
of all $n\star m$ for $n\in\ell$ and $m\in\ell'$.
In functional programming notation, this is
$$ \ell\star\ell' \;\equiv\; \FLATTEN\;
   \big(\MAP\,\ell\,(\Lamb n.\MAP\,\ell'(\Lamb m.n\star m))\big), $$
where $\MAP$ applies a function to each member of a list, returning a list,
and $\FLATTEN$ turns a list of lists into a simple list.
Categorically, $\MAP$ is the effect of the functor $\List\blank$ on morphisms,
and $\FLATTEN$ is the multiplication for the $\List$ monad.
Using the corresponding notions for $\Kur\blank$,
$\star$ can similarly be defined for finite subsets instead.

Now let
$$ \beta^{\ell\star\ell'} \;\equiv\;
    (\Some n\in\ell.\Some m\in\ell'.\beta^{n\star m}) \;=\;
   (\gamma^\ell\land\gamma^{\ell'})
$$
by distributivity in $\Sigma^X$,
whilst $\gamma^{\listof 1}\equiv\beta^1$ and $D_{\listof 1}\equiv A_1$
(where $\listof 1$ is the singleton \emph{list})
serve for $\gamma^1$ and~$D_1$.
(This uses equational induction, \cite[\S 2]{TaylorP:insema}.)
 \qed

\smallskip

\begin{r@Remark}
\label{make and basis}
By switching the quantifiers, we may similarly obtain $\land$-bases,
and turn an $\lor$-basis into a lattice basis.
In fact, this construction featured in the proof of Theorems \ref{lcpct<PN}
and~\ref{ctsdist<UKN}.
This time the filter property is not preserved,
since if $A_n\phi$ means $K^n\subset U$ then $A_\ell\phi$
means that $\Some n\in\ell.K^n\subset U$,
which is not the same as testing containment of a single compact subobject.
Proposition~\ref{ASD dual Wilker} shows how to define an $\land$-basis
for a locally compact object in which $A_{n\star m}$ actually captures
the intersection of the compact subobjects.
\end{r@Remark}

\medskip

\begin{r@Proposition}
\label{N bases}
$N$ and $\Sigma^N$ have effective bases as follows:
$$\def\arraystretch{1.2}
\begin{array}{llll}
  % \hbox{on} && \beta^\blank && A_\blank \\
  N & \hbox{prime}
  & \beta^n\equiv\setof n\equiv\Lamb m.(n=m)
  & A_n\equiv\eta_N(n)\equiv\Lamb\xi.\xi n \\
  N & \hbox{filter }\lor
  & \beta^\ell\equiv\Lamb m.m\in\ell 
  & A_\ell\equiv\Lamb\xi.\All m\in\ell.\xi m \\
  N & \hbox{lattice}
  & \beta^L\equiv\Lamb m.\All\ell\in L.m\in\ell
  & A_L\equiv\Lamb\xi.\Some\ell\in L.\All m\in\ell.\xi m \\
  \Sigma^N & \hbox{prime }\land
  & B^\ell \equiv A_\ell 
  & \A_\ell \equiv \eta_{\Sigma^N}\beta^\ell \equiv
        \Lamb \Phi.\Phi(\Lamb m.m\in\ell) \\
  \Sigma^N & \hbox{filter lattice}
  & B^L \equiv A_L 
  & \A_L \equiv \Lamb \Phi.\All\ell\in L.\Phi(\Lamb m.m\in\ell)
\end{array}
$$
indexed by $n:N$, $\ell:\Fin(N)$ or $L:\Fin\big(\Fin(N)\big)$.
The interchange of sub- and superscripts means that we're reversing
the imposed order on these indexing objects (Remark~\ref{dep bigvee}).
\qed
  \end{r@Proposition}

\medbreak

The last of these also provides a basis for $\Sigma^{\Sigma^N}$,
using Axiom~\ref{Scott} twice.

\nobreak

\begin{r@Lemma}
\label{SSN basis} $\Sigma^{\Sigma^N}$ has a
$\Fin(\Fin N)$-indexed prime $\land$-basis with
$$ \B^L\equiv\A_L\equiv\Lamb \Phi.\All\ell\in L.\Phi\beta^\ell
   \quad\hbox{and}\quad
   \AAA_L \equiv \eta_{\Sigma^2 N}A_L \equiv\Lamb\F.\F A_L,
$$
so, using $\Sigma^3 N$ as a shorthand for a tower of exponentials,
$\Sigma^{\sizenine{\Sigma^{\Sigma^N}}}\!\equiv
  \big(((N\to\Sigma)\to\Sigma)\to\Sigma\big)$,
$$ \F:\Sigma^3 N,\; \Phi:\Sigma^{\Sigma^N} \ \proves\
   \F \Phi \;\eq\;
   \Some L:\Fin(\Fin N).\; \F A_L \land \All\ell\in L. \Phi\beta^\ell. $$
\end{r@Lemma}

\Proof The prime $\land$-basis on $\Sigma^N$
makes $\Sigma^{\Sigma^N}\retract\Sigma^{\Fin(N)}$ by
$$ \Phi\mapsto\Lamb\ell.\Phi\beta^\ell 
   \quad\hbox{and}\quad
   \Lamb\xi.\Some\ell.\Psi\ell\land\All n\in\ell.\xi n\mapsfrom \Psi, $$
so $\Sigma^3 N\retract\Sigma^2\Fin(N)$ by
$$ \F\mapsto\Lamb \Psi.\F(\Lamb\xi.\Some\ell.\Psi\ell\land\All n\in\ell.\xi n)
   \quad\hbox{and}\quad
   \Lamb \Phi.\G(\Lamb\ell.\Phi\beta^\ell)\mapsfrom\G. $$
For any $\F:\Sigma^3 N$ and $\Psi:\Sigma^{\Fin(N)}$,
let $\G$ be obtained from $\F$ in this way,
and use the prime $\land$-basis on $\Sigma^{\Fin(N)}$:
\begin{eqnarray*}
  \G \Psi  &\eq&
 \Some L:\Fin\Fin(N).\G(\Lamb\ell.\ell\in L)\land\All\ell\in L.\Psi\ell \\
  &\eq& \Some L.
    \F\big(\Lamb\xi.\Some\ell.
     (\Lamb\ell.\ell\in L)\ell\land\All n\in\ell.\xi n\big)
     \land  \All\ell\in L.\Psi\ell \\
  &\eq& \Some L. \F(\Lamb\xi.\Some\ell\in L.\All n\in\ell.\xi n) 
     \land  \All\ell\in L.\Psi\ell \\
  &\eq& \Some L. \F A_L \land  \All\ell\in L.\Psi\ell 
\end{eqnarray*}
so $\F \;=\; \Lamb \Phi.\G(\Lamb\ell.\Phi\beta^\ell) \;=\;
   \Lamb \Phi.\Some L. \F A_L \land  \All\ell\in L.\Phi\beta^\ell$.
\qed

\begin{r@Lemma}
\label{SX lbasis}
If an object $X$ has an effective basis $(\beta^n,A_n)$ then
its topology $\Sigma^X$ has a prime lattice basis $(D^L,\Lamb F.F\gamma_L)$,
indexed by $L:\Fin(\Fin N)$, where
$$ \gamma_L\;\equiv\;\Lamb x.\Some\ell\in L.\All n\in\ell.\beta^n x
   \quad\hbox{and}\quad
   D^L\;\equiv\;\Lamb\phi.\All\ell\in L.\Some n\in\ell.A_n\phi.
$$
\end{r@Lemma}

\Proof By Lemma~\ref{base->subsp}, there is an embedding
$i:X\rSplitinto\Sigma^N$ with $\Sigma$-splitting $I$ by
$$ i x \;\equiv\; \Lamb n.\beta^n x
   \quad\hbox{and}\quad
   I\phi\;\equiv\; \Phi \;\equiv\; \Lamb\xi.\Some n.A_n\phi\land\xi n,
   \quad\hbox{so}\quad
   \phi \;=\; \Lamb x.\Phi(i x).
$$
For $F:\Sigma^{\Sigma^X}$, let
$\F \equiv F\cdot\Sigma^i \equiv \Lamb\Psi.F\big(\Lamb x.\Psi(i x)\big)
 :\Sigma^3 N$,
so $\F\cdot I=F\cdot\Sigma^i\cdot I=F$.
Then
\begin{eqnarray*}
  F\phi \;\eq\; \F(I\phi)
  &\eq& \Some L. \F A_L \land  \All\ell\in L.I\phi\beta^\ell 
  &Lemma~\ref{SSN basis}\\
  &\eq& \Some L. F(\Sigma^i A_L) \land  \All\ell\in L.I\phi\beta^\ell \\
  &\eq& \Some L. F\gamma_L \land  D^L\phi
\end{eqnarray*}
\vadjust{\nobreak}%
since the given formulae are $\gamma_L=\Sigma^i A_L$ and 
$D^L=\Lamb\phi.\All\ell\in L.I\phi\beta^\ell$. \qed

\goodbreak

\begin{r@Theorem}
\label{all loc cpct}
Every definable (or locally compact, \cf Remark~\ref{ext def loc cpct})
object $X$ has a lattice basis, indexed by $\Fin(\Fin\natno)$,
and is a $\Sigma$-split subobject of $\Sigma^\natno$.
Moreover, $\Sigma^X$ is \emph{stably} locally compact,
as it has a prime lattice basis (\cf Proposition~\ref{slc=latt filt}). \qed
\end{r@Theorem}

\begin{r@Remark}
\label{nfthm} This is a ``normal form'' theorem, and, like all such
theorems, it can be misinterpreted.
It is a bridge over which we may pass in \emph{either} direction
between $\lambda$-calculus and a discrete encoding of topology,
not an intention to give up the very pleasant synthetic results
that we saw in \cite{TaylorP:geohol}.
In particular, we make no suggestion that either arguments in topology
or their computational interpretations need go \emph{via}
the list or subset representation
(though \cite{JungA:mullsc} seems to have this in mind).
Indeed, subsets may instead be represented by $\lambda$-terms
\cite{TaylorP:insema}.
It is simply a method of proof, and is exactly what we need
to connect synthetic abstract Stone duality
with the older lattice-theoretic approaches to topology,
as we shall now show.
\end{r@Remark}

%============================================================================
\section{Basic corollaries}\label{corolls}

Making use of the availability of bases for all definable objects, this
section establishes the basic properties that justify the claim that
abstract Stone duality is an account of domain theory and general
topology, at least in so far as its morphisms are continuous functions.
We first prove something that was claimed in Remark~\ref{N enough}.

\begin{r@Proposition}
\label{ov disc=subquot N}
Let $M$ be an overt discrete object with effective basis
$(\beta^n,A_n)$ indexed by an overt discrete object $N$.
Then $M$ is the subquotient of $N$ by an open partial equivalence relation.
\end{r@Proposition}

\Proof Write $n\realises x$ for
$n:N,\; x:M \proves A_n\setof x\land \beta^n x$ (using discreteness of $M$)
and $N'=\collect n{\Some x.n\realises x}\subset N$, which is open,
using overtness.

Then, using the basis expansion of open $\setof x$,
$$ % x:M\ \proves\
   \top\;\eq\;(x=_M x)
   \;\eq\;\setof x(x)
   \;\eq\;\Some n.A_n\setof x\land\beta^n x
   \;\eq\;\Some n.n\realises x,
$$
so every point $x:M$ has some code $n:N'$.
The latter belongs only to $x$ since
\begin{eqnarray*}
  n\realises x\land n\realises y 
  &\eq& A_n\setof x\land\beta^n x\land A_n\setof y\land\beta^n y \\
  &\Implies& A_n\setof x\land \beta^n y  \\
  &\equiv& (\Some n.A_n\setof x\land\beta^n) y \\
  &\eq& \setof x y \;\eq\; (x=_M y).
\end{eqnarray*}
Hence $N'\to\Sigma^M$ by $n\mapsto\Lamb x.{n\realises x}$ factors through
$\setof{}:M\rSplitinto\Sigma^M$, and $M$ is $N/{\sim}$ where
$(m\sim n)\equiv(\Some x.m\realises x\land n\realises x)$
\cite{TaylorP:geohol}. \qed

\medskip

We leave it as an exercise to define partial equivalence relations on
$\List(N)$ whose quotients are $\List(M)$ and $\Kur(M)$.

\begin{r@Corollary}
\label{ovdisc=subquot} Every definable overt discrete object
is a subquotient of $\natno$ by an open partial equivalence relation. \qed
\end{r@Corollary}

\smallskip

This does not restrict how ``big'' overt discrete objects can be
in \emph{general} models of ASD:
for example, $\aleph_1$ still belongs to the classical model.
It simply says that, having required certain \emph{base} types
to be overt discrete, as we did with $\natno$ in Axiom~\ref{N axiom},
the additional overt discrete objects that can be \emph{constructed} from them
are no bigger.

\begin{r@Corollary}
 In the free model, if an object $X$ has a basis indexed by any
overt discrete object $M$ then it has one indexed by $\natno$.
\end{r@Corollary}

\Proof Let $n:\natno,\; m:M\proves n\realises m$ be the relation defined
in the Proposition and $(\beta^m,A_m)$ the basis on $X$. Define
$$ \gamma^n \;\equiv\; \Lamb x.\Some m.n\realises m\land\beta^m x
   \quad\hbox{and}\quad
   D_n \;\equiv\; \Lamb\phi.\Some m.n\realises m\land A_m\phi.
$$
So $\gamma^n=\beta^m$ and $D_n=A_m$ if $n\realises m$,
but $\gamma^n=\bot$ and $D_n=\bot$ if $n\notin N'$.
Then, using the properties of $\realises$,
$(\gamma^n,D_n)$ is an effective basis because
\begin{eqnarray*}
  \Some n. D_n\phi \land \gamma^n x &\eq&
  \Some n m m'.n\realises m\land A_m\phi\land n\realises m'\land\beta^{m'} x \\
  &\eq& \Some m. A_m\phi\land\beta^m x \;\eq\;\phi x &\qeds
\end{eqnarray*}

\medskip

The next goal is Scott continuity, \ie preservation of directed joins.
Recall from Definition~\ref{DJ} that these are defined in terms
of a structure $(S,0,{+})$ that indexes two families
$$ s:S\ \proves\ \alpha_s:\Sigma\quad\hbox{and}\quad \phi^s:\Sigma^X. $$

\begin{r@Lemma}
\label{Kchoice} $\Gamma,\;\ell:\Fin(M)\ \proves\
(\All m\in\ell.\Some s:S.\alpha_s\land\phi^s m) \;\eq\;
(\Some s:S.\alpha_s\land\All m\in\ell.\phi^s m) $.
\end{r@Lemma}

\Proof We have to show $\Implies$, as $\Impliedby$ is easy.
For the base case, $\ell\equiv 0$, put $s\equiv 0$.
For the induction step, $\ell'=m::\ell$,
suppose by the induction hypothesis that
$$ \alpha_t\land\phi^t m  \quad\land\quad
   \alpha_s\land\All m\in\ell.\phi^s m
$$
Put $u\equiv s+t:S$, so $\alpha_u\eq\alpha_{s+t}\eq\alpha_s\land\alpha_t$
and $\phi^s,\phi^t\leq\phi^u$, so we have
$\alpha_u\land\All m\in{m::\ell}.\phi^u m$.

This proof is for lists, but a similar induction scheme works for
subsets too, \cf \cite[Section 6.5]{TaylorP:prafm}.
Since this equational hypothesis is of the form $\sigma\eq\top$,
it can be eliminated in favour of an open subobject of the context
\cite[\S 2]{TaylorP:insema}.
\qed

\begin{r@Lemma}
\label{F:S3Ncts}
Any $\Gamma\proves\F:\Sigma^3 N$ preserves the directed join
$\Gamma\proves\Some s.\alpha_s\land \Phi^s:\Sigma^{\Sigma^N}$.
\end{r@Lemma}

\Proof Using the previous lemma for $X\equiv \Sigma^N$ and $M\equiv\Fin N$,
\begin{eqnarray*}
  \F(\Some s.\alpha_s\land \Phi^s)
  &\eq&\Some L.\F A_L\land\All\ell\in L.\Some s.\alpha_s\land\Phi^s\beta^\ell
    &Lemma~\ref{SSN basis}\\
  &\eq&\Some L.\F A_L\land\Some s.\alpha_s\land\All\ell\in L.\Phi^s\beta^\ell\\
  &\eq&\Some s.\alpha_s\land\Some L.\F A_L\land\All\ell\in L.\Phi^s\beta^\ell\\
  &\eq&\Some s.\alpha_s\land \F \Phi^s &Lemma~\ref{SSN basis} \qEd
\end{eqnarray*}

\begin{r@Theorem}
\label{F:SSX cts} Any $\Gamma\proves F:\Sigma^{\Sigma^X}$
preserves the directed join $\Some s.\alpha_s\land\phi^s:\Sigma^X$.
\end{r@Theorem}

\Proof Making the same substitutions as in Lemma~\ref{SX lbasis},
\begin{eqnarray*}
  F(\Lamb x.\Some s.\alpha_s\land\phi^s x)
  &\eq& \F\cdot I\big(\Lamb x.\Some s.\alpha_s\land \Phi^s(i x)\big)\\
  &\eq& (\F\cdot I\cdot\Sigma^i)(\Some s.\alpha_s\land \Phi^s) \\
  &\eq& \Some s.\alpha_s\land(\F\cdot I\cdot\Sigma^i) \Phi^s
      &Lemma~\ref{F:S3Ncts}\\
  &\eq& \Some s.\alpha_s\land F\phi^s &similarly \qEd
  \vadjust{\nobreak}
\end{eqnarray*}
\closeupaline

\nobreak

\begin{r@Corollary}
 All $F:\Sigma^Y\to\Sigma^X$ preserve directed joins. \qed
\end{r@Corollary}

\goodbreak

We can now see the construction in Lemma~\ref{SX lbasis} as a composite,
of Lemma~\ref{make or basis} twice and the following result.
In order to apply $F:\Sigma^{\Sigma^X}$ to the basis decomposition of
$\phi:\Sigma^X$, the decomposition must be a directed join.

\begin{r@Lemma}
\label{SX basis} If $X$ has a $\lor$-basis then
$\Sigma^X$ has a prime $\land$-basis.
\end{r@Lemma}

\Proof For $F:\Sigma^{\Sigma^X}$ and $\phi:\Sigma^X$,
since the $\lor$-basis gives a directed join,
$$ F\phi \;\eq\; F(\Some n.A_n\phi\land\beta^n) \;\eq\;
  \Some n.A_n\phi\land F\beta^n. $$
This is the decomposition of $F$ with respect to the prime $\land$-basis with
$$ B^n\;\equiv\;A_n \quad\hbox{and}\quad \A_n\;\equiv\;\Lamb F.F\beta^n. $$
Notice that it reverses the imposed order on the indexing object. \qed

\smallskip

We can also go up to the second exponential with just one level of lists.

\begin{r@Lemma}
\label{SSX basis}  Let $(\beta^n,A_n)$ be a $\lor$-basis for $X$.
Then
$$ \B^\ell\equiv\Lamb F.\All n\in\ell.F\beta^n \quad\hbox{and}\quad
   \AAA_\ell\equiv\Lamb\F.\F(\Lamb\phi.\Some n\in\ell.A_n\phi) $$
define a prime $\land$-basis for $\Sigma^{\Sigma^X}$, \ie
we have a basis decomposition
$$ \F:\Sigma^3 X \ \proves\ \F \;=\; \Lamb F.\Some\ell.
    \F(\Lamb\phi.\Some n\in\ell.A_n\phi)\land\All n\in\ell.F\beta^n.  $$
\end{r@Lemma}

\Proof $(A_n,\Lamb F.F\beta^n)$ is a basis and
$ (\Lamb\phi.\Some n\in\ell.A_n\phi,\; \Lamb F.\All n\in\ell.F\beta^n) $
an $\lor$-basis for $\Sigma^X$ by
Lemmas \ref{SX basis} and~\ref{make or basis}.
Finally, Lemma~\ref{SX basis} gives the prime $\land$-basis for
$\Sigma^{\Sigma^X}$.  \qed

\medskip

As a corollary, we have another result similar to
Corollary~\ref{ovdisc=subquot}, 
although this one does restrict the size of overt discrete compact objects.
Such objects are called \textdf{Kuratowski finite}:
see \cite{TaylorP:insema} for more detail.
Beware that even finiteness is not the same as being a numeral,
for example a M\"obius band has a finite set of edges
that is locally equal to~2, but this set is not a numeral.

\begin{r@Theorem}
\label{Kfinsp}
An object $N$ is overt discrete compact iff it is listable,
and therefore a quotient of a numeral by an open equivalence relation.
\end{r@Theorem}

\Proof Expanding the quantifier $\forall_N$ using the filter $\lor$-basis
for $N$ indexed by $\List N$,
$$ \forall_N\phi \;\eq\; \Some\ell.\All n.n\in\ell\land\All m\in\ell.\phi m, $$
so $\top\eq\forall_N\top\eq\Some\ell.\All n.n\in\ell$.
The numeral is (the length of) $\ell$,
though we require the Existence Property
(\cf Remark~\ref{exists for bigvee}) to pin down its order and
so make it a quotient of a numeral.
\qed

\medbreak

As another corollary, we can complete the business of
Proposition~\ref{cpct=closed} and Lemma~\ref{finite AD=DA}.

\begin{r@Theorem}
\label{AD=DA} Let $A:\Sigma^Y\to\Sigma^X$ and $B:\Sigma^V\to\Sigma^U$.
\begin{letterlist}
\item If $A$ and $B$ both preserve $\top$ and $\land$ then the squares
  \begin{diagram}
    & \Sigma^X & \lTo^A & \Sigma^Y && X & \rHto^{\hayo A} & Y \\
    \Sigma^U & \Sigma^{X\times U} & \lTo^{A^U} & \Sigma^{Y\times U}
      && X\times U \shift-.5em & \rHto^{A\times U} & Y\times U\shift.5em && U\\
    \uTo<B & \uTo<{B^X} &?& \uTo>{B^Y} 
      && \dHto<{B\times X} &?& \dHto>{B\times Y} && \dHto>{\hayo B}\\
    \Sigma^V & \Sigma^{X\times V} & \lTo^{A^V} & \Sigma^{Y\times V} 
      && X\times V\shift-.5em & \rHto^{A\times V} & Y\times V\shift.5em && V
  \end{diagram}
  commute. Hence $\times$ is a symmetric monoidal structure
  in the category whose objects
  are those of $\S$ and whose morphisms $\hayo A:X\Hto Y$
  are the $A$ that preserve $\top$ and $\land$,
  \cf \cite[Section 3]{TaylorP:sobsc}.
  In particular, these maps commute with $\forall_K$,
  for any compact $K$.
\item If instead $A$ and $B$ both preserve $\bot$ and $\lor$
  then again the squares commute,
  so $\times$ is a symmetric monoidal structure in that category too,
  and these maps commute with $\exists_N$, for any overt $N$.
\end{letterlist}
\end{r@Theorem}

\Proof Lemma~\ref{finite AD=DA} and Theorem~\ref{F:SSX cts}. \qed

\begin{r@Corollary}
 Closed and compact subobjects coincide in any compact Hausdorff
object, whilst open and overt ones agree in any overt discrete object. \qed
\end{r@Corollary}

%============================================================================
\section{Primes and nuclei}\label{prime/nucleus}

In this section we prove the finitary lattice-theoretic characterisations
of primes and nuclei that we announced in Remark~\ref{formulate prime/nucleus}.

Since homomorphisms of the four lattice connectives also preserve
directed joins, we might expect them to be frame homomorphisms.
However, they only preserve joins indexed by \emph{overt} objects,
not ``arbitrary'' ones.
On the other hand, objects and inverse image maps are defined
in Abstract Stone Duality as algebras and homomorphisms
\emph{for the monad} corresponding to $\Sigma^\blank\adjoint\Sigma^\blank$,
rather than for an infinitary algebraic theory.

We actually consider ``Curried'' homomorphisms.

\begin{r@Definition}
\label{def prime} The term
$\Gamma\proves P:\Sigma^{\Sigma^X}$ is \textdf{prime}
\cite[Definition~8.1]{TaylorP:sobsc} if
$$\Gamma,\F:\Sigma^3 X \proves \F P \eq
    P \big(\Lamb x.\F(\Lamb\phi.\phi x)\big).
$$
Axiom~\ref{sober} says that we may then introduce $a=\focus P$
such that $P=\Lamb\phi.\phi a$.

\smallskip

Plainly if $P\equiv\Lamb\phi.\phi a$ for some $a$ then
$$  P\top \;\eq\top, \qquad
    P(\phi\land\psi) \;\eq\; \phi a\land\psi a \;\eq\; P\phi\land P\psi
$$
and similarly for $\bot$ and $\lor$.
We know that frame homomorphisms (as defined \emph{externally}
using infinitary lattices) agree with Eilenberg--Moore homomorphisms
in the case of the classical models \cite{TaylorP:sobsc,TaylorP:subasd}.
Now we can use bases to prove a similar result for
the \emph{internal finitary} lattice structure in our category.
This means that we can import at least some of the familiar lattice-theoretic
arguments about topology into our category.
\end{r@Definition}

\begin{r@Theorem}
\label{tfae prime} $P:\Sigma^{\Sigma^X}$
is prime iff it preserves $\top$, $\bot$, $\land$ and $\lor$.
\end{r@Theorem}

\Proof $[\Implies]$ 
Put $\F\equiv\Lamb F.F\phi\land F\psi$, so $\F P\eq P\phi\land P\psi$, whilst
$$ P\big(\Lamb x.\F(\Lamb\phi.\phi x)\big) \;\eq\;
   P(\Lamb x.\phi x\land\psi x) \;\eq\;
   P(\phi\land\psi).
$$
The other connectives are handled in the same way.

$[\Impliedby]$ First note that, using Axiom~\ref{Phoa} for $P$ and $\alpha$,
and $P\bot\eq\bot$,
$$ P(\alpha\land\phi)
   \;\eq\;  P\bot\lor\alpha\land P\phi
   \;\eq\; \alpha\land P\phi.
$$
Second, $P$ commutes with $\exists$ by Theorem~\ref{AD=DA}(b), so
\begin{eqnarray*}
  \F P &\eq&
  \Some\ell. \F(\Lamb\phi.\Some n\in\ell.A_n\phi)\land\All n\in\ell.P\beta^n
   &Lemma~\ref{SSX basis}\\
  &\eq&
  \Some\ell. \F(\Lamb\phi.\Some n\in\ell.A_n\phi)\land P(\All n\in\ell.\beta^n)
    &$P$ preserves $\land$, $\top$\\
  &\eq& \Some\ell.P\big(\F(\Lamb\phi.\Some n\in\ell.A_n\phi)
  \land\All n\in\ell.\beta^n\big) &above\\ 
  &\eq& P\big(\Some\ell.\F(\Lamb\phi.\Some n\in\ell.A_n\phi)
  \land\All n\in\ell.\beta^n\big) \\
  &\eq& P\big(\Lamb x.\Some\ell.\F(\Lamb\phi.\Some n\in\ell.A_n\phi)
  \land\All n\in\ell.(\Lamb\phi.\phi x)\beta^n\big) \\
  &\eq& P\big(\Lamb x.\F(\Lamb\phi.\phi x)\big). &Lemma~\ref{SSX basis} \qEd
\end{eqnarray*}

\begin{r@Corollary}
\label{tfae hom}
$H:\Sigma^X\to\Sigma^Y$ is an Eilenberg--Moore homomorphism
iff it is a lattice homomorphism.  In this case, it is of the form
$H=\Sigma^f$ for some unique $f:Y\to X$. \qed
\end{r@Corollary}

\medskip

We have in particular a way of introducing \emph{points} of an object by
finitary lattice-theoretic arguments.  By this method we can derive a
familiar domain-theoretic result, but beware that, as the objects of
our category denote locally compact objects and not merely domains,
the order by no means characterises the objects and morphisms.

\begin{r@Theorem}
\label{dirsup}
 Every object has and every morphism preserves directed joins.
\end{r@Theorem}

\Proof \cite[Lemma II 1.11]{JohnstonePT:stos}
Let $\Gamma,\; s:S\ \proves\ a_s:X$ be a directed family in $X$
with respect to the intrinsic order (Definition~\ref{order}), then
$$ \Gamma,\; s:S \ \proves\ \Lamb\phi.\phi a_s :\Sigma^{\Sigma^X} $$
is a directed family of primes.
But primes are characterised lattice-theoretically, and the property
is preserved by directed joins (Theorem~\ref{F:SSX cts} for $\F$), so
$$ \Gamma \ \proves\ P\;\equiv\;\Lamb\phi.\Some s.\phi a_s :\Sigma^{\Sigma^X}$$
is also prime, and $P=\Lamb\phi.\phi a$,
where $a\equiv\focus P:X$.

I claim that this is the join of the given family.
By Definition~\ref{order}, the order relation $a_s\leq_X a$ means 
$\Lamb\phi.\phi a_s\leq\Lamb\phi.\phi a\equiv P$,
which we have by the definition of $P$, 
whilst similarly if $\Gamma,\; s:S\ \proves\ a_s\leq b$ then
$\Lamb\phi.\phi a_s\leq\Lamb\phi.\phi b$ so $P\leq\Lamb\phi.\phi b$
as $P$ is the join, and then $a\leq b$.

By a similar argument, given $f:X\to Y$, we have $f(a_s)\leq_Y f a$
since $\psi:\Sigma^Y$ provides $\psi\cdot f:\Sigma^X$, and
if $f(a_s)\leq b$ then $f(a)\leq b$, so $f$ preserves the join. \qed

\begin{r@Corollary}
Suppose that an object $X$ has a least element $\bot$ in its intrinsic order.
Then any $f:X\to X$ has a least fixed point, namely
$$ a\;\equiv\;\focus(\Lamb\phi.\Some n.\phi a_n),
   \quad\hbox{where}\
   a_0\;\equiv\;\bot
   \ \hbox{and}\
   a_{n+1}\;\equiv\; f a_n.
   \eqno\qEd
$$
\end{r@Corollary}

\bigskip

Finally, we show that the equations in Lemma~\ref{E@} that
characterised the Scott-continuous functions $E$ that are of the
form $I\cdot\Sigma^i$ for a $\Sigma$-split sub\emph{locale}
are also valid for $\Sigma$-split subobjects in abstract Stone duality.

\begin{r@Definition}
\label{def nucleus}
Recall from \cite[Definition~4.3]{TaylorP:subasd}
that $E:\Sigma^X\to\Sigma^X$ is called a \textdf{nucleus} if
$$ \F:\Sigma^3 X\ \proves\
 E\big(\Lamb x.\F(\Lamb\phi. \phi x)\big) \;=\;
 E\big(\Lamb x.\F(\Lamb\phi. E\phi x)\big):\Sigma^X.
$$
(This equation arises from Beck's monadicity theorem,
and is applicable without assuming any lattice structure on $\Sigma$.)
\end{r@Definition}

\begin{r@Lemma}
\label{nucleus->E@} If $E$ is a nucleus then
$$ \phi,\psi:\Sigma^X \ \proves\
   E(\phi\land\psi)\;=\;E(E\phi\land E\psi)
   \quad\hbox{and}\quad
   E(\phi\lor\psi)=E(E\phi\lor E\psi).
$$
\end{r@Lemma}

\Proof Putting $\F\equiv\Lamb F.F\phi\land F\psi$,
$$ E\big(\Lamb x.\F(\Lamb\phi.\phi x)\big)
   \;=\; E(\Lamb x.\phi x\land \psi x)
   \;=\; E(\phi\land\psi) \eqno LHS%
$$
$$ E\big(\Lamb x.\F(\Lamb\phi.E\phi x)\big) \;=\;
    E(\Lamb x.E\phi x\land E\psi x) \;=\;
    E(E\phi\land E\psi), \eqno RHS
$$
are equal. The argument for $\lor$ is the same. \qed

\medskip

For the converse, first observe that the equations allow us to insert
or remove $E$s as we please in any sub-term of a lattice expression,
so long as $E$ is applied to the whole expression.
In particular, $E\phi=E(E\phi)$ and
$E(\phi\lor\psi\lor\theta)=E(E\phi\lor E\psi\lor\theta)$ (\emph{sic}).

\begin{r@Lemma}
\label{Elist} Although we needn't have $E\top=\top$ or $E\bot=\bot$,
we may extend the binary $\lor$-formula to finite (possibly empty) sets
$\ell:\Fin(N)$:
$$ E(\Some n\in\ell.\alpha_n\land\phi^n) \;=\; 
   E(\Some n\in\ell.\alpha_n\land E\phi^n),
$$
where $ n:N\proves \alpha_n:\Sigma$ and $\phi^n:\Sigma^X$.
Similarly but more simply, from the $\land$-equation,
$$ E(\All n\in\ell.\phi^n) \;=\; E(\All n\in\ell. E\phi^n). $$
\end{r@Lemma}

\Proof The base case of the induction, $\ell\equiv 0$, is $E\bot=E\bot$,
and the next case is
$$ E(\alpha\land\phi) \;=\; E\bot\lor\alpha\land E\phi
    \;=\; E\bot\lor\alpha\land E(E\phi)
    \;=\; E(\alpha\land E\phi)
$$
by the Phoa principle.
For the induction step \cite[\S 2]{TaylorP:insema},
\begin{eqnarray*}
  E(\Some n\in{m::\ell}.\alpha_n\land\phi^n) \hidewidth \\
  \qquad &=&
  E\big((\Some n\in\ell.\alpha_n\land\phi^n)\lor(\alpha_m\land\phi^m)\big) \\
  &=&
  E\big(E(\Some n\in\ell.\alpha_n\land\phi^n)\lor E(\alpha_m\land\phi^m)\big)
    &$\lor$-equation\\
  &=& E\big(E(\Some n\in\ell.\alpha_n\land E\phi^n)\lor
     E(\alpha_m\land E\phi^m)\big) &ind.~hyp.\\
  &=& E\big((\Some n\in\ell.\alpha_n\land E\phi^n)\lor
     (\alpha_m\land E\phi^m)\big) &$\lor$-equation\\
  &=& E(\Some n\in{m::\ell}.\alpha_n\land E\phi^n) \qEd
\end{eqnarray*}

\begin{r@Lemma}
\label{Esome} The $\exists$ equation extends by Scott continuity
(Proposition~\ref{F:SSX cts}).
\end{r@Lemma}

\Proof\closeupaline
\begin{eqnarray*}
   E(\Some n:N.\alpha_n\land\phi^n) 
   &=& E(\Some\ell:\Fin N.\Some n\in\ell.\alpha_n\land\phi^n) \\
   &=& \Some\ell.E(\Some n\in\ell.\alpha_n\land\phi^n)
      &Proposition~\ref{F:SSX cts}\\
   &=& \Some\ell.E(\Some n\in\ell.\alpha_n\land E\phi^n)
      &Lemma~\ref{Elist}\\
   &=& E(\Some n:N.\alpha_n\land E\phi^n) &similarly \qEd
\end{eqnarray*}
 
\begin{r@Theorem}
\label{lattice nucleus} $E$ is a nucleus iff it satisfies
$$ \phi,\psi:\Sigma^X \ \proves\
   E(\phi\land\psi)\;=\;E(E\phi\land E\psi)
   \quad\hbox{and}\quad
   E(\phi\lor\psi)=E(E\phi\lor E\psi).
$$
\end{r@Theorem}

\Proof We expand $\F:\Sigma^3 X$ in the defining equation for a nucleus
\begin{eqnarray*}
   E\big(\Lamb x.\F(\Lamb\phi.E\phi x)\big) 
   &=& E\big(\Some L.\F A_L \land\All\ell\in L. E\beta^\ell\big)
       &Proposition~\ref{SX basis}\\
   &=& E\big(\Some L.\F A_L \land E(\All\ell\in L.E\beta^\ell)\big)
       &Lemma~\ref{Esome}\\
   &=& E\big(\Some L.\F A_L \land E(\All\ell\in L.\beta^\ell)\big)
       &Lemma~\ref{Elist}\\
   &=& E\big(\Some L.\F A_L \land\All\ell\in L.\beta^\ell\big)
       &Lemma~\ref{Esome}\\
   &=& E\big(\Lamb x.\F(\Lamb\phi.\phi x)\big) &Proposition~\ref{SX basis} \qEd
\end{eqnarray*}

%============================================================================
\section{The way-below relation}\label{way-below}

We now introduce the last (g) of our abstract characterisations
of local compactness in the Introduction. %(Remark~\ref{summary loc cpct}(g)).
Recall from Section~\ref{intro locales}
that any continuous distributive lattice carries a binary relation 
(written $\ll$ and called ``way-below'') such that
$$
\bot\ll\gamma
\qquad
\begin{prooftree}
  \alpha\ll\gamma \quad \beta\ll\gamma \Justifies \alpha\lor\beta\ll\gamma
\end{prooftree}
\qquad
\begin{prooftree}
  \alpha'\leq\alpha\ll\beta\leq\beta' \justifies \alpha'\ll\beta'
\end{prooftree}
\qquad
\begin{prooftree}
  \alpha\ll\gamma \Justifies \Some\beta.\alpha\ll\beta\ll\gamma 
\end{prooftree}
$$ $$
\begin{prooftree}
  \alpha\ll\beta \quad \beta\ll\phi \quad \beta\ll\psi
  \justifies
  \alpha\ll(\phi\land\psi)
\end{prooftree}
\qquad\qquad
\begin{prooftree}
  \alpha\ll\beta\lor\gamma
  \Justifies
  \Some\beta'\gamma'.\ 
    \alpha\ll\beta'\lor\gamma'\quad \beta'\ll\beta\quad\gamma'\ll\gamma
\end{prooftree}
$$

\begin{r@Notation}
\label{<< notn}
We define a new relation $n\waybelow m$ as $A_n\beta^m$.

\emph{This is an open binary relation (term of type $\Sigma$ with two
  free variables) on the overt discrete object $N$ of indices of a
  basis, not on the lattice $\Sigma^X$.}  It is an ``imposed''
structure on $N$ in the sense of Remark~\ref{use of N}.

For most of the results of this section, 
we shall require $(\beta^n,A_n)$ to be an $\lor$-basis for~$X$.
\end{r@Notation}

\begin{r@Examples}
(Not all of these are $\lor$-bases.)
\begin{letterlist}
\item Let $\beta^n$ classify $U^n\subset X$,
  and $A_n=\Lamb\phi.(K^n\subset\phi)$ in a locally compact sober space.
  Then $n\waybelow m$ means that $K^n\subset U^m$.
  This is consistent with Notation~\ref{spatial <<}
  if we identify the basis element $n$ with the pair $(U^n\subset K^n)$.
\item Let $A_n=\Lamb\phi.(\beta^n\ll\phi)$ in a continuous lattice.
  Then $n\waybelow m$ means that $\beta^n\ll\beta^m$,
  \cf Definition~\ref{rec cts latt}.
\item In the interval basis on $\realno$ in Example~\ref{eg R2},
  $<q,\delta>\waybelow<p,\epsilon>$
  means that $[q,\delta]\subset(p,\epsilon)$,
  \ie $p-\epsilon\lt q-\delta\leq q+\delta\lt p+\epsilon$.
\item In the basis of disjoint pairs of opens, $(U^n\disjoint V_n)$,
  for a compact Hausdorff space in Example~\ref{cpct Hdf base},
  $n\waybelow m$ means that $V_n\cup U^m=X$.
\item In the prime basis $(\setof n,\eta n)$ for $N$ in Example~\ref{N basis},
  $n\waybelow m$ just when $n=m$.
\item In the $\Fin(N)$-indexed filter $\lor$-basis on $N$ 
  in Proposition~\ref{N bases},
  $\ell\waybelow\ell'$ iff $\ell\subset\ell'$.
\item In the prime $\land$-basis on $\Sigma^N$ in Example~\ref{SN basis},
  on the other hand, $\ell\waybelow\ell'$ iff $\ell'\subset\ell$.
\item More generally, in the prime $\land$-basis on $\Sigma^X$ 
  derived from an $\lor$-basis on $X$ in Lemma~\ref{SX basis},
  $(n\waybelow_{\Sigma^X}m)$ iff $(m\waybelow_X n)$.  
\item\label{KKN<<}
  In the $\Fin\big(\Fin(N)\big)$-indexed filter lattice basis on $\Sigma^N$
  in Proposition~\ref{N bases},
  $$ L\waybelow R \;\equiv\; \A_L B^R \;\eq\; R\shp\subset L \;\equiv\;
   \All\ell\in L.\Some\ell'\in R.(\ell'\subset\ell),
  $$
  where $R\shp\subset L$ is known as the \textdf{upper order} on subsets
  induced by the relation $\subset$ on elements.
\item In the prime $\land$-basis on $\Sigma^{\Sigma^N}$
  in Lemma~\ref{SSN basis},
  $L\waybelow R$ iff $L\shp\subset R$.
\item In the $\lor$-basis on $\Sigma^{\Sigma^X}$
  derived from an $\lor$-basis on $X$ in Lemma~\ref{SSX basis},
  $$ (\ell\waybelow_{\Sigma^2 X}\ell') \iff (\ell\shp\waybelow_X\ell'). $$  
\item Any stably locally compact object $X$
  has a filter lattice basis $(N,0,1,{+},{\star},{\waybelow})$
  such that the opposite $(N,1,0,{\star},{+},{\wayabove})$
  is the basis of another such object, known as its \textdf{Lawson dual}
  (Proposition~\ref{ASD dual Wilker}).
\end{letterlist}
\end{r@Examples}

\medbreak

\noindent
Our first result just restates the assumption of an $\lor$-basis,
\cf Lemmas \ref{classical 0+<<} and~\ref{locale 0+<<}:

\begin{r@Lemma}
\label{0+ll} $0\waybelow p$, whilst
$n+m\waybelow p$ iff both $n\waybelow p$ and $m\waybelow p$.
\end{r@Lemma}

\Proof $A_0\beta^p\eq\top$ and $A_{n+m}\beta^p\eq A_n\beta^p\land A_m\beta^p$.
\qed

\smallskip

In a continuous lattice, $\alpha\ll\beta$ implies $\alpha\leq\beta$,
but we have no similar property relating $\waybelow$ to~$\baseleq$.
We shall see the reasons for this in Section~\ref{DLmonad}.
But we do have two properties that carry most of the force of
$\alpha\ll\beta\Implies\alpha\leq\beta$.
We call them \textdf{roundedness}.
The second also incorporates many of the uses of directed joins
and Scott continuity into a notation that will become increasingly
more like discrete mathematics than
it resembles the technology of traditional topology.

\begin{r@Lemma}
\label{<<b}\label{<<A} $\beta^n=\Some m.(m\waybelow n)\land\beta^m$ and 
 $A_m=\Some n.A_n\land(m\waybelow n)$.
\end{r@Lemma}

\Proof The first is simply the basis expansion of $\beta^n$.
For the second, we apply $A_m$ to the basis expansion of $\phi$, so
$$ A_m\phi \;\eq\; \Some n.A_n\phi\land A_m\beta^n \;\eq\;
   (\Some n.A_n\land A_m\beta^n)\phi, $$
since $A_n$ preserves directed joins (Theorem~\ref{F:SSX cts}). \qed

\begin{r@Corollary}
\label{<<-leq}
If $m\waybelow n$ then $\beta^m\leq\beta^n$ and $A_m\geq A_n$. \qed
\end{r@Corollary}

\begin{r@Corollary}
\label{<.<<.<}
If $m\waybelow n$, $A_{m'}\geq A_m$ and $\beta^n\leq\beta^{n'}$,
then $m'\waybelow n'$.
\end{r@Corollary}

\Proof
$(m\waybelow n)\equiv A_m\beta^n\Implies A_{m'}\beta^{n'}\equiv(m'\waybelow n')$.
\qed

\begin{r@Corollary}
\label{interpolation} The relation $\waybelow$ satisfies
transitivity and the \textdf{interpolation lemma}:
$$ (m\waybelow n) \;\eq\; (\Some k. m\waybelow k\waybelow n). $$
\end{r@Corollary}

\Proof $A_m\beta^n \eq (\Some k.A_k\land m\waybelow k)\beta^n \eq
\Some k.A_k\beta^n\land m\waybelow k$. \qed

\medskip

Now we consider the interaction between $\waybelow$ and the lattice
structures $(\top,\bot,{\land},{\lor})$ and $(1,0,{+},{\star})$.
Of course, for this we need a lattice basis.

\begin{r@Lemma}
\label{*=<<**} As directed joins,
$$\phi\land\psi \;=\; \Some p q.\beta^{p\star q}\land A_p\phi\land A_q\psi
  \quad\hbox{and}\quad
  \phi\lor\psi = \Some p q.\beta^{p+q}\land A_p\phi\land A_q\psi. $$
\end{r@Lemma}

\Proof The first is the Frobenius law, since
$\beta^{p\star q}=\beta^p\land\beta^q$.

The second uses Lemma~\ref{or distrib dirsup}: we obtain the expression
$$ \phi\lor\psi \;=\; \Some p.A_p\phi\land(\beta^p\lor\psi) $$
from the basis expansion $\phi=\Some p.A_p\phi\land\beta^p$
by adding $\psi$ to the $0$th term (since $A_0\phi\eq\top$ and $\beta^0=\bot$)
and, harmlessly, $A_p\phi\land\psi$ to the other terms.
Similarly,
$$ \beta^p\lor\psi \;=\; \Some q.A_q\psi\land(\beta^p\lor\beta^q) \;=\;
   \Some q.A_q\psi\land\beta^{p+q}. $$
The joins are directed because because $A_0\phi\land A_0\psi\eq\top$ and
$$ (A_{p_1}\phi \land A_{q_1}\psi) \;\land\; (A_{p_2}\phi\land A_{q_2}\psi)
    \;\eq\; (A_{p_1+p_2}\phi\land A_{q_1+q_2}\psi). \eqno\qEd $$

\begin{r@Lemma}
\label{ll square}  For a lattice basis,
$A_n\top\eq(n\waybelow1)$, $A_n\bot\eq(n\waybelow0)$ and
\begin{eqnarray*}
  A_n(\phi\land\psi) &\eq&
 \Some p q. (n\waybelow p\,\star\, q)\land A_p\phi\land A_q\psi \\
  A_n(\phi\lor\psi) &\eq&
 \Some p q. (n\waybelow p+q)\land A_p\phi\land A_q\psi
\end{eqnarray*} 
\end{r@Lemma}

\Proof The first two are $A_n\beta^1$ and $A_n\beta^0$.
The other two are $\Some p q.A_n(\beta^{p\star q})\land A_p\phi\land A_q\psi$
and $\Some p q.A_n(\beta^{p+q})\land A_p\phi\land A_q\psi$,
which are $A_n$ applied to the directed joins in Lemma~\ref{*=<<**}. \qed

\smallskip

\begin{r@Proposition}
\label{filter lattice <<*}
The lattice basis $(\beta^n,A_n)$ is a filter basis iff $1\waybelow1$ and
$$
\begin{prooftree}
  m\waybelow p\star q
  \Justifies
  m\waybelow p \qquad m\waybelow q
\end{prooftree}
$$
\end{r@Proposition}

\Proof $(n\waybelow1)\eq A_n\beta^1\eq A_n\top$,
but recall that $n\waybelow1$ for all $n$ iff $1\waybelow1$.

The displayed rule is
$A_m\beta^p\land A_m\beta^q \eq A_m\beta^{p\star q} \equiv
 A_m(\beta^p\land\beta^q)$.
Given this, by the Frobenius law and Lemma~\ref{ll square},
\begin{eqnarray*}
  A_m\phi\land A_m\psi
 &\eq& \Some p q.A_p\phi\land A_q\psi\land A_m\beta^p\land A_m\beta^q \\
 &\eq& \Some p q.A_p\phi\land A_q\psi\land (m\waybelow p)\land (m\waybelow q)\\
 &\eq& \Some p q.A_p\phi\land A_q\psi\land (m\waybelow p\star q)\\
 &\eq& A_m(\phi\land\psi) &\qeds
\end{eqnarray*}

\smallskip\goodbreak

If we don't have a filter basis, we have to let $n$ ``slip'' by
$n\waybelow m$, \cf Lemma~\ref{locale <<*}:

\begin{r@Lemma}
\label{<<* rule} For any lattice basis,
$$
\begin{prooftree}
  n\waybelow p\star q
  \Justifies
  \Some m.\ (n\waybelow m)\ \land\ (m\waybelow p)\ \land\ (m\waybelow q)
\end{prooftree}
$$
\end{r@Lemma}

\Proof Downwards, interpolate $n\waybelow m\waybelow p\star q\baseleq p,q$,
then $m\waybelow p,q$ by monotonicity.
Conversely, using Corollary~\ref{<<-leq}, if $\A_n\beta^m\eq\top$ and
$\beta^m\leq\beta^{p\star q}$ then $A_n\beta^{p\star q}\eq\top$. \qed

\medskip

The corresponding result for $\lor$ is our version of the Wilker property,
\cf Proposition~\ref{classical Wilker} and Lemma~\ref{locale Wilker}.

\begin{r@Lemma}
\label{basis Wilker}
$$
\begin{prooftree}
  n\waybelow p+q
  \Justifies
  \Some p' q'.(n\waybelow p'+q')\land(p'\waybelow p)\land(q'\waybelow q)
\end{prooftree}
$$
\end{r@Lemma}

\Proof Lemma~\ref{ll square} with $\phi\equiv\beta^p$ and
$\psi\equiv\beta^q$.\qed

\medskip

We shall summarise these rules in Definition~\ref{abs basis}.

\bigbreak

There are special results that we have in the cases of overt and
compact objects.
We already know that $1\waybelow 1$ iff the object is compact
(\cf Lemma~\ref{Anbn=1}),
but the lattice dual characterisation of overtness cannot be
$0\waybelow 0$, as that always happens.

\begin{r@Lemma}
\label{n<<0|Eybmy} If $X$ is overt then
$$ (n\waybelow m) \;\Implies\; (n\waybelow 0) \lor (\Some y.\beta^m y)
   \quad\hbox{but}\quad
   (n\waybelow 0)\land(\Some x.\beta^n x)\eq\bot.
$$
\end{r@Lemma}

\Proof $\phi x\Implies \Some y.\phi y$ so, using the Phoa principle
(Axiom~\ref{Phoa}),
$$ A_n\phi
    \;\Implies\; A_n(\Lamb x.\Some y.\phi y) 
    \;\eq\; A_n\bot\lor\Some y.\phi y\land A_n\top
    \;\Implies\; (n\waybelow 0) \lor \Some y.\phi y. $$
Putting $\phi\equiv\beta^m$,\closeupaline
$$ (n\waybelow m) \;\equiv\; A_n\beta^m \;\Implies\;
   (n\waybelow 0)\lor \Some y.\beta^m y, $$
whilst
$ (n\waybelow 0) \land \Some x.\beta^n x
   \;\eq\; \Some x.A_n\beta^0\land\beta^n x
   \;\Implies\; \Some x.\Some n.A_n\bot\land\beta^n x
   \;\eq\; \Some x.\bot x
   \;\eq\; \bot.
$
\qed

\medskip

By a similar argument, we have Johnstone's ``Townsend--Thoresen Lemma''
\cite{JohnstonePT:opele},
$$ (n\waybelow p+q)\ \Implies\ (n\waybelow p) \lor (\Some y.\beta^q y). $$

\begin{r@Corollary}
\label{cpctovdec} Given a compact overt object, it's
decidable whether it's empty or inhabited, \cf Corollary~\ref{AKbot}.
\end{r@Corollary}

\Proof As $(1\waybelow1)\eq\top$, the Lemma makes
$(1\waybelow0)$ and $\Some x.\beta^1 x\equiv\Some x.\top$ complementary.
\qed

Compact overt subobjects are, in fact, particularly well behaved
in discrete \cite{TaylorP:insema} and Hausdorff \cite{TaylorP:lamcra}
objects.

Sections~\ref{DLmonad}--\ref{matrix} are devoted to proving that
the rules above for $\waybelow$ are complete,
in the sense that from any 
\emph{abstract basis} $(N,0,1,{+},{\star},{\waybelow})$ satisfying
them we may recover an object of the category.
This proof is very technical,
so in the next section we consider a special case,
in which the lattice structure plays a much lighter role.

%============================================================================
\section{Domain theory in ASD}\label{domain theory}

In this section we provide the basic tools for implementing domain theory
within ASD.
Denotational semantics can then be developed in the usual way using this,
allowing programming languages such as Gordon Plotkin's PCF
\cite{PlotkinGD:lcfcpl} to be interpreted in ASD,
and with them the programs required in
Definitions \ref{rec loc cpct}, \ref{rec cts} and~\ref{rec cts latt}
for \emph{computable} bases and functions.
The difference is that classical domain theory leads to the dead end
of set theory, whilst ASD can be translated back into programming languages.

However, the denotational semantics of programming languages
will have to wait for a separate investigation.
In this section we merely provide a few results that connect ASD with
the \emph{topological} aspects of domain theory
such as \cite{JungA:mullsc,JungA:stacsc,JungA:duacvo}
that rely on manipulation of bases (in the sense of that subject).

We first characterise the objects that have prime bases in our sense,
and are \textdf{continuous dcpos} in the classical setting.
The following results are essentially those of \cite[\S 2]{SmythMB:effgd},
where $\trless$ is called an \emph{R-structure}.
However, whereas Michael Smyth imposed a separate theory of recursion
on the existing set-based theory of continuous dcpos,
recursion is intrinsic to ASD.
We also form the domain of rounded ideals
using our own technology, which was developed in \cite{TaylorP:subasd}.

\begin{r@Definition}
\label{fir} A \textdf{directed interpolative relation}
$n,m:N\proves(n\trless m):\Sigma$ on an overt discrete object $N$
is an open relation that satisfies:
\begin{letterlist}
\item transitivity and interpolation:
  $n,m:N\ \proves\ 
  (n\trless m) \;\eq\; \Some k.(n\trless k)\land(k\trless m)$;
\item extrapolation downwards:
  $n:N \ \proves\ \Some k.(k\trless n)\eq\top$;
\item directedness:
  $n,r,s:N \ \proves\ (r\trless n)\land(s\trless n) \;\eq\;
             \Some k.(r\trless k)\land(s\trless k)\land(k\trless n)$.
\end{letterlist}

\smallskip

\noindent A \textdf{rounded ideal} for $(N,{\trless})$ is a predicate
$\Gamma\proves\xi:\Sigma^N$
that is
\begin{letterlist}
\item rounded:
  $\Gamma,\; n:N \ \proves\ \xi n \;\eq\; \Some k.(n\trless k)\land\xi k$;
\item inhabited:
  $\Gamma \ \proves\ \Some k.\xi k \;\eq\; \top$;
\item an ideal:
  $\Gamma,\; r,s:N \ \proves\ \xi r\land\xi s \;\Implies\;
     \Some k.(r\trless k) \land (s\trless k)\land\xi k$.
\end{letterlist}
\end{r@Definition}

\begin{r@Lemma}
For each $n:N$,\quad
$\tridown n \equiv \Lamb k.k\trless n$\quad is a rounded ideal. \qed
\end{r@Lemma}

\smallskip

We shall see that the $\tridown n$ provide the ``basic points'' $p_n$
(\cf Definition~\ref{rec loc cpct}) of continuous and algebraic dcpos.
In fact, the situation is summed up by the next result,
whose converse we shall prove in Theorem~\ref{prime basis}.

\begin{r@Proposition}
\label{prime basis->dpco}
If $X$ has an $N$-indexed prime basis, so
$$ \phi x \;\eq\; \Some n. \phi p_n \land \beta^n x, $$
then $\wayabove$ is a directed interpolative relation,
and $x:X\proves\xi\equiv\Lamb n.\beta^n x$ is a rounded ideal, with
$$ x \;=\; \dirsup\collect{p_n}{\xi n}. $$
\end{r@Proposition}

\noindent
Beware that $\waybelow$ and $\trless$ have opposite senses,
following the conventions for continuous lattices,
\cf the use of $\ll$ for compact and open subspaces
in Notation~\ref{spatial <<}:
we think of $\trless$ as relating points.
Thus $m\wayabove n$ and $m\trless n$ mean that $A_m\ll A_n$ and $p_m\ll p_n$,
but $\beta^n\ll\beta^m$.

\Proof We have already proved these properties when $\waybelow$ arises
from a filter $\lor$-basis.  An $\lor$-basis was needed in order to
use Scott continuity for $A_n$, but in this case $A_n\equiv\Lamb\phi.\phi p_n$
preserves finite joins too, so directedness is redundant.
Briefly,
\begin{letterlist}
\item the basis expansion of $(n\waybelow m)\equiv\beta^m p_n$
  gives transitivity and interpolation, as in Corollary~\ref{interpolation};
\item that of $(\Lamb x.\top)p_n$ gives extrapolation; whilst
\item that of $(\beta^r\land\beta^s)p_n$ gives directedness,
  as in Proposition~\ref{filter lattice <<*}.
\end{letterlist}
The directed join is $\focus(\Lamb\phi.{-})$ applied to the
basis expansion, as in Theorem~\ref{dirsup}. \qed

\medskip

\begin{r@Examples}
Here again are the prime bases that we have encountered.
\begin{letterlist}
\item Equality or any open equivalence relation on an overt discrete object
  is a directed interpolative relation,
  whose rounded ideals are the equivalence classes.
\item $\Sigma^N$ is the object of (rounded) ideals for
  list or subset inclusion in $\Fin(N)$.
\item If $X$ has an $N$-indexed $\lor$-basis,
  $\waybelow_X$ is a directed interpolative relation on~$N$,
  and $\Sigma^X$ is the object of rounded ideals for this
  (recall from Lemma~\ref{SX basis} that $\waybelow_{\Sigma^X}$
  is $\wayabove_X$).
\item\label{upd basis} In the case of an $N$-indexed \emph{filter lattice} basis
  (so $X$ is stably locally compact),
  the opposite relation ($\wayabove_X$) is also directed interpolative, because
  $$ \Some k.A_n\beta^k \;\eq\; A_n\beta^1 \;\eq\; A_n\top \;\eq\;\top
     \qquad % \quad\hbox{and}\quad
     A_n\beta^r\land A_n\beta^s \;\eq\; A_n(\beta^r\land\beta^s)
     \;\eq\; A_n(\beta^{r+s}).
  $$
\item\label{class upd} Classically, the rounded ideals for $\wayabove$ are
  in order-reversing bijection with the compact saturated subspaces
  of~$X$ (Theorem~\ref{cpct <<}).
\qed
\end{letterlist}
\end{r@Examples}

\begin{r@Lemma}
Any \emph{reflexive} transitive relation is directed interpolative,
and then all ideals are rounded. \qed
\end{r@Lemma}

\smallskip

\begin{r@Examples}
\begin{letterlist}
\item If the way-below relation $\waybelow$ of the space $X$
  in Proposition~\ref{prime basis->dpco} is reflexive
  then $X$ is called an \textdf{algebraic dcpo},
  and the $p_n$ are its so-called \textdf{finite} or \textdf{compact} elements.
  
\item\label{eg cts latt basis}
  The space $X$ is a continuous lattice iff
  it has a basis $(N,{\trless})$ that's a semilattice,
  \ie the extrapolative and directed conditions are strengthened
  to a least element and binary joins respectively,
  instead of the extra/interpolant~$k$.

\item In a \textdf{Scott domain}, the join is \emph{conditional}:
  the join $r+s$ exists in $N$ whenever there is \emph{some} upper bound~$n$.
  Any Scott domain is a closed subspace of a continuous lattice,
  but if they are to form a cartesian closed category,
  they must also be \emph{overt} \cite{TaylorP:pcfasd}.
  Indeed, the basis $N$ is overt in the formulation of this section.
  Instead of a conditional join, we may add a formal top element to~$N$,
  or, equivalently, specify what part of an enlarged basis is to define
  the required domain by means of a \emph{decidable}
  \textdf{consistency predicate}, \cf Proposition~\ref{R-SN}
  \cite{ScottDS:domds}.

\item\label{eg arith univ basis}
  With pushouts (joins of pairs that are bounded \emph{below})
  but not necessarily a least element,
  we have a class of objects that includes the overt discrete ones
  and is closed under $\Sigma^\blank$, sums and products.
  The structure of this subcategory is simple enough to be axiomatised
  directly in terms of overt discrete objects,
  \ie from an arithmetic universe \cite{TaylorP:insema}.
  On the other hand, the whole category $\S$ may be recovered from it
  by the monadic construction in~\cite{TaylorP:subasd}.

\item In an \textdf{L-domain}, the conditional joins in (c) are replaced
  by joins within the lower subset $\tridown n$.
  
\item If you are familiar with the characterisation of SFP or bifinite
  domains, and their continuous generalisations,
  you will be able to see how to encode them and their bases in ASD.
  Maybe this could provide a way of making the results about them
  more constructive, and identifying the largest cartesian closed
  category of locally compact objects \cite{JungA:clacd}.
\qed
\end{letterlist}
\end{r@Examples}

\medskip

In the classical theory, a domain is recovered from its basis as
the (continuous or algebraic) dcpo of (rounded) ideals,
with the Scott topology.
We can prove the corresponding result within ASD.

\begin{r@Lemma}
 For any directed interpolative relation $(N,{\trless})$,
$$ \Phi:\Sigma^{\Sigma^N},\; \xi:\Sigma^N \ \proves\
   \E \Phi\xi \;\equiv\; \Some n. \Phi(\tridown n)\land\xi n $$
is a nucleus on $\Sigma^N$.
\vadjust{\nobreak}
\end{r@Lemma}

\goodbreak

\Proof We have a  directed%
  \footnote{It is directed in the existential sense of
  Definition~\ref{def dirsup}, not in the canonical one of
  Definition~\ref{DJ} needed for Theorem~\ref{F:SSX cts}.
  An $\natno$--$\natno$ choice principle is needed for open
  (recursively enumerable) relations that amounts to sequentialising
  non-deterministic programs.
  This will be justified and applied to other questions
  arising from this paper in future work \cite{TaylorP:loccbc}.}
union
$$ (\tridown n) \;\eq\; \Lamb m.\Some k.(m\trless k\trless n)
     \;\eq\; \Some k.(\tridown k) \land (k\trless n) $$
by the extrapolation and directedness conditions,
so by Theorem~\ref{F:SSX cts} and the definition of~$\E$,
$$ \Phi(\tridown n) \;\eq\; \Some k.\Phi(\tridown k) \land (k\trless n)
   \;\eq\; \E \Phi(\tridown n). $$
Hence we deduce the equations for a nucleus
in Theorem~\ref{lattice nucleus},
\begin{eqnarray*}
  \E(\Phi\sqr \Psi)\xi
    &\equiv& \Some n. \big(\Phi(\tridown n)\sqr \Psi(\tridown n)\big)\land\xi n \\
    &\equiv& 
       \Some n. \big(\E \Phi(\tridown n)\sqr\E \Psi(\tridown n)\big)\land\xi n \\
    &\equiv& \E(\E \Phi\sqr\E \Psi)\xi. &$\qeds$%
\end{eqnarray*}

\begin{r@Definition}
\label{def admit}
We have defined $\E$ from entirely abstract data
(a directed interpolative relation)
but have shown (\via Theorem~\ref{lattice nucleus})
that it satisfies the characteristic property of the composite
$\Sigma^i\cdot I$ that arises from a $\Sigma$-split subobject.

Now, it was exactly the purpose of the monadic calculus 
\cite[\S 8]{TaylorP:subasd} that, in this situation,
we do have an object, called $X\equiv\collect{\Sigma^N}{\E}$ there,
that has such a $\Sigma$-split inclusion.
That is, the calculus provides $i:\collect{\Sigma^N}{\E}\splitinto\Sigma^N$
and $I$ \emph{formally}.
Then $\Gamma\proves\xi:\Sigma^N$ is \textdf{admissible},
\ie an element of the new subobject, if
$$ \Gamma,\; \Phi:\Sigma^{\Sigma^N}\ \proves\ \E \Phi\xi \eq \Phi\xi. $$
The idea of this is that $\E$ ``normalises'' open subobjects $\Phi$
of $\Sigma^N$ by restricting and re-expanding them,
and it is only the elements of the subobject being defined
whose membership of \emph{all} such $\Phi$ remains unchanged in this process.

In this case, the calculus allows us to introduce $x\equiv\admit\xi:X$
such that $\xi=i x$.
\end{r@Definition}

\goodbreak

\begin{r@Lemma}
\label{dcpo admit}
  $\xi:\Sigma^N$ is \emph{admissible} iff it is a rounded ideal.
\end{r@Lemma}

\Proof
Putting $\Phi\equiv\Lamb\theta.\theta n$, $\Phi\equiv\Lamb\theta.\top$
and $\Phi\equiv\Lamb\theta.\theta r\land\theta s$,
we deduce that $\xi$ is a rounded ideal.
Conversely, if $\xi$ is a rounded ideal then
$$ \All m\in\ell.\xi m \;\eq\; \Some k.(\All m\in\ell.m\trless k)\land\xi k $$
by equational induction on $\ell$.
Then, using the prime $\land$-basis on $\Sigma^N$,
\begin{eqnarray*}
  \Phi\xi &\eq& \Some\ell.\Phi(\Lamb m.m\in\ell)\land \All m\in\ell.\xi m \\
  \Phi(\tridown k) &\eq&
       \Some\ell.\Phi(\Lamb m.m\in\ell)\land \All m\in\ell. m\trless k \\
  \Phi\xi &\eq& \Some\ell.\Phi(\Lamb m.m\in\ell)\land
       \Some k. (\All m\in\ell.m\trless k)\land\xi k \\
  &\eq& \Some k.\Phi(\tridown k)\land\xi k \;\equiv\; \E\Phi\xi. &$\qeds$%
\end{eqnarray*}

\vskip-.5\baselineskip

\begin{r@Corollary}
$\tridown n$ is admissible, because it is a rounded ideal. \qed
\end{r@Corollary}

\vskip-.5\baselineskip

\begin{r@Theorem}
\label{prime basis} Any directed interpolative relation
$(N,{\trless})$ is the opposite of the way-below relation for a prime
basis on its object of rounded ideals.
\end{r@Theorem}

\Proof Let $X\equiv\collect{\Sigma^N}{\E}\rSplitinto^i\Sigma^N$
in the notation of \cite[Section 8]{TaylorP:subasd}.
Now I claim that
$$ p_n \;\equiv\; \admit(\tridown n)
   \quad\hbox{and}\quad
   \beta^n  \;\equiv\; \Lamb x. i x n \;\equiv\; \Sigma^i(\Lamb\xi.\xi n)$$
define a prime basis.
Since $\xi\equiv i x$ is admissible,
we have the basis expansion,
\begin{eqnarray*}
 \phi x \;\eq\; I\phi(i x)  \;\eq\; \E(I\phi)(i x) 
  &\eq& \Some n.I\phi(\tridown n)\land i x n \\
  &\eq& \Some n.I\phi(i p_k) \land \beta^n x \\
  &\eq& \Some n.\phi p_k \land \beta^n x.
\end{eqnarray*}
Finally,
$(n\waybelow m) \;\eq\; A_n \beta^m \;\eq\; \beta^m p_n 
   \;\eq\; i\big(\admit(\tridown n)\big) m
   \;\eq\; (\tridown n)m \;\eq\; (m\trless n)$. \qed

\medskip
   
We leave domain theory aside there and return to topology.

\begin{r@Remark}
Jung and S\"underhauf characterise \emph{stably} locally compact spaces
using a system of rules that they call a \textdf{strong proximity lattice}
\cite[Section 5]{JungA:duacvo}.
These axioms are very similar to those in Section~\ref{way-below},
except that theirs are lattice dual, whereas ours are not.
In particular, they prove the dual Wilker property
(Corollary~\ref{classical dual Wilker}), albeit using Choice.
If $(N,0,1,{+},{\star},{\waybelow})$ is an abstract basis satisfying
these axioms then so too is $(N,1,0,{\star},{+},{\wayabove})$.
The corresponding object, which is known as the \textdf{Lawson dual},
classically has the same points as the given one,
but the opposite specialisation order,
whilst the open subobjects of one correspond to
the compact saturated subobjects of the other.
\end{r@Remark}

\begin{r@Example}
\label{upd}
Let $\waybelow$ be the way-below relation for a filter lattice basis
on a stably locally compact space $X$.
The object $\powerset^\sharp X$ of rounded ideals for $\wayabove$
(Example~\ref{upd basis})
is called the \textdf{Smyth} or \textdf{upper powerdomain} of~$X$
and is Lawson dual to the topology $\Sigma^X$. \qed
\end{r@Example}
\bigbreak

Lawson duality goes way beyond the purposes of this paper,
but we can achieve the dual Wilker property by defining a new basis.

\begin{r@Proposition}
\label{ASD dual Wilker} Any stably locally compact object,
\ie one that has a filter $\land$-basis $(N,0,1,{+},{\star},{\waybelow})$,
has another such basis indexed by $\Fin(N)$ that also
satisfies the dual Wilker property,
$$
\begin{prooftree}
  p\star q\waybelow n
  \Justifies
  \Some p' q'.
  (p'\star q'\waybelow n)\land(p\waybelow p')\land(q\waybelow q').
\end{prooftree}
$$
This can be extended to a new filter lattice basis with the same property
by Lemma~\ref{make or basis}ff.
\end{r@Proposition}

\Proof We define a new version of $\star$, called $\times$, by
$$ A_{p\times q} \;\equiv\;
   \Some p' q'.A_{p'\star q'}\land(p\waybelow p')\land(q\waybelow q').
$$
This construction, unlike Remark~\ref{make and basis}, preserves the filter
property, and is also idempotent.
\begin{eqnarray*}
  A_{p\times q}\top
  &\eq& \Some p' q'.A_{p'\star q'}\top\land(p\waybelow p')\land(q\waybelow q') \\
  &\Impliedby& (p\waybelow 1)\land(q\waybelow 1) \;\eq\; \top\\
  A_{p\times q}
  &=&\Some p'' q''.A_{p''\star q''}\land(p\waybelow p'')\land(q\waybelow q'')
     &def $\times$\\
  &=&\Some p' q' p'' q''.A_{p''\star q''}\land
     (p\waybelow p'\waybelow p'')\land(q\waybelow q'\waybelow q'')
     &interpolation\\
  &=& \Some p' q'.A_{p'\times q'}\land(p\waybelow p')\land(q\waybelow q')
     &def $\times$%
\end{eqnarray*}
We deduce the dual Wilker property by applying the last equation to $\beta^n$.
Now, as the $A_n$ preserve $\land$,
\begin{eqnarray*}
  A_{p\times q}(\phi\land\psi)
  &\eq& \Some p''' q'''.A_{p'''\star q'''}(\phi\land\psi)
      \land(p\waybelow p''')\land(q\waybelow q''')\\
  &\eq& \Some p''' q'''.A_{p'''\star q'''}\phi\land A_{p'''\star q'''}\psi
      \land(p\waybelow p''')\land(q\waybelow q''')\\
  \Implies\;
  A_{p\times q}\phi\land A_{p\times q}\psi
  &\eq& \Some p' q' p'' q''.A_{p'\star q'}\phi\land A_{p''\star q''}\psi \\
  &&\quad \land(p\waybelow p')\land(p\waybelow p'')
           \land(q\waybelow q')\land(q\waybelow q'').
\end{eqnarray*}
For the converse, we put
$p'''\equiv p'\star p''$ and $q'''\equiv q'\star q''$ and need
$$ A_{p'\star q'}\phi \;\Implies\; A_{p'\star p''\star q'\star q''}\phi
   \quad\hbox{and}\quad
   (p\waybelow p')\land(p\waybelow p'') \;\Implies\; (p\waybelow p'\star p''),
$$
which we have from contravariance of $A_p$ and the fact that we had an
$\land$-basis.
The new operation $\times$ can be defined for more factors in the same way,
using recursion and equational induction.
\qed

\medskip

General locally compact objects are a good deal more complicated than domains.
The extension of the earlier results of this section rely heavily on
the ``lattice'' structure $(0,1,{+},{\star})$, which we investigate next.

%============================================================================
\ifpdf\section{The lattice basis on Sigma N}
\else \section{The lattice basis on $\Sigma^N$}
\fi
\label{DLmonad}

In the next section we shall show that the properties of $\waybelow$
listed in Section~\ref{way-below} are sufficient to reconstruct the object
$X$ and its basis.
To do this, however, we need some more technical information
about the lattice basis on $\Sigma^N$,
and about the free distributive lattice on $N$.

\begin{r@Remark}
 When an object $X$ has an $N$-indexed basis $(\beta^n,A_n)$
there is an embedding $X\rSplitinto\Sigma^N$,
given by Lemma~\ref{base->subsp}.
This structure may be summed up by the diagram,
\begin{diagram}[height=1.5em,width=5em]
  \Fin(\Fin N) & \rTo^{B^\blank} & \Sigma^{\Sigma^N} & \rTo^\E & 
     \Sigma^{\Sigma^N} \\
  & \ruSplitinto(2,4)^{\Lamb\xi.\xi\blank} && \ruSplitinto(2,4)^I \\
  \uTo<{\listof{\listof\blank}} && \dOnto<{\Sigma^i} && \uSplitinto>J \\
  \\
  N & \rOnto^{\beta^\blank} & \Sigma^X & \rSplitinto^{\Lamb n.A_n\blank}
     & \Sigma^N \\
  \uLine &&&& \uTo \\
  \HmeetV && \rLine^{m\mapsto\ddown m\equiv(\Lamb n.n\waybelow m)} && \HmeetV
\end{diagram}
in which $\listof{-}$ means the singleton \emph{list},
the map $J$ on the right is defined by
$$ J\phi \;\equiv\; \Lamb\xi.\Some n.\phi n\land \xi n, $$
and $J\phi$ preserves $\lor$, $\bot$ and $\exists$.
We define $\E$ from $\waybelow$ using the lattice basis $(B^L,\A_L)$
on $\Sigma^N$.
\end{r@Remark}

\begin{r@Remark}
\label{DL presentation}
These results may be seen as \textdf{presentation} of the algebra $\Sigma^X$.
In this, $N$ is the set of generators and
$\DL(N)$ is the free (imposed) distributive lattice on them,
which we construct in Proposition~\ref{DL cong} below.
Then there are homomorphisms
\begin{diagram}[small]
  \Fin(\Fin N) && \rOnto^{B^\blank} &&
  \Sigma^{\Sigma^N} && \rOnto^{\Sigma^i} && \Sigma^X \\
  & \rdOnto && \ruTo \\
  && \DL(N)
\end{diagram}
as $B^\blank$ takes the imposed list operations $0$, $1$, $+$ and $\star$
to the intrinsic structure $\bot$, $\top$, $\lor$ and $\land$
in $\Sigma^{\Sigma^N}$,
whilst $\Sigma^i$ preserves the latter in~$\Sigma^X$.
The relation $\waybelow$ encodes the system of ``equations''
that distinguishes the particular algebra $\Sigma^X$
from the generic one $\Sigma^{\Sigma^N}$ that is freely generated by~$N$.

This explains why $n\waybelow m$ does not imply $n\baseleq m$ in our
system, whereas $\alpha\ll\beta$ implies $\alpha\leq\beta$ in a
continuous lattice.
Topologically, we already saw the point in
Remarks \ref{basis dist latt} and~\ref{beta not inj}:
many distinct codes may in principle represent the same open or compact
subobject.
To put this the other way round, since equality (or containment) of
open subobjects is not computable, we cannot deduce equality (or
comparison) of codes from semantic coincidence of subobjects.
\end{r@Remark}

\medbreak

\begin{r@Remark}
\label{DL triangle}
 The triangle
\begin{diagram}
  \DL(N) && \rTo^{B^\blank} && \Sigma^{\Sigma^N} \\
  & \luTo<{\listof{\listof-}} && \ruTo>{\Lamb\xi.\xi\blank} \\
  && N
\end{diagram}
illustrates the comparison between the monad
that captures the imposed $(0,1,{+},{\star})$ distributive lattice structure
and the one in Axiom~\ref{monad axiom}
based on $\Sigma^\blank\adjoint\Sigma^\blank$.
The upward maps are the units of these monads.
We leave the interested student to construct the multiplication map
of the $\DL$-monad as a list program,
\cf $\FLATTEN$ in Lemma~\ref{make lattice basis}.

We are not quite justified in saying that the $\Sigma^{\Sigma^\blank}\!$
monad defines the intrinsic $(\bot,\top,{\lor},{\land})$ distributive
lattice structure.
Corollary~\ref{tfae hom} said that the homomorphisms are the same,
but I~have not been able to show that every object whose intrinsic
order is that of a distributive lattice is an algebra for the monad,
\ie of the form $\Sigma^X$ for some object~$X$.
\end{r@Remark}

\medbreak

\begin{r@Notation}
 Recall from Proposition~\ref{N bases} that the basis on $\Sigma^N$ is
$$ B^L \xi \;\equiv\; \Some\ell\in L.\All m\in\ell.\xi m
   \quad\hbox{and}\quad
   \A_L \Phi \;\equiv\; \All\ell\in L.\Phi(\Lamb m.m\in\ell),
$$
from which we obtain the way-below relation
$$ \A_L B^R \;\equiv\; (R\shp\subset L) \;\equiv\;
   \All\ell\in L.\Some\ell'\in R.(\ell'\subset\ell).
   \vadjust{\nobreak}
$$
The list of lists $R+S$ is given by concatenation, whilst $R\star S$
and $B^{R\star S}$ were defined in Lemma~\ref{make lattice basis}.
\end{r@Notation}

\goodbreak

\begin{r@Lemma}
\label{B+*} $(B^L,\A_L)$ is a lattice basis:
$$B^1=\top\quad B_0=\bot\quad B^{R\star S}=B^R\land B^S
    \quad B^{R+S}=B^R\lor B^S$$
$$ \A_0 \Phi\eq\top\quad \A_L\bot\eq(L=0) \quad \A_{R+S}=\A_R\land\A_S. 
   \eqno\qEd
$$
\end{r@Lemma}

\begin{r@Lemma}
\label{D*} $(B^L,\A_L)$ is a filter basis:
$A\top\eq\top$ and 
$ \A_L(\Phi\land \Psi) \eq \A_L \Phi \land \A_L \Psi$. 
\end{r@Lemma}

\Proof\closeupaline
\begin{eqnarray*}
  \A_L(\Phi\land \Psi)
  &\eq& \All\ell\in L.(\Phi\land \Psi)(\Lamb m.m\in\ell) \\
  &\eq& \All\ell\in L.\Phi(\Lamb m.m\in\ell)\land \Psi(\Lamb m.m\in\ell) \\
  &\eq& \big(\All\ell\in L.\Phi(\Lamb m.m\in\ell)\big)\ \land\
      \big(\All\ell\in L.\Psi(\Lamb m.m\in\ell)\big)\\
  &\eq& \A_L \Phi\land \A_L \Psi &$\qeds$%
\end{eqnarray*}

\smallskip

The Wilker condition says that we can split the list into the two
parts that satisfy the respective disjuncts.

\begin{r@Lemma}
\label{split}\quad
$\A_L(\Phi\lor \Psi) \;\eq\;
 \Some L_1 L_2.(L=L_1+L_2) \land \A_{L_1} \Phi \land \A_{L_2} \Psi$.
\qed
\end{r@Lemma}

\begin{r@Proposition}
\label{DL cong}
 We write $L\isomo R$ if both $R\shp\subset L$ and
$L\shp\subset R$.  This is an open congruence for the imposed
structure on $\Fin(\Fin N)$, and the free imposed distributive lattice
$\DL(N)$ is its quotient \cite[Section~10]{TaylorP:geohol}. \qed
\end{r@Proposition}

\bigskip

\begin{r@Remark}
 So far we have not used any of the structure on $N$ itself.
Since we have a lattice basis for $X$, by definition
$$ \beta^\blank : N \to \Sigma^X $$
takes the imposed structure $(0,1,{+},{\star})$ on $N$ to the
intrinsic structure $(\bot,\top,{\lor},{\land})$ on $\Sigma^X$.
Associated with this imposed structure is an imposed order relation
$\baseleq$, which $\beta^\blank$ takes to $\leq$,
but with respect to which the dual basis $A_\blank$ is contravariant.
\end{r@Remark}

\begin{r@Definition}
\label{leq from +*} We define $\baseleq$ from $+$ and
$\star$ as the least relation such that
$$ 0\baseleq n\baseleq n\baseleq1 \qquad
    (k\star n)+(k\star m) \baseleq k*(n+m)$$
$$
\begin{prooftree}
    n\baseleq k\baseleq m
    \justifies
    n\baseleq m
  \end{prooftree}
  \qquad
  \begin{prooftree}
    k\baseleq m \quad k\baseleq n
    \Justifies
    k\baseleq m\star n
  \end{prooftree}
  \qquad
  \begin{prooftree}
    n\baseleq k \quad m\baseleq k
    \Justifies
    n+m\baseleq k
  \end{prooftree}
$$
and again we write $n\cong m$ when both $n\baseleq m$ and $m\baseleq n$.
  \end{r@Definition}

\begin{r@Proposition}
 The relation $\cong$ is an open congruence on $N$ whose
quotient is an imposed distributive lattice. \qed
\end{r@Proposition}

\medskip

\begin{r@Notation}
\label{evL} Returning to $\Fin(\Fin N)$, $L$ is regarded as a
\emph{formal} sum of products of elements of $N$ (additive normal form).
This may be ``evaluated'' by means of the operation
$$ \ev:\Fin(\Fin N)\to N. $$
This is defined for lists by a generalisation of
Lemma~\ref{make lattice basis}.
A~similar construction works for Kuratowski-finite subsets instead,
except that then $N$ has actually to satisfy the equations for
a distributive lattice up to equality, and not just up to~$\cong$
(\cf Notation~\ref{N notn}).
The map
$$ \DL(N)\equiv\Fin(\Fin N)/({\cong})
   \rTo^{\quad\ev/({\cong})\quad}
   N/({\cong}) $$
is the structure map of the distributive lattice,
regarded as an algebra for the $\DL$-monad.
\end{r@Notation}

\begin{r@Proposition}
\label{ev +x hom}
The map $\ev:\Fin(\Fin N)\to N$ is a homomorphism in the sense that
$\ev 0=0$ and $\ev 1=1$ by construction, whilst
$$ \ev(R+S) \;\cong\; (\ev R) + (\ev S)
  \quad\hbox{and}\quad
  \ev(R\star S) \;\cong\; (\ev R) \star (\ev S).
$$
\end{r@Proposition}

\Proof This is a standard piece of universal algebra,
which again we leave as a student exercise.
The $+$ equation is proved by list induction,
using associativity and commutativity of $+$ up to $\isomo$.
The equation for $\star$ is more difficult,
as we have to take apart the inner lists, and use distributivity. \qed

\begin{r@Lemma}
\label{DL+*} $\A_L B^L \eq \top$ but
$$
\begin{array}[t]{rclcl}
  \A_L B^R \;\eq\; (R\shp\subset L)
   &\eq& \All\ell\in L.\Some\ell'\in R.(\ell'\subset\ell)
  &\Implies& (\ev L\baseleq\ev R) \\
  \A_L(\Lamb\xi.\xi n) &\eq&
  \All\ell\in L.n\in\ell &\Implies& (\ev L\baseleq n) \\
  \A_L(\Lamb\xi.\xi n\land\xi m) &\eq&
  \All\ell\in L.n\in\ell\land m\in\ell &\Implies& (\ev L\baseleq n\star m) \\
  \A_L(\Lamb\xi.\xi n\lor\xi m) &\eq&
  \All\ell\in L.n\in\ell\lor m\in\ell &\Implies& (\ev L\baseleq n+m) \\
\end{array}
$$
with equality in the cases $L\equiv\listof{\listof k}$ and $L\equiv R$.
\end{r@Lemma}

\Proof In the expansion of $\A_L B^R$, the products in $L$ are of
longer strings than those in~$R$. The other three results follow by
putting $R\equiv\listof{\listof n}$, $\listof{\listof{n,m}}$ and
$\listof{\listof n,\listof m}$, \qed

\medskip

\begin{r@Remark}
 The foregoing discussion of $\cong$ is the price that we pay
for not requiring $(N,0,1,{+},{\star})$ to satisfy the equations for
a distributive lattice in Remark~\ref{basis dist latt}.
If, like \cite{JungA:duacvo}, we had done so,
we would have instead paid the same price
to construct the basis for $\Sigma^X$.
This is indexed by the free distributive lattice
on $\opp N$ \qua $\star$-semilattice,
\ie with new joins but using the old ones as meets.
What we have bought for this price is the flexibility of
switching amongst various kinds of basis (Definition~\ref{basis jargon}),
when we needed or \emph{didn't} need them, throughout the paper.

The reason why it is unnecessary to form the quotient of $N$ or
$\Fin(\Fin N)$ by the congruence $\cong$ is that we never deal with
their elements up to equality. The things that matter are the rules
$$
\begin{prooftree}
  n\isomo n' 
  \qquad
  n \waybelow k
  \Justifies 
  n' \waybelow k
\end{prooftree}
\qquad\qquad
\begin{prooftree}
  k \waybelow n
  \qquad
  n\isomo n'
  \Justifies 
  k \waybelow n'
\end{prooftree}
$$
which are examples of Corollary~\ref{<<-leq}.
Indeed the relation $\baseleq$ itself is only needed to avoid
the extra rules that relate $+$ and $\star$ to $\waybelow$,
which appear in \cite[Lemma~7]{JungA:duacvo}.
\end{r@Remark}

%============================================================================
\section{From the basis to the space}\label{X from <<}

We are now able to show that any ``abstract'' basis satisfying
the conditions of Section~\ref{way-below} actually arises from some
definable object.
As in Section~\ref{domain theory},
we have to define a nucleus $\E:\Sigma^{\Sigma^N}\to\Sigma^{\Sigma^N}$
and (in the next section) characterise its admissible terms
in the sense of Definition~\ref{def admit}.

\begin{r@Definition}
\label{abs basis} An \textdf{abstract basis} is an overt
discrete object $N$ with elements $0,1\in N$, binary operations
${+},{\star}:N\times N\to N$ and an open binary relation
${\waybelow}:N\times N\to\Sigma$ such that
$$
  0\waybelow 0\qquad
\begin{prooftree}
    n\waybelow p\quad m\waybelow p
    \Justifies
    n+m\waybelow p
  \end{prooftree}
  \qquad
  \begin{prooftree}
    m'\baseleq m \quad m\waybelow n \quad n\baseleq n'
    \justifies
    m'\waybelow n'
  \end{prooftree}
  \qquad
  \begin{prooftree}
    n\waybelow m
    \Justifies
    n\waybelow k\waybelow m
  \end{prooftree}
$$
$$
\begin{prooftree}
  n\waybelow m \quad m\waybelow p \quad m\waybelow q
  \justifies
  n\waybelow p\star q
\end{prooftree}
\qquad
\begin{prooftree}
  n\waybelow p+q
  \Justifies
  n\waybelow p'+q' \qquad  p'\waybelow p \qquad  q'\waybelow q
\end{prooftree}
$$
where $\baseleq$ is defined from $+$ and $\star$ by
Definition~\ref{leq from +*}.
  \end{r@Definition}

\begin{r@Notation}
\label{E from <<}
\leavevmode\vadjust{\vskip-\baselineskip}
\begin{eqnarray*}
  \E &\equiv&
\Lamb \Phi.\Lamb\xi.J\big(\Lamb n.\Some L.(n\waybelow\ev L)\land\A_L \Phi\big)\xi\\
 &\equiv&
 \Lamb \Phi.\Lamb\xi.\Some n.\Some L.\xi n\land(n\waybelow\ev L)\land\A_L \Phi.
\end{eqnarray*}
\end{r@Notation}

\noindent
We have to show that this satisfies the equations in
Theorem~\ref{lattice nucleus},
$$ \Phi,\Psi:\Sigma^{\Sigma^N} \ \proves\
   \E(\Phi\land \Psi)\;=\;\E(\E \Phi\land \E \Psi)
   \quad\hbox{and}\quad
   \E(\Phi\lor \Psi)\;=\;\E(\E \Phi\lor \E \Psi).
$$
This is made a little easier by the fact that we need only test these
equations for \emph{basic} $\Phi\equiv B^R$ and $\Psi\equiv B^S$:

% ---------------------------------------------------------------------

\begin{r@Lemma}
\label{FCRGCS}
\begin{prooftree}
  R,S:\Fin(\Fin(N)\ \proves\ \E(\E B^R\sqr\E B^S)=\E(B^R\sqr B^S)
  \justifies
  \Phi,\Psi:\Sigma^{\Sigma^N}\ \proves\ \E(\E\Phi\sqr\E\Psi)=\E(\Phi\sqr\Psi).
\end{prooftree}
\end{r@Lemma}

\Proof We use the lattice basis expansion $\Phi = \Some R.\A_R \Phi\land B^R$.

Note first that the combined expansion using distributivity
(Lemma~\ref{*=<<**}),
$$ \Phi \sqr \Psi \;=\; \Some R S. \A_R \Phi \land \A_S \Psi \land (B^R\sqr B^S), $$
is directed in $<R,S>$, so
$\E$ preserves the join by Theorem~\ref{F:SSX cts}
and $\E \Phi \;=\; \Some R.\A_R \Phi \land \E B^R$.
Using distributivity, directedness and Scott continuity again, we have
\begin{eqnarray*}
  \E \Phi\sqr\E \Psi
   &=& \Some R S. \A_R \Phi \land \A_S \Psi \land (\E B^R\sqr \E B^S) 
   &distributivity\\
  \E(\E \Phi\sqr\E \Psi)
   &=& \Some R S. \A_R \Phi \land \A_S \Psi \land \E(\E B^R\sqr \E B^S)
   &directedness\\
   &=& \Some R S. \A_R \Phi \land \A_S \Psi \land \E(B^R\sqr B^S) &hypothesis\\
   &=& \E(\Phi\sqr \Psi) &$\qeds$%
\end{eqnarray*}

% ---------------------------------------------------------------------

\noindent
Next we need to evaluate the expression $\A_L(\E B^R)$.

\begin{r@Lemma}
\label{ECR}\quad $\E B^R \;=\; J(\Lamb n.n\waybelow\ev R)$.
\end{r@Lemma}

\Proof \closeupaline
\begin{eqnarray*}
  \E B^R
  &=& J(\Lamb n.\Some L.n\waybelow\ev L\land \A_L B^R)
     &Definition \ref{E from <<}\\
  &\leq& J(\Lamb n.\Some L.n\waybelow\ev L\baseleq\ev R) &Lemma~\ref{DL+*}\\
  &=& J(\Lamb n.n\waybelow\ev R) &Definition~\ref{abs basis}%
\end{eqnarray*}
but the $\leq$ is an equality,
as we may put $L\equiv R$ in the other direction. \qed

\begin{r@Lemma}
\label{dlecr} Since, by hypothesis,
$\waybelow$ satisfies $0\waybelow k$ and (\cf Lemma~\ref{0+ll})
$$
\begin{prooftree}
  n \waybelow k \qquad m \waybelow k 
  \Justifies
  n+m \waybelow k
\end{prooftree}
\qquad\qquad
\begin{prooftree}
  m\star n\baseleq n \qquad n \waybelow r
  \justifies
  m\star n\waybelow r
\end{prooftree}
$$
we have $\A_L(\E B^R) \;\Implies\; (\ev L\waybelow\ev R)$,
with equality in the case $L\equiv\listof{\listof k}$.
\end{r@Lemma}

\Proof The reason for the inequality is that, whereas
$B^\blank:\Fin(\Fin N)\to\DL(N)\to\Sigma^{\Sigma^N}$ sends $+$ to
$\lor$ and $\star$ to $\land$,\quad $\A_\blank$ only takes $+$ to $\land$,
not necessarily $\star$ to $\lor$.
\begin{eqnarray*}
  \A_L(\E B^R)
  &\eq& \All\ell\in L.\E B^R(\Lamb m.m\in\ell) &def $\A_L$\\
  &\eq& \All\ell\in L.J(\Lamb n.n\waybelow\ev R)(\Lamb m.m\in\ell)
     &Lemma~\ref{ECR}\\
  &\eq& \All\ell\in L.\Some n\in\ell.n\waybelow\ev R &def.~$J$\\
  &\Implies& \All\ell\in L.\mu\ell\waybelow\ev R &$\star$ rule\\
  &\eq& \ev L\waybelow\ev R &$0$, $+$ rules%
\end{eqnarray*}
Here $\mu\ell$ is the ``product'' of $\ell$, in the sense of $1$ and $\star$
(written $\FOLD\,{\star}\,1\,\ell$ in functional programming notation),
so $n\in\ell\Implies\mu\ell\baseleq n$.
Then $\ev L$ is the sum of these products (Definition~\ref{leq from +*}).
Equality holds when $L\equiv\listof{\listof k}$ since $\ell=\listof k$
and $\ev L=\mu\ell=k$. \qed

% ---------------------------------------------------------------------
\bigskip

Equipped with formulae for $\E B^R$ and $\A_L(\E B^R)$,
we can now verify the two equations. 
Their proofs are \emph{almost} the same,
illustrating once again the lattice duality that we get
by putting directed joins into the background.
Unfortunately, they're not quite close enough for us to use $\sqr$
and give just one proof. First, however, we give the similar but
slightly simpler argument for idempotence, although it is easily seen
to be implied by either of the other results.

\begin{r@Proposition}
 If $\waybelow$ satisfies the transitive and interpolation rules,
$$
\begin{prooftree}
  n\waybelow r
  \Justifies
  n\waybelow m\waybelow r
\end{prooftree}
$$
\cf Corollary~\ref{interpolation},
then $\E$ is idempotent: $\E(\E \Phi)=\E \Phi$.
\end{r@Proposition}

\Proof By (a simpler version of) Lemma~\ref{FCRGCS},
it's enough to consider $\Phi\equiv B^R$,
\begin{eqnarray*}
  \E(\E B^R)
  &=& J\big(\Lamb n.\Some L.(n\waybelow\ev L)\land \A_L(\E B^R)\big)
     &Notation \ref{E from <<}\\
  &\leq& J\big(\Lamb n.\Some m.(n\waybelow m)\land (m\waybelow\ev R)\big)
    &Lemma~\ref{dlecr}\\
  &=& J(\Lamb n.n\waybelow\ev R) &hypothesis\\
  &=& \E B^R &Lemma~\ref{ECR}%
\end{eqnarray*}
where $m\equiv\ev L$.
However, the $\leq$ is an equality as we may use
$L\equiv\listof{\listof m}$ in Lemma~\ref{ECR} to prove $\geq$.  \qed

\begin{r@Proposition}
 Since $\waybelow$ obeys the rule linking it with $\star$,
$$
\begin{prooftree}
  n\waybelow r\star s
  \Justifies
  n\waybelow m \quad m\waybelow r \quad m\waybelow s
\end{prooftree}
$$
(\cf Lemma~\ref{<<* rule})
$\E$ satisfies the $\land$-equation,
$ \E(\Phi\land \Psi) \;=\; \E(\E \Phi\land \E \Psi)$.
\end{r@Proposition}

\Proof By Lemma~\ref{FCRGCS},
it's enough to consider $\Phi\equiv B^R$ and $\Psi\equiv B^S$.
With $m\equiv\ev L$, $r\equiv\ev R$ and $s\equiv\ev S$,
\begin{eqnarray*}
  \qquad\llap{$\E(\E B^R\land \E B^S)$}
  &=& J\big(\Lamb n.\Some L.(n\waybelow\ev L)\land\A_L(\E B^R\land \E B^S)\big)
     &Notation \ref{E from <<}\\
  &=&
 J\big(\Lamb n.\Some L.(n\waybelow\ev L)\land\A_L(\E B^R)\land\A_L(\E B^S)\big)
     &Lemma~\ref{D*}\\
  &\leq&
J\big(\Lamb n.\Some m.(n\waybelow m)\land(m\waybelow r)\land(m\waybelow s)\big)
    &Lemma~\ref{dlecr}\\
  &=& J(\Lamb n.n\waybelow r\star s) &hypothesis\\
  &=& J\big(\Lamb n.n\waybelow\ev(R\star S)\big) &Lemma~\ref{ev +x hom}\\
  &=& \E B^{R\star S} \;=\; \E(B^R\land B^S) &Lemmas~\ref{ECR} \& \ref{B+*}%
\end{eqnarray*}
To make the $\leq$ an equality,
we use $L\equiv\listof{\listof m}$ in Lemma~\ref{DL+*}.
\qed

\begin{r@Proposition}
 If $\waybelow$ satisfies the Wilker rule linking it with $+$,
$$
\begin{prooftree}
  n\waybelow r+s
  \Justifies
  n\waybelow p+q \quad p\waybelow r \quad q\waybelow s
\end{prooftree}
$$
(\cf Lemma~\ref{basis Wilker})
then $\E$ satisfies the $\lor$-equation,
$ \E(\Phi\lor \Psi) \;=\; \E(\E \Phi\lor \E \Psi)$. 
\end{r@Proposition}

\Proof With
$r\equiv\ev R$, $s\equiv\ev S$, $p\equiv\ev L_1$, $q\equiv\ev L_2$
and $L=L_1+L_2$,
\begin{eqnarray*}
  \E(\E B^R\lor \E B^S)
  &=& J\big(\Lamb n.\Some L.(n\waybelow\ev L)\land \A_L(\E B^R\lor \E B^S)\big)
     &Notation \ref{E from <<}\\
  &=& J\big(\Lamb n.\Some L_1 L_2.
  (n\waybelow\ev L_1+\ev L_2)\land \A_{L_1}(\E B^R)\land \A_{L_2}(\E B^S)\big)
     &\ref{split}\\
  &\leq& J\big(\Lamb n.\Some p q.
      (n\waybelow p+q)\land (p\waybelow r)\land (q\waybelow s)\big)
     &Lemma \ref{dlecr}\\
  &=& J(\Lamb n.n\waybelow r+s) &hypothesis\\
  &=& J\big(\Lamb n.n\waybelow\ev(R+S)\big) &Lemma~\ref{ev +x hom}\\
  &=& \E B^{R+S} \;=\; \E(B^R\lor B^S) &Lemmas \ref{ECR} \& \ref{B+*}%
\end{eqnarray*}
Again the $\leq$ becomes an equality, using
$L_1\equiv\listof{\listof p}$, $L_2\equiv\listof{\listof q}$
and $L\equiv L_1+L_2$ in Lemma~\ref{DL+*}. \qed

\begin{r@Theorem}
 $\E$ is a nucleus on $\Sigma^N$ in the sense of 
\cite[Section 8]{TaylorP:subasd}. \qed
\end{r@Theorem}

%============================================================================
\section{The points of the new space}\label{new points}

Now we can characterise the (parametric) points of the object
that we constructed from an abstract basis in the previous section.
Then we define a lattice basis on it,
just as we did for continuous dcpos in Theorem~\ref{prime basis}.

\begin{r@Notation}
We are again in the situation of Definition~\ref{def admit},
that we have defined a nucleus $\E$ from an abstract basis,
so the monadic calculus \cite[\S 8]{TaylorP:subasd}
provides a subobject $i:\collect{\Sigma^N}{\E}\splitinto\Sigma^N$
with $\Sigma$-splitting~$I$.
As in Lemma~\ref{dcpo admit}, the next task is to characterise its elements,
\ie those $\Gamma\proves\xi:\Sigma^N$ that are \textdf{admissible}:
$$ \Gamma,\; \Phi:\Sigma^{\Sigma^N} \ \proves\ \Phi\xi \;\eq\; \E \Phi\xi. $$
\end{r@Notation}

\begin{r@Lemma}
\label{lattice rounded}
If $\Gamma\proves\xi:\Sigma^N$ is admissible then it is \textdf{rounded}
in the sense that
$$ \Gamma,\; n:N \ \proves\ \xi n \;\eq\; \Some m.\xi m\land m\waybelow n. $$
\end{r@Lemma}

\Proof Consider $\Phi\equiv\Lamb\xi.\xi n$, so $\A_L \Phi \Implies (\ev L\baseleq n)$
by Lemma~\ref{DL+*}. Then
\begin{eqnarray*}
   \xi n \;\eq\; \Phi\xi \;\eq\; \E \Phi\xi 
   &\eq& \Some m.\Some L.\xi m\land (m\waybelow\ev L)\land \A_L \Phi 
      &Notation \ref{E from <<}\\
   &\Implies& \Some m.\Some L.\xi m\land (m\waybelow\ev L\baseleq n) &above\\
   &\Implies& \Some m.\xi m\land (m\waybelow n) &monotonicity%
\end{eqnarray*}
where $\Implies$ is actually equality,
as we may put $L\equiv\listof{\listof n}$ in Lemma~\ref{DL+*}
to obtain $\Impliedby$. \qed

\begin{r@Lemma}
\label{adm=rhom} If $\Gamma\proves\xi:\Sigma^N$ is admissible
then it is a lattice homomorphism in the sense that
$$ \xi 0 \;\eq\; \bot
   \qquad
   \xi 1 \;\eq\; \top
   \qquad
   \xi(n+m) \;\eq\; \xi n\lor\xi m
   \qquad
   \xi(n\star m) \;\eq\; \xi n\land\xi m.
$$
\end{r@Lemma}

\Proof Consider $\Phi\equiv\Lamb\xi.\xi n\sqr\xi m$,
so $\A_L \Phi \;\Implies\; (\ev L\baseleq n\sqr m)$ by Lemma~\ref{DL+*}.
Then
\begin{eqnarray*}
 \xi n\sqr\xi m
 &\eq&  \Phi\xi \;\eq\; \E \Phi\xi \\
 &\eq& \Some k.\Some L.\xi k\land (k\waybelow\ev L)\land \A_L \Phi
    &Notation \ref{E from <<}\\
 &\Implies& \Some k.\Some L.\xi k\land (k\waybelow\ev L\baseleq n\sqr m)
    &above\\
 &\Implies& \Some k.\xi k\land (k\waybelow n\sqr m)
  \;\eq\; \xi(n\sqr m) &roundedness%
\end{eqnarray*}
with equality by $L\equiv\listof{\listof{n\sqr m}}$ in Lemma~\ref{DL+*}.
Similarly, for the constants,
consider $\Phi\equiv\Lamb\xi.\top$, so $\A_L \Phi\eq\top\eq(\ev L\baseleq 1)$,
and $\Phi\equiv\Lamb\xi.\bot$, so $\A_L \Phi\eq\All\ell\in L.\bot\eq(L=0)$. \qed

\smallskip

Notice that it is the fact that $\xi$ is \emph{rounded} (for $\waybelow$),
rather than a homomorphism (for $0,1,{+},{\star}$),
that distinguishes the particular object $X$ from the ambient $\Sigma^N$
into which it is embedded, \cf Remark~\ref{DL presentation}.

\begin{r@Lemma}
\label{rounded lattice->adm}
If $\Gamma\proves\xi:\Sigma^N$ is a rounded lattice homomorphism
then it is admissible.
\end{r@Lemma}

\Proof \closeupaline
\begin{eqnarray*}
  \Phi\xi &\eq& \Some L.\A_L \Phi \land B^L\xi &lattice basis on $\Sigma^N$\\
  &\eq& \Some L.\A_L \Phi \land\xi(\ev L)\\
  &\eq& \Some L.\Some n.\A_L \Phi \land\xi n\land (n\waybelow\ev L) &rounded\\
  &\eq& \E \Phi\xi &Definition~\ref{E from <<}%
\end{eqnarray*}
where the second step exploits the definition of $B^L$ (Lemma~\ref{B+*})
and $\ev L$ when $\xi$ is a homomorphism. \qed

\begin{r@Remark}
\label{upd rk}
Hence the rounded \emph{lattice} homomorphisms $\xi:\Sigma^N$
are the \emph{points} of $X$,
whilst the rounded \emph{filters} ($\land$-homomorphisms)
correspond to its \emph{compact saturated subspaces},
at least when $X$ is a classical stably locally compact space
(Example~\ref{upd}).
\end{r@Remark}

\begin{r@Lemma}
\label{basis from <<} $\beta^n\;\equiv\;\Lamb x. i x n$ and
$A_n\;\equiv\;\Lamb\phi.I\phi\setof n$ provide an effective basis for $X$.
\end{r@Lemma}

\Proof First note that $\xi\mapsto J\phi\xi$ preserves joins in $\xi$,
and therefore so do $\xi\mapsto\E \Phi\xi$ and $\xi\mapsto I\phi\xi$,
as required by Lemma~\ref{embed->basis}, so we recover
$$ i x \;=\; \Lamb n.\beta^n x
   \quad\hbox{and}\quad
   I\phi \;=\; \Lamb\xi.\Some n.A_n\phi\land\xi n. $$
For the basis expansion, let $\phi:\Sigma^X$ and $x:X$.
Then $\xi\equiv i x:\Sigma^N$ is admissible,
and $\phi=\Sigma^i \Phi$, where $\Phi\equiv I\phi:\Sigma^{\Sigma^N}$.
\begin{eqnarray*}
  \Some n.A_n\phi\land\beta^n x
  &\eq& \Some n.I\phi\setof n\land i x n &defs $\beta^n$, $A_n$\\
  &\eq& \Some n.(I\cdot\Sigma^i) \Phi\setof n\land\xi n &defs $\Phi$, $\xi$\\
  &\eq& \Some n.\E \Phi\setof n\land\xi n &defs $i$, $I$\\
  &\eq& \E \Phi\xi &$\E \Phi$ preserves $\xi\eq\Some n.\setof n\land\xi n$\\
  &\eq& \Phi\xi &$\xi$ admissible\\
  &\eq& (I\phi)(i x) &defs $\Phi$, $\xi$\\
  &\eq& \phi x &$\setof{}\eta$ \cite[Section 8]{TaylorP:subasd}. \qEd
\end{eqnarray*}

\begin{r@Corollary}
  $x \;=\; \focus(\Some n. A_n \land \xi n)$. \qed
\end{r@Corollary}

This is the generalisation of $x=\dirsup\collect{p_n}{\xi n}$
in domain theory (Proposition~\ref{prime basis->dpco}),
and the version in ASD of
$\setof x=\dirinter{K}{\setof x\waybelow K}$ in Theorem~\ref{cpct <<}.
The predicate $\xi n$ means $x\in U^n\subset K^n$ in the spatial setting,
\cf Definition~\ref{pt loc cpct}.

\begin{r@Theorem}
 $(\beta^n,A_n)$ is a lattice basis, and its way-below relation
is $\waybelow$.
\end{r@Theorem}

\Proof Since $x:X\ \proves\ \xi\equiv i x\equiv\Lamb n.\beta^n x:\Sigma^N$
is admissible, and therefore a homomorphism by Lemma~\ref{adm=rhom}, we have
$$ \beta^0 \;=\; \Lamb x.\bot,\quad \beta^1 \;=\; \Lamb x.\top
   \quad\hbox{and}\quad \beta^{n\sqr m} \;=\; \beta^n\sqr\beta^m.
$$
Next we check the equations on $A_n$ for an $\lor$-basis.
Let $\phi:\Sigma^X$ and $\Phi\equiv I\phi:\Sigma^{\Sigma^N}$.
Then $\E \Phi=I\cdot\Sigma^i\cdot I\phi=I\phi=\Phi$, so
$$ A_n\phi\;\eq\;I\phi\setof n \;\eq\; \Phi\setof n\;\eq\;\E \Phi\setof n. $$
Hence \closeupaline
\begin{eqnarray*}
  A_0\phi  &\eq& \E \Phi\setof 0 \\
  &\eq& \Some L.(0\waybelow\ev L)\land\A_L \Phi &Notation \ref{E from <<}\\
  &\Impliedby& (0\waybelow 0)\land\A_0 \Phi \;\eq\; \top \\
  A_n\phi\land A_m\phi 
  &\eq& \E \Phi\setof n\land \E \Phi\setof m \\
  &\eq& \Some L_1 L_2.
 (n\waybelow\ev L_1)\land(m\waybelow\ev L_2)\land\A_{L_1}\Phi\land\A_{L_2}\Phi
     &Def.~\ref{E from <<}\\
  &\Implies& \Some L_1 L_2.
    (n+m\waybelow\ev L_1+\ev L_2)\land\A_{L_1+L_2}\Phi &Lemma~\ref{B+*}\\
  &\Implies& \Some L.(n+m\waybelow\ev L)\land\A_L \Phi &Lemma~\ref{ev +x hom}\\
  &\eq& \E \Phi\setof{n+m} \;\eq\; A_{n+m}\phi &Notation \ref{E from <<}%
\end{eqnarray*}
using distributivity, $L=L_1+L_2$ and Lemma~\ref{D*}.
But $A_{n+m}\phi\Implies A_n\phi$, so we have equality.
Finally,\closeupaline
\begin{eqnarray*}
  A_n\beta^m
  &\eq& I(\Lamb x.i x m)\setof n 
  \;\eq\; (I\cdot\Sigma^i)(\Lamb\xi.\xi m)\setof n\qquad
     &Lemma~\ref{basis from <<}\\
  &\eq& \E(\Lamb\xi.\xi m)\setof n \\
  &\eq& \Some L.(n\waybelow\ev L)\land\A_L(\Lamb\xi.\xi m)
     &Notation \ref{E from <<}\\
  &\Implies& \Some L.(n\waybelow\ev L\baseleq m) & Lemma~\ref{DL+*}\\
  &\Implies& (n\waybelow m)  &monotonicity%
\end{eqnarray*}
where we also obtain $\Impliedby$ by putting $L\equiv\listof{\listof m}$
in Lemma~\ref{DL+*}. \qed

\medskip 

\begin{r@Corollary}
 If the object $X$ and its lattice basis $(\beta^n,A_n)$
had been given, and $\waybelow$ and $\E$ derived from them,
this construction would recover $X$ and $(\beta^n,A_n)$ up to unique
isomorphism.
In particular, if we had started with a filter lattice basis,
we would get it back.
\qed
\end{r@Corollary}

\medskip

We have shown that the notion of ``abstract basis''
is a complete axiomatisation of the way-below relation,
and is therefore the formulation of the consistency requirements
in Definitions \ref{rec loc cpct} and~\ref{rec cts latt},
without using a classically defined topological space or locale as a reference.

This completes the proof of the equivalence amongst the characterisations
of local compactness that we listed in the Introduction
%Remark~\ref{summary loc cpct}.

In Lemmas~\ref{lattice rounded}ff, we have also characterised
points $x:X$ as rounded lattice homomorphisms $\xi:\Sigma^N$.
We shall replace the predicate $\xi m$ by a binary relation $\hayo H^m_n$
in the next section,
in order to generalise from points of $X$ to continuous functions $Y\to X$.

%============================================================================
\section{Morphisms as matrices}\label{matrix}

The analogy between bases for topology and bases for linear algebra in
Section~\ref{lcpct} can also be applied to morphisms.
In this section we identify the abstract conditions satisfied by the relation
$$ f K^n \;\subset\; U^m  \quad\hbox{or}\quad   K^n \;\subset\; f^*U^m  $$
that was used in Definitions \ref{rec cts} and~\ref{rec cts latt}.
Like $\waybelow$,
this is a binary relation between the two overt discrete objects of codes.

As this condition is an observable property of the function,
it determines an open subobject of the set of functions $X\to Y$,
and such properties form a \emph{sub}-basis
of the \textdf{compact--open topology}.
If $X$ is locally compact then this object is the exponential $Y^X$
in both traditional topology and locale theory,
and the conditions below are those listed in
\cite[Lemma VII 4.11]{JohnstonePT:stos}.
Unfortunately, there need not be a corresponding dual basis of compact
subobjects to make $Y^X$ locally compact.

In our notation, the relation $f K^n\subset U^m$ is $A_n(\Sigma^f\beta^m)$.
Clearly this question can easily be generalised to $A_n(H\beta^m)$,
in which we replace $f$ by the ``second class'' map $\hayo H$
(Notation~\ref{Hayo notation}).

\begin{r@Notation}
\label{matrix notn}
Let $(\beta^n,A_n)$ and $(\gamma^m,D_m)$ be $\lor$-bases for objects
$X$ and $Y$ respectively.
Then any first or second class morphism $\hayo H:X\Hto Y$ in $\Hayo\S$
(that is, $H:\Sigma^Y\to\Sigma^X$ in $\S$) has a \textdf{matrix},
$$ \hayo H^m_n \;\equiv\; A_n(H\gamma^m):\Sigma, $$
which generalises the representation of $y:Y$ \via $H\psi\equiv\psi x$
as $\xi\equiv i y\equiv\Lamb y.\gamma^m y$ in Lemmas~\ref{lattice rounded}ff.
\end{r@Notation}

\begin{r@Lemma}
\label{H from Hmn} $H$ is recovered from $\hayo H^m_n$ as 
$H\psi=\Some m n.D_m\psi\land\hayo H^m_n\land\beta^n$.
\end{r@Lemma}

\Proof \closeupaline
\begin{eqnarray*}
  H\psi &=& H(\Some m.D_m\psi\land\gamma^m) &$(\gamma^m,D_m)$ basis for $Y$\\
  &=& \Some m.D_m\psi\land H\gamma^m &indeed $\lor$-basis\\
  &=& \Some m.D_m\psi\land\Some n.A_n(H\gamma^m)\land\beta^n
       &$(\beta^n,A_n)$ basis for $X$\\
  &=& \Some m n.D_m\psi\land \hayo H^m_n\land\beta^n &definition \qEd
\end{eqnarray*}

\closeupaline

\begin{r@Lemma}
\label{Hmn0+} The matrix is directed in $n$, \cf Lemma~\ref{0+ll}:
$$ \hayo H^m_0\eq\top \quad\hbox{and}\quad
   \hayo H^m_{n+p}\eq\hayo H^m_n\land\hayo H^m_p. $$
\end{r@Lemma}

\Proof These are $A_0\phi\eq\top$ and $A_{n+p}\phi\eq A_n\phi\land A_p\phi$
with $\phi=H\gamma^m$. They hold because $(\beta^n,A_n)$ is an $\lor$-basis
(Definition~\ref{basis jargon}(a)). \qed

\begin{r@Lemma}
\label{<Hmn<} The matrix is monotone in $m$, \cf Lemma~\ref{<.<<.<}: 
$$ (n'\baseleq_X n)\land\hayo H^m_n\land(m\baseleq_Y m')
   \;\Implies\; \hayo H^{m'}_{n'}. $$
\end{r@Lemma}

\Proof $A_n\leq A_{n'}$ (though we already had this from directedness)
and $\gamma^m\leq\gamma^{m'}$. \qed

\begin{r@Lemma}
\label{<<Hmn<<} The matrix is rounded, respecting $\waybelow$
on both sides, \cf Lemmas \ref{interpolation} \&~\ref{lattice rounded}:
$$ \Some m'. (m'\waybelow_Y m)\land\hayo H^{m'}_n
   \;\eq\; \hayo H^{m}_{n} \;\eq\;
   \Some n'. \hayo H^m_{n'}\land (n\waybelow_X n').$$
\end{r@Lemma}

\Proof\closeupaline
\begin{eqnarray*}
   \hayo H^{m'}_n &\eq& A_n(H\gamma^{m'}) \\
   &\eq& A_n\big(H(\Some m.D_m\gamma^{m'}\land\gamma^m)\big)
      &$(\gamma^m,D_m)$ basis for $Y$\\
   &\eq& \Some m.D_m\gamma^{m'}\land A_n(H\gamma^m) &indeed, $\lor$-basis\\
   &\eq& \Some m.(m\waybelow m')\land\hayo H^m_n \\
   \hayo H^m_{n'}
   &\eq& A_{n'}(H\gamma^m) \\
   &\eq& A_{n'}\big(\Some n.A_n(H\gamma^m)\land\beta^n\big)
      &$(\beta^n,A_n)$ basis for $X$\\
   &\eq& \Some n.A_n(H\gamma^m)\land A_{n'}\beta^n &indeed, $\lor$-basis\\
   &\eq& \Some n.\hayo H_n^m\land (n'\waybelow n) &$\qeds$%
\end{eqnarray*}

\begin{r@Lemma}
\label{recover Hmn} Suppose $n,m:N\proves\rho(n,m):\Sigma$
satisfies the foregoing properties, \ie
$$
   \rho(0,k)
   \qquad
\begin{prooftree}
     \rho(n,k) \quad \rho(m,k)
     \Justifies
     \rho(n+m,k)
   \end{prooftree}
   \qquad
   \begin{prooftree}
     n\waybelow n' \quad \rho(n',m)
     \Justifies
     \rho(n,m)
   \end{prooftree}
   \qquad
   \begin{prooftree}
     \rho(n,m') \quad m'\waybelow m
     \Justifies
     \rho(n,m)
   \end{prooftree}
$$
and define $H:X\Hto Y$ by
$H\psi\;\equiv\;\Some m n.D_m\psi\land\rho(n,m)\land\beta^n$.
Then $\rho(n,m) \eq A_n(H\gamma^m)$.
   \end{r@Lemma}

\Proof Since $\rho$ respects $\waybelow$ on the right,
\begin{eqnarray*}
  H\gamma^m
  &=& \Some m' n.D_{m'}\gamma^m\land\rho(n,m')\land\beta^n \\
  &=& \Some m' n.(m'\waybelow m)\land\rho(n,m')\land\beta^n \\
  &=& \Some n. \rho(n,m)\land\beta^n.
\end{eqnarray*}
Then, since $\rho$ also respects $\waybelow$ on the left
and $A_n$ preserves the join,
which is directed because $\rho$ respects $0$ and $+$,
\begin{eqnarray*}
 A_n(H\gamma^m)
  &\eq& A_n(\Some n'.\rho(n',m)\land\beta^{n'}) \\
  &\eq& \Some n'.\rho(n',m)\land A_n\beta^{n'} \\
  &\eq& \Some n'.\rho(n',m)\land n\waybelow n' 
  \;\eq\; \rho(n,m). &$\qeds$%
\end{eqnarray*}

\begin{r@Lemma}
 $\id^m_n \;\eq\; A_n(\id\beta^m) \;\eq\; (n\waybelow m) \;\eq\; \hayo\E^m_n$.
\end{r@Lemma}

\Proof The relationship with $\E$ follows from Lemma~\ref{dlecr}.  The
unit laws were given by Lemma~\ref{<<Hmn<<}, \cf the Karoubi completion,
which splits idempotents in any category. \qed

\begin{r@Lemma}
 $\hayo{K\cdot H}^k_n \eq \Some m.\hayo K^k_m\land\hayo H^m_n$.
\end{r@Lemma}

\Proof\closeupaline
\begin{eqnarray*}
  \hayo{K\cdot H}^k_n 
  &\eq& A_n\big(H(K\epsilon^k)\big) \\
  &\eq&
  A_n\big(\Some m n'. D_m(K\epsilon^k)\land \hayo H^m_{n'}\land\beta^{n'}\big)
     &Lemma~\ref{H from Hmn}\\
  &\Impliedby&
      \Some m n'. D_m(K\epsilon^k)\land \hayo H^m_{n'}\land A_n\beta^{n'}\\
  &\eq& \Some m n'. D_m(K\epsilon^k)\land \hayo H^m_{n'}\land (n\waybelow n')\\
  &\eq& \Some m. D_m(K\epsilon^k)\land \hayo H^m_n &Lemma~\ref{<<Hmn<<}\\
  &\eq& \Some m. \hayo K_m^k\land \hayo H^m_n,    &def $K$% 
\end{eqnarray*}
although $\Impliedby$ is actually $\eq$ as we have an $\lor$-basis.
(We draw attention to this because \cite{TaylorP:dedras}
uses the natural $\land$-basis on $\realno$,
and therefore the weaker result about composition.) \qed

\begin{r@Theorem}
 $\Hayo\S$ (Notation~\ref{Hayo notation}) is equivalent to the
category whose
\begin{letterlist}
\item objects are abstract bases $(N,0,1,{+},{\star},{\waybelow})$
   (Definition~\ref{abs basis});
\item morphisms are $\rho(n,m)$ satisfying the conditions in
  Lemma~\ref{recover Hmn};
\item identity is $\waybelow$;
\item composition is relational.
\qed\end{letterlist}
\end{r@Theorem}

\smallskip

The definition of $\Hayo\S$ in \cite{TaylorP:sobsc}
was essentially taken from Hayo Thielecke's work on ``computational effects''
\cite{ThieleckeH:phd},
which was in turn based on the Kleisli category for the monad.
It was therefore motivated by more syntactic considerations than ours.
From a semantic point of view,
it would have been more natural to have split the idempotents.
In the classical models, the category would then be
the opposite of that of
all continuous (but not necessarily distributive) lattices
and Scott-continuous maps (\cf Example~\ref{eg cts latt basis}).
In this result, we would drop the $\waybelow+$ and $\waybelow\star$ rules
from Definition~\ref{abs basis} of an abstract basis.

\bigskip

In order to characterise first class maps, by Corollary~\ref{tfae hom}
we have to consider preservation of the lattice connectives,
\cf Lemmas~\ref{adm=rhom}ff.
For this, the target object $Y$ must have a lattice basis.

\goodbreak

\begin{r@Lemma}
 $H\top=\top$ iff $\hayo H^1_n\eq(n\waybelow 1)$, and
$H\bot=\bot$ iff $\hayo H^0_n\eq(n\waybelow 0)$.
\end{r@Lemma}

\nobreak

\Proof If $H\top=\top$ then
$\hayo H^1_n\equiv A_n(H\gamma^1)\equiv A_n(H\top) \eq A_n\top
\equiv A_n\beta^1\equiv(n\waybelow1)$
by Notation~\ref{matrix notn}, Definition~\ref{basis jargon},
and Notation~\ref{<< notn}.

Conversely, $H\top\equiv H\gamma^1\equiv\Some n.\hayo H^1_n\land\beta^n \eq
\Some n.A_n\top\land\beta^n\equiv \top$ 
by Definition~\ref{basis jargon},
Lemma~\ref{<<Hmn<<} and Definition~\ref{eff basis}.

We may substitute $\bot$ and $0$ for $\top$ and $1$ in the same argument.
\qed

\medbreak

Similarly we are able on this occasion to handle $\land$ and $\lor$
simultaneously (Remark~\ref{+=x}).

\begin{r@Lemma}
\label{Hmn+*} $H(\phi\sqr\psi)=H\phi\sqr H\psi$ iff
$\hayo H^{s\sqr t}_n \eq
 \Some m p.\hayo H^s_m\land\hayo H^t_p\land (n\waybelow m\sqr p)$,
both when $\sqr$ is $\land$ or $\star$ and when it is $\lor$ or $+$.
\end{r@Lemma}

\Proof If $H(\phi\sqr\psi)=H\phi\sqr H\psi$ then
\begin{eqnarray*}
  \hayo H^{s\sqr t}_n &\eq&
    A_n(H\beta^{s\sqr t}) \;\eq\; A_n\big(H(\beta^s\sqr\beta^t)\big) \\
   &\eq& A_n(H\beta^s\sqr H\beta^t)
    &hypothesis\\
    &\eq& \Some m p.A_m(H\beta^s)\land A_p(H\beta^t)\land (n\waybelow m\sqr p)
    &Lemma \ref{ll square}\\
    &\eq& \Some m p.\hayo H_m^s\land \hayo H_p^t\land (n\waybelow m\sqr p) \\
\noalign{Conversely, using distributivity,}
 H\phi\sqr H\psi
    &=& (\Some m u.D_m\phi\land\hayo H^m_u\land\beta^u)
     \sqr(\Some p v.D_p\psi\land\hayo H^p_v\land\beta^v) 
    &Lemma~\ref{H from Hmn}\\
    &=& \Some m p u v.D_m\phi\land D_p\psi\land\hayo H^m_u\land\hayo H^p_v
     \land\beta^{u\sqr v}  &Lemma~\ref{or distrib dirsup}\\
    &=& \Some k m p u v.D_m\phi\land D_p\psi\land
       \hayo H^m_u\land\hayo H^p_v\land (k\waybelow u\sqr v)
     \land\beta^k &L.~\ref{<<b}\\
    &=& \Some k m p.D_m\phi\land D_p\psi\land\hayo H^{m\sqr p}_k
     \land\beta^k &hypothesis\\
    &=& \Some n k m p.D_m\phi\land D_p\psi\land
      (n\waybelow m\sqr p) \land \hayo H^n_k \land\beta^k
      &Lemma~\ref{<<Hmn<<}\\
    &=& \Some n k.D_n(\phi\sqr\psi)\land \hayo H^n_k\land\beta^k
     &Lemma~\ref{ll square}\\
   &=& H(\phi\sqr\psi) &Lemma~\ref{H from Hmn} \qEd
\end{eqnarray*}

\begin{r@Definition}
\label{def abs matrix} An \textdf{abstract matrix} is a
binary relation $\hayo H^m_n$ such that
$$
   \hayo H^m_0
   \qquad
\begin{prooftree}
     \hayo H^m_n\quad\hayo H^m_p
     \Justifies
     \hayo H^m_{n+p}
   \end{prooftree}
   \qquad
   \begin{prooftree}
      n'\baseleq_X n  \quad  \hayo H^m_n  \quad  m\baseleq_Y m'
      \justifies
      \hayo H^{m'}_{n'}
   \end{prooftree}
   \qquad
   \begin{prooftree}
     \hayo H^{m'}_{n'} 
     \Justifies
     n'\waybelow_X n  \quad   \hayo H^m_n  \quad  m\waybelow_Y m'
   \end{prooftree}
$$
$$
\begin{prooftree}
  n\waybelow 0
  \Justifies
  \hayo H^0_n
\end{prooftree}
\qquad
\begin{prooftree}
  n\waybelow 1
  \Justifies
  \hayo H^1_n
\end{prooftree}
\qquad
\begin{prooftree}
  \hayo H^{s\star t}_n 
  \Justifies
  \hayo H^s_m  \quad  \hayo H^t_p  \quad n\waybelow m\star p
\end{prooftree}
\qquad
\begin{prooftree}
  \hayo H^{s + t}_n 
  \Justifies
  \hayo H^s_m  \quad  \hayo H^t_p  \quad n\waybelow m + p
\end{prooftree}
$$
   \end{r@Definition}

\begin{r@Theorem}
 $\S$ is equivalent to the category whose
\begin{letterlist}
\item objects are abstract lattice bases
   $(N,0,1,{+},{\star},{\waybelow})$;
\item morphisms are abstract matrices;
\item identity is $\waybelow$;
\item composition is relational.
\qed\end{letterlist}
\end{r@Theorem}

\medskip

Jung and S\"underhauf characterised continuous functions between
stably locally compact spaces in a similar way \cite{JungA:duacvo}.

\begin{r@Remark}
 We shall speculate on the possible computational applications
of matrices for continuous functions in the next section,
but let's say something here about the analogy
with linear algebra that we have used.
Really, this has been much more useful that we had any right to expect,
since neither $\S$ nor $\Hayo\S$ (the ``first'' and ``second class'' maps)
is a symmetric monoidal closed category.

However, by Theorem~\ref{AD=DA}, stably locally compact objects
and the $\hayo A:X\Hto Y$ for which $A$ preserves $\top$ and $\land$
do define such a category.
In fact, they provide a model of linear logic with involutive negation
and an ``of course'' operator given by the Smyth powerdomain $\powerset^\sharp$
(Example~\ref{upd}):
\end{r@Remark}

\begin{r@Proposition}
\label{upd prop}
Let $X$ be a stably locally compact space.
The $\Hayo\S$-maps $\hayo A:\Gamma\Hto X$ for which $A:\Sigma^X\to\Sigma^\Gamma$
preserves $\top$ and $\land$
correspond bijectively and naturally to $\S$-maps $\Gamma\to\powerset^\sharp X$.
\end{r@Proposition}

\Proof Both correspond to terms $\Gamma\proves\xi\equiv\Lamb n.A\beta^n:\Sigma^N$
that are rounded filters ($\land$-homomorphisms) for $\waybelow_X$
or ideals for $\wayabove_X$. \qed

\smallskip

Since they also preserve directed joins,
such $A$ are known as \textdf{preframe homomorphisms};
for more on their relationship to the Smyth powerdomain, 
see \cite[11.2.5]{VickersSJ:topvl} and \cite{JungA:stacsc}.
For a similar investigation of join-preserving maps and
the Hoare or lower powerdomain $\powerset^\flat X$ in ASD,
presumably we would first need to identify the ``stably locally overt'' objects
to which the analogous construction may be applied.

%============================================================================
\section{Relating the classical and term models}\label{relate}

We promised to translate each step of the logical development into
topological language,
but we haven't done this since the end of Section~\ref{subobjects}.
We shall now show how the ``coding'' provides the link between
the classical and computational models (Remark~\ref{two models}),
and then how abstract matrices themselves describe exact computation
for the reals and locally compact objects in general.

\begin{r@Remark}
 On the one hand, we know from
the classical proof for $\LKSp$ in \cite[Theorem 5.12]{TaylorP:sobsc}
and the intuitionistic one for $\LKLoc$ in \cite[Theorem 3.11]{TaylorP:subasd}
that these categories are models of the calculus in
Sections~\ref{axioms I}--\ref{axioms II},
\ie that there are interpretation functors $\denote-:\S\to\LKSp$
and  $\denote-:\S\to\LKLoc$.
In the light of Theorem~\ref{all loc cpct},
that every object of $\S$ is a $\Sigma$-split subobject of $\Sigma^\natno$,
the converse part of Theorem~\ref{ctsdist<UKN}
provides another proof in the localic setting.

In this paper we have sought the ``inverse'' of this functor.
Since the classical models are of course richer,
they have to be constrained in order to obtain something equivalent
to the computational one (Remark~\ref{two models}).
This constraint was in the form of a \emph{computable basis},
as in Definitions \ref{rec loc cpct} and~\ref{rec cts latt}.
Nevertheless, as we saw in the case of $\realno$ (Example~\ref{eg R1}), 
such bases may already be familiar to us from traditional considerations.

There is no need to verify the consistency conditions that we set out
in Sections \ref{way-below} and~\ref{X from <<}.
They follow automatically from the existence of the classical space,
which serves as a reference as in Remark~\ref{spc as ref}.
So long as $\star$, $+$ and $\waybelow$ are defined by programs,
which can be translated into our $\lambda$-calculus,
we already have an abstract basis.
\end{r@Remark}

\begin{r@Examples}
\label{concl egs}
At this point let us recall the various ways in which
a lattice basis can be defined on a locally compact sober space or locale.
\begin{letterlist}
\item In a stably locally compact sober space (which is, in
  particular, compact in the global sense), we may choose a sublattice
  of compact subspaces $K^n$ and corresponding sublattice of open
  ones $U^n$, such that $U^n\subset K^n$ and the basis expansion is
  satisfied.  Then the indexing set $N$, together with the operations
  $+$ and $\star$ on codes corresponding to unions and intersections
  of open--compact pairs, and the relation $n\waybelow m$ given by
  $K^n\subset U^m$, define an abstract basis, so long as these
  operations are computable.
  In the corresponding lattice filter basis in the $\lambda$-calculus,
  $\beta^n$ and $A_n$ classify $U^n$ and the Scott-open filter
  $\F_n\equiv\collect V{K^n\subset V}$.
\item In a compact Hausdorff space, the compact subspaces $K^n$
  are the complements of open subspaces $V_n$,
  which may be chosen from the same sublattice as the $U^n$,
  but with $U^n\disjoint V_n$.
\item In particular, finite unions of open or of closed rational intervals
  provide this structure for the closed real unit interval $[0,1]$. 
\item $\realno$ is not globally compact, though binary intersections
  of compact subspaces are compact. The lattice basis may be defined
  in the same way as for $[0,1]$, with the single exception of
  $A_1\equiv\Lamb\phi.\bot$, which does not preserve $\top$.
\item Let $(N,{\trless})$ be a recursively enumerable directed
  interpolative relation (Definition~\ref{fir}).
  Then Theorem~\ref{prime basis} defines an object in~$\S$,
  whose classical interpretation is the continuous dcpo of
  rounded ideals of $(N,{\trless})$;
  it is algebraic iff $\trless$ is reflexive.
\item The reflexive order $\baseleq$ defined from any imposed
  distributive lattice $(N,0,1,{+},{\star})$ by Definition~\ref{leq from +*}
  satisfies the conditions on $\waybelow$ for an abstract basis,
  and so defines an object of $\S$ whose classical interpretation
  is the coherent space whose compact open subspaces are indexed by$N$.
\item Given a locally compact locale,
  we choose a sublattice $N$ that provides a basis
  for the corresponding continuous distributive lattice $L$,
  with inclusion $\beta^\blank:N\to L$.
  Then define $(n\waybelow m)\equiv(\beta^n\ll\beta^m)$.
  This is a lattice basis,
  but only a filter basis in the stably locally compact case.
\item Finally, in the case of a non-stably locally compact sober
  space, we only have a $\cup$-semilattice $N$ of compact subspaces $K^n$,
  and therefore a (filter) $\lor$-basis $(\beta^n,A_n)$.
  Remark~\ref{make and basis} turned this into a lattice basis
  $(\beta^\ell,A_\ell)$ indexed by $\Fin(N)$, by defining
  $$ \beta^\ell x \;\equiv\; \All n\in\ell.x \in U^n
     \quad\hbox{and}\quad
     A_\ell\phi \;\equiv\; \Some n\in\ell.K^n\subset V, $$
  where $\phi$ classifies $V$ as usual.
  Then $\beta^\ell$ simply classifies the intersection of the basic open
  subspaces as in the stably locally compact case, 
  but $A_\ell$ is a logical disjunction, not a union of subspaces,
  \cf Lemma~\ref{KcapL}.  Then
  $$ \ell\waybelow \ell' \;\eq\; \Some n\in\ell.\All m\in\ell'.K^n\subset U^m, $$
  but this is not a filter basis.
\end{letterlist}
\end{r@Examples}

\begin{r@Remark}
 Section~\ref{X from <<} showed that the abstract basis is
inter-definable with the nucleus~$\E$. These are interpreted both in
the computational model $\S$ and in the classical ones $\LKSp$ and
$\LKLoc$.  This means that the idempotent $\denote\E$ on
$\Upsilon\Kur\natno$ is the one that defines the $\Sigma$-split
embedding of the original space in $\powerset(\natno)$,
as in Theorems \ref{lcpct<PN} and~\ref{ctsdist<UKN}.

\newarrow{ToFrom}<--->
\newdiagramgrid{123}{.8,.8,1.1,1.1,1.1,1.1}{1,1,1,1}

Hence any classically defined locally compact sober space or locale
that has a computable basis may be ``imported'' into abstract Stone duality
as an object, whose interpretation in the classical model is homeomorphic to
the given space.
Summing this up diagrammatically,
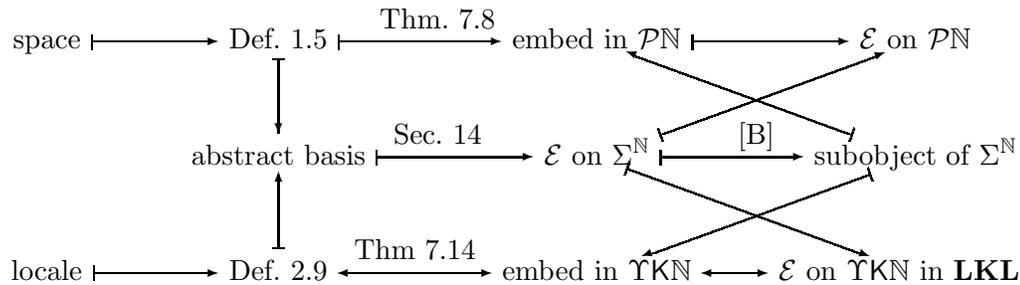
\begin{diagram}[width=5em,height=2em,grid=123]
  \hbox{space} & \rMapsto & \hbox{Def.~\ref{rec loc cpct}}
     & \rMapsto^{\hbox{Thm.~\ref{lcpct<PN}}} &
     \hbox{embed in }\powerset\natno & \rMapsto
      & \E \hbox{ on } \powerset\natno\\
   && \dMapsto &&& \ruMapsto\luMapsto \\
   &&\hbox{abstract basis} & \rMapsto^{\hbox{Sec.~\ref{X from <<}}}
     & \E \hbox{ on } \Sigma^\natno
     & \rMapsto^{\hbox{\cite{TaylorP:subasd}}}
     & \hbox{subobject of } \Sigma^\natno \\
   && \uMapsto &&& \rdMapsto\ldMapsto \\
  \hbox{locale} & \rMapsto & \hbox{Def.~\ref{stab loc cpct loc}}
     & \rToFrom^{\hbox{Thm~\ref{ctsdist<UKN}}} &
     \hbox{embed in }\Upsilon\Kur\natno & \rToFrom
     & \E \hbox{ on } \Upsilon\Kur\natno \hbox{ in } \LKLoc
\end{diagram}
\end{r@Remark}

\smallskip

\begin{r@Remark}
 Having fixed computational bases for two classically defined
spaces or locales, $X$ and~$Y$, we may look at continuous functions
$f:X\to Y$.
By the basis property, any such map is determined by the relation
$$ f K^n\;\subset\; V^m $$
as $n$ and $m$ range over the bases for $X$ and $Y$ respectively.
If this relation is recursively enumerable (Definition~\ref{rec cts})
then the corresponding program may be translated into our $\lambda$-calculus.
Just as we saw for abstract bases,
the resulting term satisfies Definition~\ref{def abs matrix}
for an abstract matrix,
because its interpretation agrees with a continuous function.
Between the objects whose denotations are $X$ and $Y$ there is
therefore a term whose denotation is~$f$.

In particular, computationally equivalent bases for the same space
give rise to an isomorphism between the objects of~$\S$.
In the case of morphisms, extensionally equivalent terms give rise
to the same $(f K^n\subset V^m)$-relations,
and therefore to the same continuous functions.
However, programs may be extensionally equivalent for some deep
mathematical reason, or as a result of the stronger logical principles
in the classical situation, without being provably equivalent
within our calculus.
This is the reason why we required the computable aspects of the
definitions in the Introduction to be accompanied by actual programs.
\end{r@Remark}

\medbreak

\goalbreak{5\baselineskip}

This completes the proof of our main result:

\begin{r@Theorem}
 Abstract Stone duality, \ie the free model $\S$
of the axioms in Sections \ref{axioms I} and~\ref{axioms II},
is equivalent to the category of computably based locally compact locales
and computably continuous functions. \qed
\end{r@Theorem}

\smallskip

\begin{r@Remark}
\label{lacuna}
 An obvious lacuna in this result arises from the difference
between sober spaces and locales: we are relying on the axiom of
choice within the classical models to say that the two are the same.
Recall \cite[Theorem VII 4.3]{JohnstonePT:stos}
that, using excluded middle, the crucial requirement is to find,
for any $\phi\not\leq\psi:L$ in a distributive continuous lattice,
\begin{letterlist}
\item some (Scott open) filter $A:\Sigma^L$ for which $A\phi\eq\top$
  but $A\psi\eq\bot$, and then
\item using Zorn's lemma, enlarge $A$ to a maximal such filter $P$,
  which is prime.
\end{letterlist}
\smallskip

For stage (a), we would like to show, more generally,
that every object (which we know carries a lattice basis)
also has a filter $\lor$-basis.
The idea, due to Jimmie Lawson \cite[\S I~3.3]{GierzGK:comcl},
is to iterate the interpolation property \cite{TaylorP:loccbc}.
This can be shown with the aid of a weaker choice principle that
merely extracts a total function from a non-deterministic one.

For stage (b), since any space has a $\natno$-indexed lattice basis,
we only need to consider $\Sigma^\natno$,
and not \emph{general} distributive continuous lattices.
The filter $A$ provides a rounded filter (\emph{semi}lattice homomorphism)
$\xi_0:\Sigma^\natno$ such that $I\phi\xi\eq\top$ but $I\psi\xi\eq\bot$,
which we enlarge to a rounded \emph{lattice} homomorphism $\xi$
(Lemmas~\ref{lattice rounded}ff) with the same values on $\phi$ and $\psi$.
To do this, we may replace Zorn's Lemma by an argument that uses
excluded middle but not Choice:
we build up $\xi$ from $\xi_0$ by adding each number 0, 1, 2, ... in turn
if it generates a proper filter. \qed
\end{r@Remark}

\medskip

\begin{r@Remark}
 Let's review what we've achieved by way of
a \textdf{type theory} for topology.
\begin{letterlist}
\item The original idea of abstract Stone duality was that the
  non-computable unions could be eliminated from general topology by
  expressing the category of ``frames'' by a monadic adjunction over
  its opposite category of ``spaces'' rather than over sets.
\item Beck's theorem says that this is equivalent to the, perhaps less
  friendly, condition that the functor $\Sigma^\blank$ reflect
  invertibility and ``create $\Sigma$-split coequalisers''.
\item In \cite{TaylorP:subasd} we saw that the latter can be
  interpreted as (certain) subobjects, and that the data for such
  subobjects could be encapsulated in a single morphism $\E$,
  called a ``nucleus''.
\item This was further developed, in Section~8 of that paper, into a
  $\lambda$-calculus similar to comprehension in set theory. 
  However, the data defining a subobject remained arcane:
  terms satisfying the equation defining a nucleus
  could only be found with considerable expert ingenuity.
\item In Theorem~\ref{lattice nucleus} of this paper we reformulated the
  $\lambda$-equation in a much simpler way using the lattice connectives.
\item The notion of abstract basis in this paper puts the construction
  within the grasp of anyone who has a knowledge of open and compact
  subspaces in topology.
\end{letterlist}
\end{r@Remark}

% ---------------------------------------------------------------------
\bigbreak

\begin{r@Remark}
\label{matrix compute} Let us consider
the computational meaning of the matrix $\hayo H^m_n$
when $H\equiv\Sigma^f$ for some continuous function $f:\realno\to\realno$.
In order to have a lattice basis,
the indices $n$ and $m$ must range over (finite) \emph{unions} of intervals,
although by Lemma~\ref{Hmn0+}, $n$ need only denote a single closed interval.
The matrix therefore encodes the predicate
$$ f[x_0\pm\delta] \;\subset\; \bigcup_j(y_j\pm\epsilon_j), $$
as an ASD term of type $\Sigma$, the union being finite.
Suppose that we have a real input value $x$ that we know to lie in the
interval $[x_0\pm\delta]$, and we require $f(x)$ to within $\epsilon$.

We substitute the rational values $x_0$, $\delta$ and
$\epsilon_j\equiv\epsilon$ in the predicate, leaving $(y_j)$ indeterminate.
Recall from Remark~\ref{normalisation} and \cite[Remark~11.3]{TaylorP:sobsc}
that any such term may be translated into a \LPROLOG\ program.
Such a program permits substitution of values for any subset of
the free variables, and is executed by resolving unification problems,
which result in values of (or at least constraints on) the remaining variables.
In this case, we obtain (nondeterministically) some finite set $(y_j)$.

In the language of real analysis, we are seeking to cover the compact
interval $f[x_0\pm\delta]$ with (finitely many) open intervals of size
$\epsilon$, centred on the $y_j$.
The Wilker property (Lemma~\ref{Hmn+*}) then provides
$$ [x_0\pm\delta] \;\subset\; \bigcup_j(x_j\pm\delta')
   \quad\hbox{with}\quad
   f[x_j\pm\delta'] \;\subset\; (y_j\pm\epsilon). $$
Responsibility now passes back to the supplier of the input value $x$
to choose which of the $x_j$ is nearest, and the corresponding $y_j$
is the required approximation to the result $f(x)$.
This discussion is taken up again in~\cite{TaylorP:dedras,TaylorP:lamcra}. \qed
\end{r@Remark}

\begin{r@Remark}
 This illustrates the way in which we would expect to use
abstract Stone duality for computations with objects such as $\realno$
that we regard, from a mathematical point of view, as ``base types''
(though of course only $\terminalobj$, $\natno$ and $\Sigma$ are
actually base types of our $\lambda$-calculus).
Where higher types, such as continuous or differentiable function-spaces,
can be shown to be locally compact, they too have bases and matrices,
but it would be an example of the mis-use of normalisation theorems
(Remark~\ref{nfthm}) to insist on reducing everything to matrix form.
We would expect to use higher-type $\lambda$-terms of our calculus
to encode higher-order features of analysis, for example in the
calculus of variations.
Of course, much preliminary work with $\realno$ itself needs to be
done before we see what can be done in such subjects
\cite{TaylorP:dedras,TaylorP:lamcra}.
\end{r@Remark}

\begin{r@Remark}
 The manipulations that we have done since introducing the
abstract way-below relation have all required \emph{lattice} bases.
The naturally occurring basis on an object such as $\realno$,
on the other hand, is often just an $\land$-basis.
This didn't matter in Sections~\ref{X from <<}--\ref{matrix},
as they were only concerned with the \emph{theoretical} issue
of the consistency of the abstract basis.
We have just seen, however, that the ``matrices'' in Section~\ref{matrix}
encapsulate actual computation, in which the base types are those
of the indices of the bases. The lists used in 
Lemma~\ref{make or basis} would then be a serious burden.

There is a technical issue here that is intrinsic to topology.
In locale theory, which is based on an algebraic theory of finite meets
distributing over arbitrary unions of ``opens'',
it is often necessary to specify when two such expressions are equal,
which may be reduced to the question of when an intersection is
contained in a union of intersections.
This \emph{coverage relation} is an important part of the technology
of locale theory \cite[Section II 2.12]{JohnstonePT:stos},
whilst it was chosen as the focus of the axiomatisation of
Formal Topology.
\end{r@Remark}

\begin{r@Remark}
 The question of whether, using abstract Stone duality,
we can develop a technically more usable approach than these
warrants separate investigation, led by the examples.
Since, in a locally compact object,
we may consider coverages of \emph{compact} subobjects,
the covering families of open subobjects need only be finite.
Jung, Kegelmann and Moshier have exploited this idea to develop a
Gentzen-style sequent calculus \cite{JungA:mullsc}.

As this finiteness comes automatically,
maybe we don't need to force it by using lists.
To put this another way, as the lists act disjunctively,
we represent them by their membership predicates (Remark~\ref{embed in SSN}).
That is, we replace the matrix $\hayo H^n_m$ with the predicate
$$ m:M,\; \xi:\Sigma^N \ \proves\
   A_m\cdot H(\Lamb y.\Some n.\xi n\land\beta^n y):\Sigma,
$$
where $N$ need no longer have~$+$.
In Remark~\ref{matrix compute} above,
we could add the constraint that
$\xi$ only admits intervals of size $\lt\epsilon$.
\end{r@Remark}

\begin{r@Remark}
\label{concl rk}
Another striking feature of the matrix description is that it
reduces the topological theory to an entirely discrete one.
The latter may be expressed in an \emph{arithmetic universe},
which is a category with finite limits, stable disjoint coproducts,
stable effective quotients of equivalence relations and a $\List$ functor
\cite{TaylorP:insema}.

Once again, we need to see this normalisation theorem in reverse.
It appears that any arithmetic universe 
may conversely be embedded as the full subcategory of overt discrete objects
in a model of abstract Stone duality.
Example~\ref{eg arith univ basis} shows that we can do this by modelling
a much simpler structure than that of abstract bases and matrices
for general locally compact spaces.
The differences in logical strength amongst models of ASD
are therefore measured by the overt discrete objects that they contain.

Such a construction would enable
topological, domain-theoretic and $\lambda$-calculus reasoning
to be applied to problems in discrete algebra and logic.
Topologically, it would strongly vindicate Marshall Stone's dictum,
\emph{always topologise},
whilst computationally it would provide continuation-passing translations
of discrete problems,
and of type theories for inductive types.
\end{r@Remark}

%============================================================================
%\iflmcs\newpage\else\allowlines{5}\fi

{\small\raggedright
\bibliographystyle{alpha}\def\longestlabel{ABC+99}%
\bibliography{baire}

\medskip

%\iflmcs\newpage\else\allowlines{5}\fi

\input{asdbib.sty}
}% end of {\small\raggedright

\ifx\href\undefined\def\href#1#2{#2}\fi

\parindent 1em

\section*{Acknowledgements}

{\small
This paper evolved from \cite{TaylorP:loccbc},
which included roughly Sections
\ref{cpct subsp}--%, \ref{lcpct}, \ref{subobjects},
\ref{bases} and \ref{way-below} of the present version.
It was presented at
\href{http://www.mathstat.uottawa.ca/lfc/ctcs2002}%
{\emph{Category Theory and Computer Science}~\textbf{9}},
in Ottawa on 17 August 2002,
and at \emph{Domains Workshop}~\textbf{6} in Birmingham a month later.
The characterisation of locally compact objects using effective bases
(Section~\ref{lcpct}) had been announced on \texttt{categories}
on 30 January 2002.

The earlier version also showed that any object has a filter basis,
and went on to prove Baire's category theorem, that the intersection
of any sequence of dense open subobjects (of any locally compact overt object)
is dense.
These arguments were adapted from the corresponding ones in the theory
of continuous lattices \cite[Sections I~3.3 and~3.43]{GierzGK:comcl}.

I would like to thank
Andrej Bauer,
Mart{\iacute}n Escard\'o,
Peter Johnstone,
Achim Jung,
Jimmie Lawson,
Graham White
and the CTCS and LMCS referees for their comments.
Graham White has given continuing encouragement
throughout the abstract Stone duality project,
besides being an inexhaustible source of mathematical ideas.

This research is now supported by UK EPSRC project GR/S58522,
but this funding was obtained in part on the basis of the work in this paper.
Apart from
\ref{classical dual Wilker},
\ref{upd}, \ref{ASD dual Wilker},
\ref{adm=rhom}, \ref{upd rk}, \ref{upd prop} and~\ref{lacuna},
the work here was carried out during a period of unemployment,
supported entirely from my own savings.
However, I would have been unable to do this without
the companionship and emotional support of my partner, Richard Symes.

%\copyright\ 2002--5 Paul Taylor.

}% end of {\small

%============================================================================
\end{document}